\documentclass[a4paper,12pt]{article}
\usepackage{bm}

\usepackage{hyperref}
\usepackage[english]{babel}
\usepackage{tikz}
\usetikzlibrary{intersections,calc,arrows.meta}
\usepackage{graphicx}
\makeatletter 

\@addtoreset{figure}{section}
\@addtoreset{table}{section}

\makeatother 

\usepackage{epstopdf}
\usepackage {svg}
\usepackage{here}

\usepackage{amsmath}
\usepackage{amssymb}
\usepackage{amsthm}
\usepackage{mathrsfs}
\usepackage[all]{xy}
\usepackage{float}

\newtheorem{theorem}{Theorem}[section]
\newtheorem{proposition}[theorem]{Proposition}
\newtheorem{lemma}[theorem]{Lemma}
\newtheorem{definition}[theorem]{Definition}

\newtheorem{corollary}[theorem]{Corollary}

\newtheorem{remark}[theorem]{Remark}

\newcommand{\C}{{\mathbb C}}

\newcommand{\SW}{{\mathcal {SW}}}

\numberwithin{equation}{section}

\usepackage{braket}

\newcommand\restr[2]{{
  \left.\kern-\nulldelimiterspace 
  #1 
  \vphantom{\big|} 
  \right|_{#2} 
  }}

\usepackage{comment}

\usepackage{lipsum} 
\title{Tensor structure on the module category of the triplet superalgebra $\mathcal{{SW}}(m)$}
\author{Hiromu Nakano}
\date{}

\begin{document}

\maketitle
\begin{abstract}
We discuss the tensor structure on the category of modules of the $N=1$ triplet vertex operator superalgebra $\mathcal{SW}(m)$ introduced by Adamovi\'{c} and Milas. Based on the theory of vertex tensor supercategories, we determine the structure of fusion products between the simple and projective $\mathcal{SW}(m)$-modules and show that the tensor supercategory on $\mathcal{SW}(m)$-mod is rigid. Technically, explicit solutions of a fourth-order Fuchsian differential equation are important to show the rigidity of $\mathcal{SW}(m)$-modules. We construct solutions of this Fuchsian differential equation using the theory of the Dotsenko-Fateev integrals developed by Sussman.
\end{abstract}

\section{Introduction}
In recent years, comprehensive studies of logarithmic vertex operator algebras have developed, where the adjective ``logarithmic'' comes from
the property that \(L_0\) the generator of scale transformations is no longer diagonalizable.
In the logarithmic conformal field theories, this non-diagonalizability leads to interesting examples of physical systems such as polymers, spin chains, percolations and sand-pile models \cite{CR,GK,Gu,JPR,MR,Ni,PR,R,RS}.
In the representation theories of logarithmic vertex operator algebras, indecomposable modules appear such that $L_0$ acts non-semisimply, and this non-semisimplicity made it challenging to formulate the category theories associated with vertex operator algebras.

Among the logarithmic vertex operator algebras, the triplet $W$-algebra $\mathcal{W}_{p}$ \cite{FGST1,FF2,FF3,Ka} is particularly famous and is known for satisfying a cetain finiteness called $C_2$-cofinite condition. 
By the rigorous studies \cite{AM,McRae,NT,TWFusion}, the structure of the category of $\mathcal{W}_{p}$-modules is completely determined. 
Furthermore, by \cite{CLR,GN}, it is shown that the category of $\mathcal{W}_{p}$-modules is equivalent to a quasi-Hopf modification of the category of finite-dimensional representations for the restricted quantum group $\overline{U}_q(sl_2)$ at $q=e^{\frac{\pi i}{p}}$.

As a natural super analogue of the triplet algebra $\mathcal{W}_{p}$, a $\frac{1}{2}\mathbb{Z}_{\geq 0}$-graded (logarithmic) vertex operaor superalgebra, usually denoted by $\mathcal{SW}(m)$, was introduced by Adamovi\'{c} and Milas in the paper \cite{AM2}. This vertex operaor superalgebra $\mathcal{SW}(m)$ is called $N=1$ triplet vertex operaor superalgebra, and contains conformal and superconformal vectors 
at the central charge
\begin{align*}
c_{1,2m+1}=\frac{15}{2}-3(2m+1+\frac{1}{2m+1}),\qquad m\in \mathbb{Z}_{\geq 1}.
\end{align*}
Adamovi\'{c} and Milas proved the $C_2$-cofiniteness of $\mathcal{SW}(m)$, classified all simple $\mathcal{SW}(m)$-modules and conjectured the equivalence between the module category of $\mathcal{SW}(m)$ and the category of finite dimensional modules over the small quantum group $U^{small}_{q}(sl_2)$ at $q=e^{\frac{2\pi i}{2m+1}}$. Furthermore they showed that the characters of the simple $\mathcal{SW}(m)$-modules can be expressed in terms of the characters of the simple $\mathcal{W}_{p}$-modules.

For a $C_2$-cofinite, $\frac{1}{2}\mathbb{Z}_{\geq 0}$-graded vertex operator superalgebra $V$, let $V$-mod be the category of grading-restricted generalized modules. It is shown by Creutzig, Genra, Nakatsuka and Sato \cite{CGNS} that the category $V$-mod admits the structure of the vertex tensor supercategory developed by Huang-Lepowsky-Zhang 
\cite{HLZ1}-\cite{HLZ8} 
and Creutzig-Kanade-McRae \cite{CKM}. 
Recently, as a deeper result, Creutzig, McRae, Orosz Hunziker and Yang \cite{CMOY} have shown that the category of $C_1$-cofinite grading-restricted generalized modules for a vertex operator superalgebra has the structure of the above vertex tensor supercategory when the module category satisfy some appropriate conditions.
By these important results, in particular, $\mathcal{SW}(m)$-mod admits the structure of the vertex tensor supercategory.

When studying the structure of $\mathcal{SW}(m)$-mod, in addition to general theories of vertex operator superalgebras, the rigidity of modules are also important. In the case of $\mathcal{W}_{p}$ \cite{CMY,McRae,TWFusion}, a second-order Fuchsian differential equation, called BPZ differential equation \cite{BPZ}, becomes important for examining the rigidity of modules. 
In the case of $\mathcal{SW}(m)$, the investigation of the rigidity of simple modules leads to a fourth-order Fuchsian differential equation.
This Fuchsian differential equation is derived from the four point function with the minimal conformal weight vectors of $N=1$ super Virasoro simple modules $L(c_{1,2m+1},h_{2,2})$ inserted, where
$h_{2,2}$ is the minimal conformal weight defined by
\begin{equation*}
h_{2,2}=\frac{3}{8}(2m+1)-\frac{3}{4}+\frac{3}{8(2m+1)}.
\end{equation*}
In this paper, we construct the fundamental system of solutions of this fourth-order Fuchsian differential equation using the theory of Dotsenko-Fateev integrals \cite{DF1,DF2} developed by Sussman \cite{S,S2}, and determine some monodromy properties of the solutions.
Those monodromy data are important to show the rigidity of a simple $\mathcal{SW}(m)$-module whose minimal conformal weight is $h_{2,2}$ (Theorem \ref{rigid12}).

The following are the main results and goals of this paper:
\begin{itemize}
\item We determine the structure of the projective covers of the simple $\mathcal{SW}(m)$-modules (Propositions \ref{Proj0}-\ref{projstr}).
\item We show that the tensor supercategory on $\mathcal{SW}(m)$-mod is rigid.
\item We determine the structure of fusion products between the simple and projective $\mathcal{SW}(m)$-modules (Subsection \ref{fusion ring}).
\end{itemize}
We show the above results in Subsections \ref{SecFusion}-\ref{fusion ring}, using properties of vertex tensor supercategories and the Dotsenko-Fateev integrals.
These results are partially based on our thesis \cite{Nak}.

Recently, in \cite{CMOY}, Creutzig, McRae, Orosz Hunziker and Yang obtain important results for the rigidity and the structure of the $C_1$-cofinite module category of the $N=1$ super Virasoro vertex operator superalgebra. 
They show the rigidity and self-duality of the simple module $L(c_{1,2m+1},h_{2,2})$ using a certain embedding technique for Virasoro vertex operator algebras. 
Just as the structure of $\mathcal{W}_{p}$-mod can be determined from the Virasoro tensor category at central charge $13-6p-6p^{-1}$ $(p\in \mathbb{Z}_{\geq 2})$ \cite[Section 7]{McRae},
we expect that their results will rederive the rigidity and structure of $\mathcal{SW}(m)$-mod and further develop the theory of the super singlet algebra \(\mathcal{SM}(m)\) introduced in \cite{AM2}.

\vspace{3mm}
This paper is organized as follows.

In Section \ref{NSalgebra}, we review some facts of the representation theory of $\frac{1}{2}\mathbb{Z}_{\geq 0}$-graded vertex operator superalgebras and $N=1$ Neveu-Schwarz algebra in accordance with \cite{W,CKM,IK1,IK2,KW}.

In Section \ref{NStriplet}, we review some important properties for $\mathcal{SW}(m)\mathchar`-{\rm mod}$, such as the structure of the simple $\mathcal{SW}(m)$-modules and the Zhu-algebra $A(\mathcal{SW}(m))$ determined by Admovi\'c and Milas in 
\cite{AM2,AM3}.

In Section \ref{DF}, 
we introduce a $\mathcal{SW}(m)$ correlation function and examine some monodromy properties of this correlation function. In Subsection \ref{FourthO}, we show that this correlation function satisfies a complicated fourth-order Fuchsian differential equation.
In Subsection \ref{DotsenkoFateev}, we introduce the Dotsenko-Fateev integrals \cite{DF1,DF2} and regularization results of these integrals recently established by Sussman \cite{S2}.
In Subsection \ref{DotsenkoFateev2}, we introduce some formulas for the Dotsenko-Fateev integrals given by Forrester \cite{Forrester0,Forrester} and some results for meromorphic continuations of the Dotsenko-Fateev integrals given by Sussman \cite{S}. 
We also introduce some transformation formulas between Dotsenko-Fateev integrals of special types.
In Subsection \ref{ConSol}, we construct the fundamental system of solutions of the Fuchsian differential equation given in Subsection \ref{FourthO}, using the results in Subsection \ref{DotsenkoFateev2} and a certain deformation technique by Tsuchiya and Wood \cite{TW}.
As a result, we get some monodromy properties of this Fuchsian differential equation.

In Section \ref{NSfusion}, we examine the structure of fusion products and the rigidity of $\mathcal{SW}(m)$-modules.
In Subsection \ref{X''dual}, using the results in Section \ref{DF} and certain vertex tensor categorical techniques introduced by Creutzig-McRae-Yang \cite{CMY,CMY2} and Tsuchiya-Wood \cite{TWFusion}, we show that a simple $\mathcal{SW}(m)$-module $X_2$ is rigid and self-dual.
By using self-duality of the simple module $X_2$, we show that all simple and projective modules can be obtained by repeatedly multiplying $X_2$. As a result, we can determine the structure of all projective modules and show that the tensor supercategory on $\mathcal{SW}(m)\mathchar`-{\rm mod}$ is rigid.
In Subsection \ref{fusion ring}, we introduce a non-semisimple fusion ring ${P}(\mathcal{SW}(m))$ and determine the ring structure. 
\section{Basic definitions and notation}
\label{NSalgebra}
The $N=1$ Neveu-Schwarz algebra is the Lie superalgebra
\begin{equation*}
\mathfrak{ns}=\bigoplus_{n\in\mathbb{Z}}\mathbb{C} L_n\oplus \bigoplus_{r\in\frac{1}{2}+\mathbb{Z}}\mathbb{C} G_r\oplus \bigoplus\mathbb{C} C
\end{equation*}
with the relations $(k,l\in\mathbb{Z},\ r,s\in\mathbb{Z}+\frac{1}{2})$:
\begin{align*}
&[L_k,L_l]=(k-l)L_{k+l}+\delta_{k+l,0}\frac{k^3-k}{12}C,\\
&[L_k,G_r]=(\frac{1}{2}k-r)G_{k+r},\\
&\{G_r,G_s\}=2L_{r+s}+\frac{1}{3}(r^2-\frac{1}{4})\delta_{r+s,0}C,\\
&[L_k,C]=0,\ \ [G_r,C]=0,
\end{align*}
where $\{,\}$ is the anti-commutator.
We identify $C$ with a scalar multiple of the identity, $C=c\cdot {\rm id}$, when acting on modules and refer to the number $c\in\mathbb{C}$ as the central charge. In this section, we briefly review basic facts of representation theories of the Neveu-Schwarz algebra.

\subsection{Vertex operator superalgebras}
In this subsection we briefly review the definitions of $N=1$ vertex operator superalgebras and the notion such as vertex algebra modules and intertwining operators used in later section. See \cite{CKL,CKM,CMOY,FLM,KW} for details. 

Let us recall the definition of $\frac{1}{2}\mathbb{Z}_{\geq 0}$-graded vertex operator superalgebras.
\begin{definition}
 A four pairs $(V,{\mid}{0}\rangle,T,Y)$ is called a $\frac{1}{2}\mathbb{Z}_{\geq 0}$-{\rm graded vertex operator superalgebra} where
\begin{enumerate}
\item $V$ is a $\frac{1}{2}\mathbb{Z}_{\geq 0}$-graded $\mathbb{C}$-vector space
\begin{align*}
&V=\bigoplus_{n\in \frac{1}{2}\mathbb{Z}_{\geq 0}} V[{n}].
\end{align*}
For $\bar{0},\bar{1}\in \mathbb{Z}/2\mathbb{Z}$, set
\begin{align*}
&V^{\bar{0}}:=\bigoplus_{n\in\mathbb{Z}_{\geq 0}} V[{n}],
&V^{\bar{1}}:=\bigoplus_{n\in\mathbb{Z}_{\geq 0}+\frac{1}{2}} V[{n}].
\end{align*}
\item ${\mid}{0}\rangle\in V[{0}]$ is called {\rm the vacuum vector}.
\item $T\in V[{2}]$ is called {\rm the conformal vector}.
\item $Y$ is a $\mathbb{C}$-linear map
\begin{align*}
Y: V\rightarrow {\rm End}_\mathbb{C}(V)[[z,z^{-1}]].
\end{align*}
\end{enumerate}
These data are subject to the following axioms:
\begin{itemize}
\item ${\rm dim}_{\mathbb{C}}V[{0}]=1$ and $0<{\rm dim}_{\mathbb{C}}V[{n}]<\infty$ for any $n\in\frac{1}{2}\mathbb{Z}_{\geq 0}$.
\item For each $v\in V^{\bar{i}}[{h}]$ there exists a {\rm field}
\begin{align*}
Y(v,z)=\sum_{n\in\mathbb{Z}+\frac{i}{2}}v_{n}z^{-n-h}.
\end{align*}
\item 
$Y({\mid}{0}\rangle,z)={\rm id}_V$ and 
\begin{align*}
&Y(v,z){\mid}{0}\rangle-v\in V[[z]]z
\end{align*}
for all $v\in V$.
\item 
The modes of the field
$
Y(T,z)=T^{}(z)=\sum_{n\in\mathbb{Z}}L^{}_nz^{-n-2}
$
define the commutation relations of the Virasoro algebra with fixed central charge $c=c_V$:
\begin{align*}
&[L^{}_k,L_l]=(k-l)L_{k+l}+\delta_{k+l,0}\frac{k^3-k}{12}c_{V}.
\end{align*}
\item The zero mode $L_{0}$ of $T(z)$ acts semisimply on $V$ and 
\begin{align*}
V[{h}]=\{\ v\in V\ {\mid}\ L_0v=hv\ \}.
\end{align*}
\item For all $v\in V$
\begin{align*}
Y(L_{-1}v,z)=\frac{{\rm d}}{{\rm d}z}Y(v,z).
\end{align*}
\item For $v_1\in V^{\bar{i}}$ and $v_2\in V^{\bar{j}}$, the following {\rm super Jacobi identity} holds
\begin{align*}
z^{-1}_0\delta\Bigl(\frac{z_1-z_2}{z_0}\Bigr)&Y(v_1,z_1)Y(v_2,z_2)-(-1)^{{i}{j}}z^{-1}_0\delta\Bigl(\frac{z_2-z_1}{-z_0}\Bigr)Y(v_2,z_2)Y(v_1,z_1)\\
&=z^{-1}_2\delta\Bigl(\frac{z_1-z_0}{z_2}\Bigr)Y(Y(v_1,z_0)v_2,z_2).
\end{align*}
\end{itemize}
\end{definition}
In the above defintion, we call $V^{\bar{0}}$ the $even$ $part$ of $V$ and $V^{\bar{1}}$ the $odd$ $part$ of $V$. For any $v\in V^{\bar{i}}$($i=1,0$), we call $v$ $parity$-$homogeneous$ $vector$ in $V$, and we denote by $|v|=i$ the parity of $v$.

The $N=1$ Neveu-Schwarz vertex operator superalgebras are special cases of the $\frac{1}{2}\mathbb{Z}_{\geq 0}$-graded vertex operator superalgebras, which are subject to an additional axiom:\\
There exists $G\in V[{\frac{3}{2}}]$ (super conformal vector) such that the modes of fields
\begin{align*}
&Y(T,z)=T^{}(z)=\sum_{n\in\mathbb{Z}}L^{}_nz^{-n-2},
&Y(G,z)=G^{}(z)=\sum_{r\in\mathbb{Z}+\frac{1}{2}}G^{}_rz^{-r-\frac{3}{2}},
\end{align*}
define the commutation relations of the Neveu-Schwarz algebra with fixed central charge $c=c_V$:
\begin{equation}
\label{NSrel}
\begin{split}
&[L^{}_k,L_l]=(k-l)L_{k+l}+\delta_{k+l,0}\frac{k^3-k}{12}c_{V},\\
&[L_k,G_r]=(\frac{1}{2}k-r)G_{m+r},\\
&\{G_r,G_s\}=2L_{r+s}+\frac{1}{3}(r^2-\frac{1}{4})\delta_{r+s,0}c_{V}.
\end{split}
\end{equation}
Let us recall the definition of modules of $\frac{1}{2}\mathbb{Z}_{\geq 0}$-graded vertex operator superalgebras.
\begin{definition}
\label{weak-module}
Given a $\frac{1}{2}\mathbb{Z}_{\geq 0}$-graded vertex operator superalgebra $(V,{\mid}{0}\rangle,T,G,Y)$, a {\rm grading restricted generalised} $V$-{\rm module} is a pair $(M,Y_M)$ of a vector space $M$ and a linear map $Y_M$ from $V$ to ${\rm End}M[[z,z^{-1}]]$ satisfying the following conditions
\begin{enumerate}
\item $Y_M({\mid}{0}\rangle,z)={\rm Id}_M$ and the modes of 
\begin{align*}
&Y_M(T,z)=\sum_{n\in\mathbb{Z}}L_n^Mz^{-n-2}
\end{align*}
satisfy the commutation relations of the Virasoro algebra with the central charge $c_V$.
\item For all $v\in V$, 
\begin{align*}
Y_M(L_{-1}v,z)=\frac{{\rm d}}{{\rm d}z}Y_M(v,z).
\end{align*}
\item For $v_1\in V^{\bar{i}}$ and $v_2\in V^{\bar{j}}$, the following super Jacobi identity holds
\begin{align*}
&z^{-1}_0\delta\Bigl(\frac{z_1-z_2}{z_0}\Bigr)Y_M(v_1,z_1)Y_M(v_2,z_2)\\
&\ \ \ \ -(-1)^{ij}z^{-1}_0\delta\Bigl(\frac{z_2-z_1}{-z_0}\Bigr)Y_M(v_2,z_2)Y_M(v_1,z_1)\\
&\ \ \ =z^{-1}_2\delta\Bigl(\frac{z_1-z_0}{z_2}\Bigr)Y_M(Y(v_1,z_0)v_2,z_2),
\end{align*}
where $\delta(z)$ is the formal delta function $\delta(z)=\sum_{n\in \mathbb{Z}}z^n$.
\item 
$M$ is a $\mathbb{C}$-graded superspace
\begin{align*}
M=\bigoplus_{\bar{i}\in\mathbb{Z}/2\mathbb{Z}} M^{\bar{i}}=\bigoplus_{h\in H(M)}M[{h}]
\end{align*}
such that
\begin{itemize}
\item For some finite subset $H_0(M)$ of $\mathbb{C}$, $H(M)=H_0(M)+\frac{1}{2}\mathbb{Z}_{\geq 0}$.
\item For $h\in H(M)$, $M[{h}]=\{\psi\in M: \exists n\geq 0\ {\rm s.t.}\ (L_0-h)^n\psi=0 \}$.
\item $0<{\rm dim}_{\mathbb{C}}M[{h}]<\infty$.
\item For all $v\in V$, $v_{h}M[{h'}]\subset M[h+h']$.
\item For $\bar{i}=\bar{0},\bar{1}$, $M^{\bar{i}}=\bigoplus_{h\in H(M)}M^{\bar{i}}[h]$, where $M^{\bar{i}}[{h}]=M^{\bar{i}}\cap M[{h}]$.
\item For $v\in V^{\bar{i}}$ and $\psi\in M^{\bar{j}}$ $(\bar{i},\bar{j}\in\mathbb{Z}/2\mathbb{Z})$, $Y_M(v;z)\psi\in M^{\bar{i}+\bar{j}}[[z,z^{-1}]]$.
\end{itemize}
\end{enumerate}
\end{definition}
In the above defintion, we call $M^{\bar{0}}$ the $even$ $part$ of $M$ and $M^{\bar{1}}$ the $odd$ $part$ of $M$. For any $\psi\in M^{\bar{i}}$($i=1,0$), we call $\psi$ $parity$-$homogeneous$ $vector$ in $M$, and we denote by $|\psi|=i$ the parity of $\psi$.

Given an $N=1$ Neveu-Schwarz vertex operator superalgebra $V$, a grading restricted generalised $V$-module $M$ is a special case of Definition \ref{weak-module}, which is subject to an additional axiom:
The modes of 
\begin{align*}
&Y_M(T,z)=\sum_{n\in\mathbb{Z}}L_n^Mz^{-n-2}
&Y_M(G,z)=\sum_{r\in\mathbb{Z}+\frac{1}{2}}G^{M}_rz^{-r-\frac{3}{2}}
\end{align*}
satisfy the commutation relations of the Neveu-Schwarz algebra with the central charge $c_V$.

Analogous to non-super cases, contragredient modules and intertwining operators can be defined as follows.
\begin{definition}
Let $V$ be a $\frac{1}{2}\mathbb{Z}_{\geq 0}$-graded vertex operator superalgebra and $M$ be a a grading restricted generalised $V$-module. Let
\begin{align*}
M^*=\bigoplus_{h\in H(M)}M^*[h]
\end{align*}
be the graded dual space of $M$, where $M^*[h]={\rm Hom}_{\mathbb{C}}(M[h],\mathbb{C})$ with parity decomposition 
\begin{align*}
&(M^*)^{\bar{i}}=\bigoplus_{h\in H(M)}(M^{*}[h])^{\bar{i}},
&(M^{*}[h])^{\bar{i}}={\rm Hom}_\mathbb{C}(M^{\bar{i}}[h],\mathbb{C}).
\end{align*}
Let $\langle\ ,\ \rangle$ be the natural dual pairing between $M^*$ and $M$. Then we define the $V$-module structure $Y_{M^*}$ as follows
\begin{align}
\label{cont-act}
\langle Y_{M^*}(v,z)\psi^*,\psi\rangle:=(-1)^{ij}\langle \psi^*,Y_M(e^{zL_{1}}(-z^{-2})^{L_0}v,z^{-1})\psi\rangle,
\end{align}
where $\psi^*\in (M^*)^{\bar{i}}$, $\psi\in M$ and $v\in V^{\bar{j}}$, for $\bar{i},\bar{j}\in \mathbb{Z}/2\mathbb{Z}$.
\end{definition}

\begin{definition}
Let $V$ be a a $\frac{1}{2}\mathbb{Z}_{\geq 0}$-graded vertex operator superalgebra and $M_1$, $M_2$ and $M_3$ a triple of $V$-module. 
Denote by $M_3\{z\}[{\rm log}z]$ the space of formal power series in $z$ and ${\rm log}z$ with coefficient in $M_3$, where the exponents of $z$ can be arbitrary complex numbers and with only finitely many ${\rm log}z$ terms. 
A {\rm parity}-{\rm homogeneous} {\rm intertwining operator} $\mathcal{Y}(\cdot,z)$ of type {\scriptsize{$\begin{pmatrix}
   M_3  \\
   M_1\ M_2
\end{pmatrix}
$}} is a parity-homogeneous linear map
\begin{align*}
\mathcal{Y}: &M_1\rightarrow {\rm End}(M_2,M_3)\{z\}[{\rm log}z],\\
&\psi_1\mapsto \mathcal{Y}(\psi_1,z)=\sum_{t\in\mathbb{C}}\sum_{s\geq 0}(\psi_1)^{\mathcal{Y}}_{t,s}z^{-t-1}({\rm log}z)^s
\end{align*}
satisfying the following conditions for parity-homogeneous vectors $\psi_1\in M_1$, $\psi_2\in M_2$ and $v\in V$:
\begin{enumerate}
\item $\mathcal{Y}(L_{-1}\psi_1,z)=\frac{{\rm d}}{{\rm d}z}\mathcal{Y}(\psi_1,z)$.
\item $(\psi_1)^{\mathcal{Y}}_{t,s}\psi_2=0$ for ${\rm Re(t)}$ sufficiently large.
\item The following super Jacobi identity holds
\begin{equation}
\label{JacobiInt}
\begin{split}
&(-1)^{|v||\mathcal{Y}|}z^{-1}_0\delta\Bigl(\frac{z_1-z_2}{z_0}\Bigr){Y}_{M_3}(v,z_1)\mathcal{Y}(\psi_1,z_2)\\
&-(-1)^{|v||\psi_1|}z^{-1}_0\delta\Bigl(\frac{z_2-z_1}{-z_0}\Bigr)\mathcal{Y}(\psi_1,z_2){Y}_{M_2}(v,z_1)\\
&=z^{-1}_2\delta\Bigl(\frac{z_1-z_0}{z_2}\Bigr)\mathcal{Y}(Y_{M_1}(v,z_0)\psi_1,z_2).
\end{split}
\end{equation}
\end{enumerate}
A general intertwining operator of type {\scriptsize{$\begin{pmatrix}
   M_3  \\
   M_1\ M_2
\end{pmatrix}
$}} is a sum of parity-homogeneous ones.
\end{definition}

Given a $\frac{1}{2}\mathbb{Z}_{\geq 0}$-graded vertex operator superalgebra $V$ and $V$-modules $M_1,M_2,M_3$, we denote by
$I_V${\footnotesize{$
\begin{pmatrix}
\ M_3 \\
M_2\ \ M_1
\end{pmatrix}$}}
the vector superspace of the intertwining operators of type
{\footnotesize{$\begin{pmatrix}
\ M_3 \\
M_2\ \ M_1
\end{pmatrix}$}}. This vector superspace has the parity decomposition
\begin{equation*}
I_V
\begin{pmatrix}
\ M_3 \\
M_2\ \ M_1
\end{pmatrix}
=I^0_V
\begin{pmatrix}
\ M_3 \\
M_2\ \ M_1
\end{pmatrix}
\oplus
I^1_V
\begin{pmatrix}
\ M_3 \\
M_2\ \ M_1
\end{pmatrix}
,
\end{equation*}
where $I^0_V$ and $I^1_V$ are the parity-homogeneous subspaces of $I_V$ whose parity are even and odd, respectively.  

\subsection{Free field realization of the Neveu-Scwarz algebra}
\label{FreeNS}
In this subsection, we review the free field realization of the Neveu-Scwarz algebra in accordance with the papers \cite{W,IK1,IK2}.

First, let us introduce the bosonic Fock modules and the bosonic vertex operators.
The Heisenberg algebra is the Lie algebra
\begin{equation*}
\mathfrak{h}=\bigoplus_{n\in\mathbb{Z}}\mathbb{C} a_n\oplus\mathbb{C}\bold{1} ,
\end{equation*}
with commutation relations:
\begin{align*}
[a_k,a_l]=k\delta_{k+l,0}\bold{1},\qquad[a_k,\bold{1}]=0\ \ \ (k,l\in\mathbb{Z}).
\end{align*}
The Heisenberg algebra $\mathfrak{h}$ has the triangular decomposition
\begin{align*}
&\mathfrak{h}^{\pm}=\bigoplus_{n>0}\mathbb{C} a_{\pm},& \mathfrak{h}^0=\mathbb{C} a_0\oplus \mathbb{C}\bold{1}.
\end{align*}
For $\beta\in\mathbb{C}$, let $\mathbb{C}{\mid}\beta;{\rm B}\rangle$ be the one dimensional representation of $\mathfrak{h}^{\geq}=\mathfrak{h}^+\oplus\mathfrak{h}^0$, which satisfies
\begin{align*}
a_0{\mid}\beta;{\rm B}\rangle=\beta{\mid}\beta;{\rm B}\rangle,\ \ \ \bold{1}{\mid}\beta;{\rm B}\rangle={\mid}\beta;{\rm B}\rangle,\ \ \ \mathfrak{h}^+{\mid}\beta;{\rm B}\rangle=0.\\
\end{align*}
\begin{definition}
The {\rm bosonic Fock module} is defined by induced representation
\begin{align*}
{F}^{\rm B}_{\beta}={\rm Ind}^{\mathfrak{h}}_{\mathfrak{h}^\geq}\mathbb{C}{\mid}\beta;{\rm B}\rangle.
\end{align*}
\end{definition}

Let $a(z)=\sum_{n\in\mathbb{Z}}a_nz^{-n-1}$ and we define the following bosonic energy-momentum tensor 
\begin{align*}
T^{({\rm B};\alpha)}(z)=\frac{1}{2}(:a(z)^2:+\alpha\partial a(z))=\sum_{n\in\mathbb{Z}}L^{({\rm B};\alpha)}_nz^{-n-2}
\end{align*}
where $:\ :$ is the normal ordered product. The modes $\{L^{({\rm B};\alpha)}_n\}_{n\in\mathbb{Z}}$ generate the Virasoso algebra with the central charge fixed to $
1-3\alpha^2.
$
By the energy-momentum tensor $T^{({\rm B};\alpha)}(z)$, each bosonic Fock module $F^{\rm B}_{\beta}(\beta\in\mathbb{C})$ becomes a Virasoro module with
\begin{align}
\label{hbeta}
L^{({\rm B};\alpha)}_0{\mid}\beta;{\rm B}\rangle=h_{\beta}{\mid}\beta;{\rm B}\rangle,\ \ \ \ \ C{\mid}\beta;{\rm B}\rangle=(1-3\alpha^2){\mid}\beta;{\rm B}\rangle
\end{align}
where $h_{\beta}=\frac{1}{2}\beta(\beta-\alpha)$.
Set $T^{({\rm B};\alpha)}:=\frac{1}{2}(a^2_{-1}+\alpha a_{-2}){\mid}0;{\rm B}\rangle$.
The Fock module ${F}^{\rm B}_0$ carries the structure of a vertex operator algebra, with
\begin{align*}
&Y({\mid}0;{\rm B}\rangle,z)={\rm id},\ \ \ Y(a_{-1}{\mid}0;{\rm B}\rangle,z)=a(z),\ \ \ Y(T^{({\rm B};\alpha)},z)=T^{({\rm B};\alpha)}(z).
\end{align*}
We denote by $\mathcal{F}^{\rm B}_{\alpha}$ this vertex operator algebra.

We extend the Heisenberg algebra $\mathfrak{h}$ by a generator $\hat{a}$ satisfying the relations
\begin{align*}
&[a_n,\hat{a}]=\delta_{n,0},
& [\hat{a},\bold{1}]=0.
\end{align*}
For $\beta\in\mathbb{C}$, we define the following vertex operator
\begin{align*}
V_{\beta}(z)=e^{\beta\hat{a}}z^{\beta a_0}\prod_{n\geq 1}\bigl[ {\rm exp}\bigl(\beta\frac{a_{-n}}{n}z^{n}\bigr){\rm exp}\bigl(-\beta\frac{a_{n}}{n}z^{-n}\bigr)\bigr].
\end{align*}
The composition of $k$ vertex operators is given by (cf. \cite{W})
\begin{align}
\label{eq:opeVV}
\begin{aligned}
V_{\beta_1}(z_1)\cdots V_{\beta_k}(z_k)
&=e^{\sum_{i=1}^k\beta_i\hat{a}}\prod_{1\leq i\leq j\leq k}(z_i-z_j)^{\beta_i\beta_j}\prod_{i=1}^kz_i^{\beta_i a_0}\\
&\qquad \cdot \prod_{n\geq 1}\Bigl[ {\rm exp}\bigl(\frac{a_{-n}}{n}\sum_{i=1}^k\beta_i z_i^{n}\bigr){\rm exp}\bigl(-\frac{a_{n}}{n}\sum_{i=1}^k\beta_iz^{-n}_i\bigr)\Bigr].
\end{aligned}
\end{align}
By identifying $e^{\beta\hat{a}}{\mid}\gamma;{\rm B}\rangle={\mid}\beta+\gamma;{\rm B}\rangle$, $V_{\beta}(z)$ becomes a linear map
\begin{align*}
V_{\beta}(z):F_{\gamma}\rightarrow F_{\beta+\gamma}[[z,z^{-1}]]z^{\beta\gamma}
\end{align*}
such that
\begin{equation}
\label{0202V}
\bra{\beta+\gamma;{\rm B}}V_{\beta}(z)\ket{\gamma;{\rm B}}=z^{\beta\gamma},
\end{equation}
where $\bra{\beta;{\rm B}}$ is the bra vector of $(F^{\rm B}_{\beta})^*$.

Since all $F^{\rm B}_{\beta}$ $(\beta\in \mathbb{C})$ are simple $\mathcal{F}^{\rm B}_{\alpha}$-modules, we have ${\rm dim}_{\mathbb{C}}I_{\mathcal{F}^{\rm B}_{\alpha}}${\footnotesize{$\begin{pmatrix}
\ F^{\rm B}_{{\beta+\beta'}} \\
F^{\rm B}_{{\beta}}\ \ F^{\rm B}_{{\beta}'}
\end{pmatrix}$}}$=1$, and ${\rm dim}_{\mathbb{C}}I_{\mathcal{F}^{\rm B}_{\alpha}}${\footnotesize{$\begin{pmatrix}
\ F^{\rm B}_{{\gamma}} \\
F^{\rm B}_{{\beta}}\ \ F^{\rm B}_{{\beta}'}
\end{pmatrix}$}}$=0$ for $\gamma\neq \beta+\beta'$.
For any non-zero intertwining operator $I^{\rm B}_{\beta,\beta'}\in I_{\mathcal{F}^{\rm B}_{\alpha}}${\footnotesize{$\begin{pmatrix}
\ F^{\rm B}_{{\beta+\beta'}} \\
F^{\rm B}_{{\beta}}\ \ F^{\rm B}_{{\beta}'}
\end{pmatrix}$}}, $I_{\beta,\beta'}(\ket{\beta},z)$ is equal to $V_{\beta}(z)$, up to scalar multiples.

\vspace{5mm}

Next, let us introduce the Neveu-Schwarz fermionic Fock module and some notation related to it. The Neveu-Schwarz fermion algebra $\mathfrak{f}$ is the Lie superalgebra
\begin{align*}
\mathfrak{f}=\bigoplus_{r\in\mathbb{Z}+\frac{1}{2}}\mathbb{C} b_r\oplus\mathbb{C}\bold{1}
\end{align*}
with anti-commutation relations
\(
\{b_r,b_s\}=\delta_{r+s,0},
\)
\(\{b_r,\bold{1}\}=0.
\)
The Neveu-Schwarz fermion algebra $\mathfrak{f}$ has the triangular decomposition
\begin{align*}
\mathfrak{f}^{\pm}=\bigoplus_{r>0}\mathbb{C} b_r,\ \ \ \mathfrak{f}^0=\mathbb{C}\bold{1}.
\end{align*}
Let $\mathbb{C}{\mid}{\rm NS}\rangle$ be the one dimensional representation of $\mathfrak{f}^{\geq}=\mathfrak{f}^+\oplus\mathfrak{f}^0$ defined by
\begin{align*}
\bold{1}{\mid}{\rm NS}\rangle={\mid}{\rm NS}\rangle,\ \ \ \mathfrak{f}^+{\mid}{\rm NS}\rangle=0.\\
\end{align*}

\begin{definition}
The {\rm Neveu-Schwarz fermionic Fock module} ${F}^{{\rm NS}}$ is defined by
\begin{eqnarray*}
{F}^{{\rm NS}}={\rm Ind}^{\mathfrak{f}}_{\mathfrak{f}^{\geq}}\mathbb{C}{\mid}{\rm NS}\rangle.\\
\end{eqnarray*}
\end{definition}
Let $b(z)=\sum_{n\in\mathbb{Z}+\frac{1}{2}}b_nz^{-n-\frac{1}{2}}$. Then this field satisfies the operator product expansion
\begin{align}
\label{eq:opebb}
b(z)b(w)=\frac{1}{z-w}+\cdots,
\end{align}
where ``\ $\cdots$'' denote the holomorphic parts about $z=w$.
We define the following energy-momentum tensor
\begin{align*}
T^{(\mathfrak{f})}(z)=\frac{1}{2}:\partial b(z)b(z):=\sum_{n\in\mathbb{Z}}L^{(\mathfrak{f})}_nz^{-n-2}.
\end{align*}
The modes $\{L^{(\mathfrak{f})}_n\}_{n\in\mathbb{Z}}$ generate the Virasoso algebra with the central charge fixed to $\frac{1}{2}$.
By the energy-momentum tensor $T^{(\mathfrak{f})}(z)$, the Neveu-Schwarz fermionic Fock module ${F}^{{\rm NS}}$ becomes a Virasoro module with
\begin{eqnarray*}
L^{(\mathfrak{f})}_0{\mid}{\rm NS}\rangle=0,\ \ \ \ \ C{\mid}{\rm NS}\rangle=\frac{1}{2}{\mid}{\rm NS}\rangle.\\
\end{eqnarray*}

We introduce an even field and an odd field
\begin{align}
\label{eq:TG}
\begin{aligned}
&T^{(\alpha)}(z)=T^{({\rm B};\alpha)}(z)\otimes\bold{1}+\bold{1}\otimes T^{(\mathfrak{f})}(z),\\
&G^{(\alpha)}(z)=a(z)\otimes b(z)+\alpha\bold{1}\otimes \partial b(z).
\end{aligned}
\end{align}
We see that
$T^{(\alpha)}(z)$ and $G^{(\alpha)}(z)$ satisfy the operator product expansions
\begin{align}
\label{eq:opeTG}
\begin{aligned}
&T^{(\alpha)}(z)T^{(\alpha)}(w)=\frac{c_{\alpha}/2}{(z-w)^4}+\frac{2T^{(\alpha)}(w)}{(z-w)^3}+\frac{\partial T^{(\alpha)}(w)}{z-w}+\cdots,\\
&T^{(\alpha)}(z)G^{(\alpha)}(w)=\frac{\frac{3}{2}G^{(\alpha)}(w)}{(z-w)^2}+\frac{\partial G^{(\alpha)}(w)}{z-w}+\cdots,\\
&G^{(\alpha)}(z)G^{(\alpha)}(w)=\frac{2c_{\alpha}/{3}}{(z-w)^3}+\frac{2T^{(\alpha)}(w)}{z-w}+\cdots,
\end{aligned}
\end{align}
where $c_{\alpha}:=\frac{3}{2}-3\alpha^2$.  For the Fourier mode expansions of fields
\begin{align}
\label{eq:TGex}
T^{(\alpha)}(z):=\sum_{n\in\mathbb{Z}}L^{(\alpha)}_nz^{-n-2},\ \ \ G^{(\alpha)}(w):=\sum_{r\in\mathbb{Z}+\frac{1}{2}}G^{(\alpha)}_rz^{-r-\frac{3}{2}},
\end{align}
the modes $\{L^{(\alpha)}_n\}$ and $\{G^{(\alpha)}_r\}$ define the following commutation and anti-commutation relations
\begin{align*}
&[L^{(\alpha)}_k,L^{(\alpha)}_l]=(k-l)L^{(\alpha)}_{k+l}+\delta_{k+l,0}\frac{k^3-k}{12}c_{\alpha},\\
&[L^{(\alpha)}_k,G^{(\alpha)}_r]=(\frac{1}{2}k-r)G^{(\alpha)}_{k+r},\\
&\{G^{(\alpha)}_r,G^{(\alpha)}_s\}=2L^{(\alpha)}_{r+s}+\frac{1}{3}(r^2-\frac{1}{4})\delta_{r+s,0}c_{\alpha}.
\end{align*}
Thus the modes of the fields $T^{(\alpha)}(z)$ and $G^{(\alpha)}(z)$ generate the Neveu-Schwarz algebra with the central charge fixed to $c_{\alpha}=\frac{3}{2}-3\alpha^2$.

\begin{definition}
For $\beta\in\mathbb{C}$, we set 
\begin{align*}
F_{\beta}:={F}^{\rm B}_{\beta}\otimes {F}^{{\rm NS}}
\end{align*}
and call this tensor product {\rm Fock module} simply.
\end{definition}

We set
\begin{align*}
{\mid}\beta\rangle={\mid}\beta;{\rm B}\rangle\otimes {\mid}{\rm NS}\rangle.
\end{align*}
In the following, we omit the tensor product in $a_n\otimes\bold{1}$ and $\bold{1}\otimes b_{r}$ and simply denote them as $a_n$ and $b_r$.
We also use the following shorthand notation
\begin{align*}
V_{\beta}(z)=Y(\ket{\beta},z)=Y({\mid}\beta;{\rm B}\rangle,z)\otimes Y({\mid}{\rm NS}\rangle,z).
\end{align*}

We define the following two vectors in $F_0$
\begin{align*}
&T^{(\alpha)}=\frac{1}{2}(a^2_{-1}+\alpha a_{-2}+b_{-\frac{1}{2}}b_{-\frac{3}{2}}){\mid}0\rangle,
&G^{(\alpha)}=(a_{-1}b_{-\frac{1}{2}}+\alpha b_{-\frac{3}{2}}){\mid}0\rangle.
\end{align*}
The Fock module ${F}_0$ carries the structure of an $N=1$ Neveu-Schwarz vertex operator superalgebra, with
\begin{align*}
&Y({\mid}0\rangle,z)={\rm id},\ \ Y(a_{-1}{\mid}0\rangle,z)=a(z),\ \ Y(b_{-\frac{1}{2}}{\mid}0\rangle,z)=b(z),\\
&Y(G^{(\alpha)},z)=G^{(\alpha)}(z),\ \ \ Y(T^{(\alpha)},z)=T^{(\alpha)}(z).
\end{align*}
We denote by $\mathcal{F}_{\alpha}$ this vertex operator superalgebra.
We use the shorthand notation $T(z)=T^{(\alpha)}(z)$, $G(z)=G^{(\alpha)}(z)$, $L^{(\alpha)}_n$ and $G^{(\alpha)}_{r}$, unless otherwise stated.

Before introducing the structure of Fock modules, let us review the construction of the intertwining operators between Fock modules.
Given $\alpha,\beta,\beta'\in \mathbb{C}$ and $\gamma\in \mathbb{C}$ satisfying $\gamma\neq \beta+\beta'$, we have
\begin{equation*}
\begin{split}
{\rm dim}I^0_{\mathcal{F}_{\alpha}}
\begin{pmatrix}
\ F_{\beta+\beta'} \\
F_{\beta}\ \ F_{\beta'}
\end{pmatrix}
=1,\ \ \ 
{\rm dim}I^1_{\mathcal{F}_{\alpha}}
\begin{pmatrix}
\ F_{\beta+\beta'} \\
F_{\beta}\ \ F_{\beta'}
\end{pmatrix}
=1,\ \ \ 
{\rm dim}I_{\mathcal{F}_{\alpha}}
\begin{pmatrix}
\ F_{\gamma} \\
F_{\beta}\ \ F_{\beta'}
\end{pmatrix}
=0
\end{split}
.
\end{equation*}
For $I^0_{\beta,\beta'}\in I^0_{\mathcal{F}_{\alpha}}
{\footnotesize{\begin{pmatrix}
\ F_{\beta+\beta'} \\
F_{\beta}\ \ F_{\beta'}
\end{pmatrix}}}\setminus\{0\}$ and
$I^1_{\beta,\beta'}\in I^1_{\mathcal{F}_{\alpha}}
{\footnotesize{\begin{pmatrix}
\ F_{\beta+\beta'} \\
F_{\beta}\ \ F_{\beta'}
\end{pmatrix}}}\setminus\{0\}$, we have $I^0_{\beta,\beta'}(\ket{\beta},z)=V_{\beta}(z)$ and $I^1_{\beta,\beta'}(\ket{\beta},z)=b(z)V_{\beta}(z)$ up to scalar multiples.
Thus the non-trivial ${\mathcal{F}_{\alpha}}$-intertwining operators of type 
{\footnotesize{$\begin{pmatrix}
\ F_{{\beta+\beta'}} \\
F_{{\beta}}\ \ F_{{\beta}'}
\end{pmatrix}$}}
can be obtained by the tensor product of the bosonic intertwining operators of type {\footnotesize{$\begin{pmatrix}
\ F^{\rm B}_{{\beta+\beta'}} \\
F^{\rm B}_{{\beta}}\ \ F^{\rm B}_{{\beta}'}
\end{pmatrix}$}} and the Neveu-Schwarz fermion vertex operator.

\subsection{Structure of Fock modules}
Let $m\in \mathbb{Z}_{\geq 1}$. In this subsection we review the structure of Fock modules whose central charges are
\begin{align*}
c=c_{1,2m+1}=\frac{15}{2}-3(2m+1+\frac{1}{2m+1})
\end{align*}
in accordance with the papers \cite{W,IK1,IK2}.

We set
\begin{align}
\label{eq:lattip}
\alpha_{+}=\sqrt{2m+1},\qquad \alpha_-=-\sqrt{\frac{1}{2m+1}},\qquad\alpha_0=\alpha_++\alpha_-.
\end{align}
Note that $c_{\alpha_0}=\frac{3}{2}-3\alpha^2_0=c_{1,2m+1}$.
For $r,s,n\in\mathbb{Z}$, we set
\begin{align}
\label{betars}
&\beta_{r,s;n}:=\frac{1-r}{2}\alpha_++\frac{1-s}{2}\alpha_-+\frac{n}{2}\alpha_+,
&\beta_{r,s}=\beta_{r,s;0},
\end{align}
and we use the shorthand notation $F_{r,s;n}=F_{\beta_{r,s;n}}$ and $F_{r,s}=F_{\beta_{r,s}}$.
For $r,s,n\in\mathbb{Z}$, we introduce the notation
\begin{equation*}
\begin{aligned}
&h_{r,s}:=h_{\beta_{r,s}}=\frac{1}{8}(r^2-1)(2m+1)-\frac{1}{4}(rs-1)+\frac{1}{8}(s^2-1)\frac{1}{2m+1},\\
&h_{r,s;n}:=h_{\beta_{r,s;n}}=h_{r-n,s}=h_{r,s+(2m+1)n}
\end{aligned}
\end{equation*}
and denote by $L(h)$ the simple $\mathfrak{ns}$-module whose minimal conformal weight and central charge are $h$ and $c_{1,2m+1}$.

Before describing the structure of Fock modules, let us introduce the notion of socle series. 
\begin{definition}
Let $V$ be a vertex operator superalgebra or the $\mathfrak{ns}$ algebra, and let $M$ be a finite length $V$-module. 
We denote by ${\rm Soc}(M)$ the socle of $M$, that is ${\rm Soc}(M)$ is the maximal semisimple submodule of $M$. 
Then we have the sequence of the submodules
\begin{align*}
{\rm Soc}_1(M)\subsetneq {\rm Soc}_2(M)\subsetneq \cdots \subsetneq {\rm Soc}_n(M)=M
\end{align*}
such that ${\rm Soc}_1(M)={\rm Soc}(M)$ and ${\rm Soc}_{i+1}(M)/{\rm Soc}_{i}(M)={\rm Soc}(M/{\rm Soc}_{i}(M))$. We call such a sequence of the submodules of $M$ the {\rm socle series} of $M$.
\end{definition}

\begin{proposition}[\cite{IK2}]
\label{socleFockmodule}
For $(r,s)\in \mathbb{Z}^2$ such that $r-s\in2\mathbb{Z}$, the Fock modules $F_{r,s}\in \mathcal{F}_{\alpha_0}\mathchar`-{\rm mod}$ have the following socle series as the $\mathfrak{ns}$-modules:
\begin{enumerate}
\item For each $F_{1,s;n}\ (1\leq s<2m+1,n\in\mathbb{Z},s-n\in2\mathbb{Z}+1)$, we have
\begin{align*}
&{\rm Soc}(F_{1,s;n})=\bigoplus_{k\geq 0}L(h_{1,2m+1-s;{\mid}n{\mid}+2k+1}),\\
&F_{1,s;n}/{\rm Soc}(F_{1,s;n})=\bigoplus_{k\geq a}L(h_{1,s;{\mid}n{\mid}+2k}),
\end{align*}
where $a=0$ if $n\geq 0$, $a=1$ if $n<0$.
\item For each $F_{1,2m+1;2n} (n\in\mathbb{Z})$, we have
\begin{align*}
{\rm Soc}(F_{1,2m+1;2n})=F_{1,2m+1;2n}=\bigoplus_{k\geq 0}L(h_{1,2m+1;{\mid}2n{\mid}+2k}).\\
\end{align*}
\end{enumerate}
\end{proposition}

We introduce the following two fields
\begin{align}
\label{screencurrent}
&Q_+(z):=b(z)V_{\alpha_+}(z),
&Q_-(z):=b(z)V_{\alpha_-}(z).
\end{align}
These fields are the so-called screening currents, which satisfy 
the operator product expansions
\begin{align}
&T(z)Q_{\pm}(w)=\partial_{w}\frac{Q_{\pm}(w)}{z-w}+\cdots,
&G(z)Q_{\pm}(w)=\frac{1}{\alpha_{\pm}}\partial_{w}\frac{V_{\alpha_{\pm}}(w)}{z-w}+\cdots.
\label{sc}
\end{align}
By (\ref{sc}), the operators  
\begin{align*}
&Q_+:=\oint_{z=0}Q_+(z){\rm d}z:F_{{1,2k+1}}\rightarrow F_{{-1,2k+1}}\ \ (k\in\mathbb{Z}),\\
&Q_-:=\oint_{z=0}Q_-(z){\rm d}z:F_{{2k+1,1}}\rightarrow F_{{2k+1,-1}}\ \ (k\in\mathbb{Z})
\end{align*}
become commutative with the $\mathfrak{ns}$-action of $\mathcal{F}_{\alpha_0}\mathchar`-{\rm Mod}$.
These zero-modes of $Q_\pm(z)$ are called screening operators. 

We define fields 
\begin{align*}
Q^{[s]}_-(z):F_{{s+2k,s}}\rightarrow F_{{s+2k,-s}}[[z,z^{-1}]]\ \ (s\geq 2,\ k\in\mathbb{Z}),
\end{align*}
as follows
\begin{align*}
{Q}^{[s]}_{-}(z)=\int_{{\Gamma}_{s}}Q_-(z)Q_-(zx_{1})Q_-(zx_2)\cdots Q_-(zx_{s-1})z^{s-1}{\rm d}x_1\cdots{\rm d}x_{s-1},
\end{align*}
where ${\Gamma}_{s}$ are certain regularized cycles constructed from the simplexes
\begin{align*}
\Delta_{s-1}=\{\ (x_1,\dots,x_{s-1})\in \mathbb{R}^{s-1}\ {\mid}\ 1>x_1>\cdots>x_{s-1}>0\ \}
\end{align*}
(see \cite{TK} for the detatiled construction of the cycles ${\Gamma}_{s}$). Then by the results in \cite{IK2,TK},
the zero-modes
\begin{eqnarray*}
Q^{[s]}_-:=\oint_{z=0}Q^{[s]}_-(z){\rm d}z:F_{{s+2k,s}}\rightarrow F_{{s+2k,-s}}\ \ (s\geq 2,\ k\in\mathbb{Z})
\end{eqnarray*}
are non trivial and commutative with the $\mathfrak{ns}\mathchar`-$action of $\mathcal{F}_{\alpha_0}\mathchar`-{\rm Mod}$.
These fields $Q^{[s]}_{-}(z)$ are called screening currents and the zero-modes $Q^{[s]}_{-}$ are called screening operators. 

We set $Q^{[1]}_{-}:=Q_-$. The structure of the kernels of the screening operators $Q^{[s]}_-$ is given by the following proposition.
\begin{proposition}[\cite{IK2}]
For any $1\leq s\leq 2m$ and $n\in\mathbb{Z}$ such that $s-n$ is odd, let
\begin{align*}
K_{s;n}&={\rm ker}Q^{[s]}_-:F_{1,s;n}\rightarrow F_{1,-s;n}.
\end{align*}
Then we have $K_{s;n}={\rm Soc}(F_{1,s;n})$.  
\label{BRST}
\end{proposition}

\section{The abelian category $\mathcal{SW}(m)$-{\rm mod}}
\label{NStriplet}
In this section,
we introduce the $N=1$ triplet vertex operator superalgebra $\mathcal{SW}(m)$ and review some important results for the abelian category of untwisted $\mathcal{SW}(m)$-modules given in \cite{AM2} (for the twisted sector, see \cite{AM20}).

\subsection{The triplet vertex operator superalgebra $\mathcal{SW}(m)$}
\label{Triplet}
Let $m\in \mathbb{Z}_{\geq 1}$. Let $L=\mathbb{Z}\alpha_+=\mathbb{Z}\sqrt{2m+1}$ be an integral lattice.
\begin{definition}
The {\rm lattice vertex operator superalgebra} $\mathcal{V}_L$ is the quadruple
\begin{equation*}
\Bigl(\bigoplus_{\beta\in L}F_{\beta},{\mid}0\rangle,T,G,Y\Bigr)
\end{equation*}
where the fields corresponding to ${\mid}0\rangle$, $a_{-1}{\mid}0\rangle$, $b_{-\frac{1}{2}}{\mid}0\rangle$, $T$ and $G$ are those of $\mathcal{F}_{\alpha_0}$ and
$
Y({\mid}\beta\rangle,z)=V_{\beta}(z)\ \ (\beta\in L).
$
\end{definition}

For each $i\in\mathbb{Z}$, we introduce the following symbol
\begin{align*}
\gamma_i=\frac{i}{2m+1}\alpha_+=-i\alpha_-.
\end{align*}
It is a known fact that simple $\mathcal{V}_{L}$-modules are given by
\begin{align}
\label{simpleLattice}
\mathcal{V}_{L+\gamma_i}:=\bigoplus_{n\in\mathbb{Z}}F_{\beta_{1,1;2n}+\gamma_i}=\bigoplus_{n\in\mathbb{Z}}F_{1,1+2i;2n}\ \ \ (i=0,\dots,2m).
\end{align}
We define $2m+1$ vector spaces $X_s$ $(1\leq s\leq 2m+1)$ as follows: 
\begin{enumerate}
\item For each $i\in\{0,\dots,m-1\}$, we define
\begin{align}
\label{dfn0830-1}
&X_{2i+1}:={\rm ker}\ Q^{[2i+1]}_-{\mid}_{\mathcal{V}_{L+\gamma_i}},
&X_{2(m-i)}:={\rm ker}\ Q^{[2(m-i)]}_-{\mid}_{\mathcal{V}_{L+\gamma_{2m-i}}}.
\end{align}
\item For $s=2m+1$, we define $X_{2m+1}:=\mathcal{V}_{L+\gamma_m}$.
\end{enumerate}

By Propositions \ref{socleFockmodule}-\ref{BRST}, $X_s$ satisfy the following decomposition as $\mathfrak{ns}$-modules
\begin{equation}
\label{decomp}
\begin{split}
X_{2i+1}&\simeq \bigoplus_{n\in\mathbb{Z}_{\geq 0}}(2n+1)L(h_{1,2i+1;-2n})\ \ (i=0,\dots,m),\\
X_{2(m-j)}&\simeq \bigoplus_{n\in\mathbb{Z}_{\geq 1}}(2n)L(h_{1,2(m-j);-2n+1})\ \ (j=0,\dots,m-1).
\end{split}
\end{equation}

\begin{proposition}[\cite{AM2}]
Let
$
\mathcal{SW}(m)=X_{1}.
$ 
Then $\mathcal{SW}(m)$ has the structure of an $N=1$ Neveu-Schwarz vertex operator superalgebra.
\end{proposition}
This vertex operator superalgebra is called $N=1$ $triplet$ $vertex$ $operator$ $superalgebra$ or $N=1$ $triplet$ $superalgebra$.

\subsection{Simple $\mathcal{SW}(m)$-modules}
\begin{proposition}[\cite{AM2}]
\label{socleV}
\mbox{}
\begin{enumerate}
\item The vector spaces $X_{2i+1}$ $(0\leq i\leq m)$ and $X_{2(m-j)}$ $(0\leq j\leq m-1)$ become simple $\mathcal{SW}(m)$-modules.
\item For each $0\leq i\leq m-1$, the simple $\mathcal{V}_L$-modules $\mathcal{V}_{L+\gamma_i}$ and $\mathcal{V}_{L+\gamma_{2m-i}}$ become $\mathcal{SW}(m)$-modules by restriction, and satisfy the following exact sequences:
\begin{align*}
&0\rightarrow X_{2i+1}\rightarrow \mathcal{V}_{L+\gamma_i}\rightarrow X_{2(m-i)}\rightarrow 0,\\
&0\rightarrow X_{2(m-i)}\rightarrow \mathcal{V}_{L+\gamma_{2m-i}}\rightarrow X_{2i+1}\rightarrow 0.
\end{align*}
\end{enumerate}
\end{proposition}
\begin{remark}
For the simple $\mathcal{SW}(m)$-modules, the notation of \cite{AM2} and ours correspond as follows:
\begin{align*}
&S\Lambda(i+1)=X_{2i+1}\ (0\leq i\leq m),
&S\Pi(m-i)=X_{2(m-i)}\ (0\leq i\leq m-1).
\end{align*}
\end{remark}
We define the following three elements in $\mathcal{V}_L$
\begin{align*}
W^{-}:={\mid}\beta_{1,1;-2}\rangle,\ \ \ W^0:=Q_+W^-,\ \ \ W^+:=Q^{2}_+W^-.
\end{align*}
These elements have the same $L_0$-weight $h_{3,1}=2m+\frac{1}{2}$.
We define the following three elements 
\begin{align}
\label{widehatW}
\widehat{W}^-:=b_{-\frac{1}{2}}{\mid}\beta_{1,1;-2}\rangle,\ \ \ \widehat{W}^0:=Q_+\widehat{W}^-,\ \ \ \widehat{W}^+:=Q^2_+\widehat{W}^-.
\end{align}
These elements have the same $L_0$-weight $2m+1$.

The following three theorems are very important in examining the detailed structure of the module category of $\mathcal{SW}(m)$.
\begin{theorem}[{\cite{AM2}}]
The $N=1$ triplet vertex operator superalgebra $\mathcal{SW}(m)$ is generated by $Y(W^{\pm},z),Y(W^{0},z),G(z)$. Furthermore $\mathcal{SW}(m)$ is strongly generated by
\begin{align*}
G(z),\ T(z),\ Y(W^{\pm},z),\ Y(W^0,z),\ Y(\widehat{W}^{\pm},z),\ Y(\widehat{W}^0,z).
\end{align*}
\end{theorem}
\begin{theorem}[\cite{AM2}]
The $N=1$ triplet vertex operator superalgebra $\mathcal{SW}(m)$ is $C_2$-cofinite.
\label{C_2}
\end{theorem}

\begin{theorem}[\cite{AM2}]
All simple $\mathcal{SW}(m)$-modules are completed by $2m+1$ simple $\mathcal{SW}(m)$-modules in the set $\{X_s\ |\ 1\leq s\leq 2m+1\}$.
\label{simpleclass}
\end{theorem}
The following proposition is straightforward from (\ref{dfn0830-1}), Proposition \ref{BRST} and the definition of $Q_+$.
\begin{proposition}
\label{dfn0830}
\mbox{}
\begin{enumerate}
\item For $0\leq i\leq m$, $n\geq 0$ and $-n \leq k\leq n$, we define 
\begin{equation*}
w^{(n)}_{k}(X_{2i+1}):=
Q^{n+k}_+\ket{\beta_{1,2i+1;-2n}}.
\end{equation*}
Then the set $\{w^{(n)}_{k}(X_{2i+1})\}_{k=-n}^{n}$ gives a basis of the minimal conformal weight spaces of $(2n+1)L(h_{1,2i+1;-2n})\subset X_{2i+1}$.
\item For $0\leq i\leq m-1$, $n\geq 0$ and $-n\leq k\leq n+1$, we define
\begin{equation*}
v^{(n)}_{\frac{2k-1}{2}}(X_{2(m-i)}):=
Q^{n+k}_+\ket{\beta_{1,2(m-i);-2n-1}}.
\end{equation*}
Then the set $\{v^{(n)}_{\frac{2k-1}{2}}(X_{2(m-i)})\}_{k=-n}^{n+1}$ gives a basis of the minimal conformal weight spaces of $(2n+2)L(h_{1,2(m-i);-2n-1})\subset X_{2(m-i)}$.
\end{enumerate}
\end{proposition}
In the following, we use the shorthand notation $w^{(n)}_{k}=w^{(n)}_{k}(X_{2i+1})$ and $v^{(n)}_{\frac{2k-1}{2}}=v^{(n)}_{\frac{2k-1}{2}}(X_{2(m-i)})$.
The transitive $\mathcal{SW}(m)$-actions on the simple $\mathfrak{ns}$-modules of (\ref{decomp}) are given by the following proposition.
\begin{proposition}[\cite{AM2}]
Let $0\leq i\leq m$ and $0\leq j\leq m-1$.
Then the fields $Y(W^{\pm},z)$, $Y(W^{0},z)$, $Y(\widehat{W}^{\pm},z)$ and $Y(\widehat{W}^{0},z)$ act on the vectors $w^{(n)}_{k}\in X_{2i+1}$ and $v^{(n)}_{\frac{2k-1}{2}}\in X_{2(m-j)}$ as follows:
\begin{enumerate}
\item The vectors $w^{(n)}_{k}$ $(n\geq 0,-n\leq k\leq n)$ satisfy
\begin{align}
\label{prop1010}
W^\pm[-h]w^{(0)}_{0}=0,\ \ \ h<h_{1,2i+1;-2}-h_{1,2i+1}=\frac{1}{2}-i+2m,
\end{align}
and
\begin{align*}
\widehat{W}^\pm[0]w^{(n)}_{k}&\in \C^\times w^{(n)}_{k\pm1}+\sum_{l=0}^{n-1}U(\mathfrak{ns})w^{(l)}_{k\pm 1},\\
w^{(n+1)}_{k\pm 1}&\in \C^\times W^\pm[h_{1,2i+1;-2n}-h_{1,2i+1;-2n-2}]w^{(n)}_{k} +\sum_{l=0}^{n}U(\mathfrak{ns})w^{(l)}_{k\pm 1},\\
w^{(n+1)}_{k}&\in \C^\times W^0[h_{1,2i+1;-2n}-h_{1,2i+1;-2n-2}]w^{(n)}_{k}+\sum_{l=0}^{n}U(\mathfrak{ns})w^{(l)}_{k},
\end{align*}
where $w^{(n)}_{n+1}=w^{(n)}_{-n-1}=w^{(-1)}_{k}=0$, $W^\bullet[s]=W^\bullet_s$, $\widehat{W}^\bullet[t]=\widehat{W}^\bullet_t$ $(s,t\in \frac{1}{2}\mathbb{Z})$ and $U(\mathfrak{ns})$ is the universal enveloping algebra of $\mathfrak{ns}$. 
\item  The vectors $v^{(n)}_{\frac{2k+1}{2}}$ $(n\geq 0,-n-1\leq k\leq n)$ satisfy
\begin{equation*}
\begin{split}
W^{\delta}[-h]v^{(0)}_{\frac{\pm1}{2}}
\end{split}
\ \ 
\begin{cases}
=0 &\delta=\pm \\
\in  U(\mathfrak{ns})v^{(0)}_{\frac{\pm1}{2}} &\delta=0\\
\in U(\mathfrak{ns})v^{(0)}_{\frac{\mp1}{2}}& \delta=\mp
\end{cases}
\ \ \ 
\begin{split}
h&<h_{1,2(m-j);-3}-h_{1,2(m-j);-1}\\
&=2m+j+\frac{3}{2},
\end{split}
\end{equation*}
and
\begin{align}
&\widehat{W}^\pm[0]v^{(n)}_{\frac{2k+1}{2}}\in \C^\times v^{(n)}_{\frac{2k+1\pm 2}{2}}+\sum_{l=0}^{n-1}U(\mathfrak{ns})v^{(l)}_{\frac{2k+1}{2}\pm 1},\label{rel20240509}\\
\hspace{-5mm}&v^{(n+1)}_{\frac{2k+1}{2}\pm 1}\in \C^\times W^\pm[h_{1,2(m-j);-2n-1}-h_{1,2(m-j);-2n-3}]v^{(n)}_{\frac{2k+1}{2}}+\sum_{l=0}^{n}U(\mathfrak{ns})v^{(l)}_{\frac{2k+1}{2}\pm 1},\nonumber\\
\hspace{-5mm}&v^{(n+1)}_{\frac{2k+1}{2}}\in \C^\times W^0[h_{1,2(m-j);-2n-1}-h_{1,2(m-j);-2n-3}]v^{(n)}_{\frac{2k+1}{2}}+\sum_{l=0}^{n}U(\mathfrak{ns})v^{(l)}_{\frac{2k+1}{2}}\nonumber,
\end{align}
where $v^{(n)}_{\frac{-2n-3}{2}}=v^{(n)}_{\frac{2n+3}{2}}=v^{(-1)}_{\frac{2k+1}{2}}=0$.
\end{enumerate}
\label{sl2action2}
\end{proposition}

\begin{remark}
The above proposition 
can be shown in a similar way as for the triplet $\mathcal{W}_p$ case \cite{AM}, using the free field realizations and the screening operators (see also \cite{FFL,TW}).
\end{remark}

Let $A(\mathcal{SW}(m))$ be the Zhu-algebra of $\mathcal{SW}(m)$ (for the definition of Zhu-algebras, see \cite{KW},\cite{Zhu}).
For the structure of the Zhu-algebra $A(\mathcal{SW}(m))$, the following theorem holds.
\begin{theorem}[\cite{AM3}]
\label{AM}
The Zhu-algebra $A(\mathcal{SW}(m))$ decomposes as a sum of ideals
\begin{align*}
A(\mathcal{SW}(m))=\bigoplus_{i=2m+1}^{3m}\mathbb{M}_{h_{1,2i+1}}\oplus \bigoplus_{i=0}^{m-1}\mathbb{I}_{h_{1,2i+1}}\oplus \mathbb{C}_{h_{1,2m+1}},
\end{align*}
where $\mathbb{M}_{h_{1,2i+1}}\simeq M_2(\mathbb{C})$, ${\rm dim}(\mathbb{I}_{h_{1,2i+1}})=2$ and ${\rm dim}(\mathbb{C}_{h_{1,2m+1}})=1$.
\end{theorem}
\begin{remark}
Each $\mathbb{M}_{h_{1,2i+1}}$ $(i=2m+1,\dots,3m)$ corresponds to the minimal conformal weight space of the simple module $X_{2(3m+1-i)}$, $\mathbb{C}_{h_{1,2m+1}}$ to the minimal conformal weight space of $X_{2m+1}$ and each $\mathbb{I}_{h_{1,2i+1}}$ $(i=0,\dots,m-1)$ to the minimal conformal weight space of the projective cover of $X_{2i+1}$.
\end{remark}

\subsection{The block decomposition of $\mathcal{SW}(m)\mathchar`-{\rm mod}$}
\label{BlockSW}
Let $\mathcal{SW}(m)\mathchar`-{\rm mod}$ be the abelian category of grading restricted generalised $\mathcal{SW}(m)$-modules (for the definition of grading restricted generalised modules, see Defnition \ref{weak-module}). Since $\mathcal{SW}(m)$ is $C_2$-cofinite, all
objects of $\mathcal{SW}(m)\mathchar`-{\rm mod}$ have finite length \cite{H}.
Note that $\mathcal{SW}(m)\mathchar`-{\rm mod}$ is closed under contragredient.

We denote ${\rm Ext}^n_{\mathcal{SW}(m)}(\bullet,\bullet)$ by the $n$-th ${\rm Ext}$ groups in the abelian category $\mathcal{SW}(m)\mathchar`-{\rm mod}$.
The following proposition can be proved in a similar way as \cite[Theorem 4.4]{AM} by using Theorem \ref{simpleclass} and the results for the semisimple category of $\mathfrak{ns}$ in \cite{IK2}, so we omit the proof.
\begin{proposition}
For all $i\neq j$, we have 
\begin{align*}
&{\rm Ext}^1_{\mathcal{SW}(m)}(X_{2i+1},X_{2j+1})={\rm Ext}^1_{\mathcal{SW}(m)}(X_{2(m-i)},X_{2(m-j)})=0,\\
&{\rm Ext}^1_{\mathcal{SW}(m)}(X_{2i+1},X_{2(m-j)})=0.
\end{align*}
\label{AMsplit}
\end{proposition}
For each $0\leq i\leq m-1$ we denote by $\mathcal{B}_{i+1}$ the full abelian subcategory of $\mathcal{SW}(m)\mathchar`-{\rm mod}$ such that
\begin{align*}
M\in \mathcal{B}_{i+1}\ \Leftrightarrow\ {\rm every\ composition\ factors\ of}\ M\ {\rm are\ given\ by}\ X_{2i+1},X_{2(m-i)}.
\end{align*}
We denote by $\mathcal{B}_{m+1}$ the full abelian subcategory of $\mathcal{SW}(m)\mathchar`-{\rm mod}$ such that
\begin{align*}
M\in \mathcal{B}_{m+1}\ \Leftrightarrow\ {\rm every\ composition\ factors\ of}\ M\ {\rm are\ given\ by}\ X_{2m+1}.
\end{align*}
By Proposition \ref{AMsplit}, we have the following proposition.
\begin{proposition}
\label{Blockdecomp}
The abelian category $\mathcal{SW}(m)\mathchar`-{\rm mod}$ has the following block decomposition
\begin{align*}
\mathcal{SW}(m)\mathchar`-{\rm mod}=\bigoplus_{i=0}^m \mathcal{B}_{i+1}.
\end{align*}
\end{proposition}

\section{Correlation functions}
\label{DF}
In this section, using a certain free field realization technique for vertex operators by \cite{DF1,DF2,Felder,FG} and some properties of the Dotsenko-Fateev integrals \cite{DF1,DF2} given by \cite{S,S2}, we will construct the fundamental system of solutions of a fourth-order differential equation, and determine the connection matrix between the solutions arround $z=0$ and $z=1$. 
The results of this section will be important to show the self-duality of the simple module $X_2$
(see Subsection \ref{X''dual}).

\subsection{A fourth-order Fuchsian differential equation}
\label{FourthO}
Let $\mathcal{Y}_1$ and $\mathcal{Y}_2$ be even $\mathcal{SW}(m)$-intertwining operators of type 
{\footnotesize{$\begin{pmatrix}
\ X_2 \\
X_2\ \ M
\end{pmatrix}$}} and  {\footnotesize{$\begin{pmatrix}
\ M \\
X_2\ \ X_2
\end{pmatrix}$}}, respectively, for some $\mathcal{SW}(m)$-module $M$. 
Fix any minimal conformal weight vector $v_{X_2}\in X_2[h_{2,2}]$ and let $v^*_{X_2}$ be a minimal conformal weight vector of $X^*_2(\simeq X_2)$ such that
$
\langle v^*_{X_2}, v_{X_2}\rangle\neq 0. 
$
We define two correlation functions
\begin{equation}
\label{eq:ou1}
\begin{split}
R^{0}(z_1,z_2)&=\langle v^*_{X_2},\mathcal{Y}_1(v_{X_2},z_1)\mathcal{Y}_2(v_{X_2},z_2)v_{X_2}\rangle,\\
R^{{1}}(z_1,z_2)&=\langle v^*_{X_2},\mathcal{Y}_1(G_{-\frac{1}{2}}v_{X_2},z_1)\mathcal{Y}_2(G_{-\frac{1}{2}}v_{X_2},z_2)v_{X_2}\rangle,
\end{split}
\end{equation}
where $|z_1|>|z_2|>|z_1-z_2|>0$.
Using (\ref{NSrel}) and (\ref{JacobiInt}), we have (cf. \cite[Subsection 2.4]{CMOY})
\begin{equation}
\label{rel0205}
\begin{split}
[G_{-n-\frac{1}{2}},\mathcal{Y}_i(v_{X_2},z)]&={z^{-n}}\mathcal{Y}_i(G_{-\frac{1}{2}}v_{X_2},z),\\
[G_{-n-\frac{1}{2}},\mathcal{Y}_i(G_{-\frac{1}{2}}v_{X_2},z)]&=\bigl({z^{-n}}\partial_{z}-2nz^{-n-1}h_{2,2}\bigr)\mathcal{Y}_i(v_{X_2},z),\\
[L_{-n},\mathcal{Y}_i(v_{X_2},z)]&=\bigl({z^{-n+1}}\partial_{z}+(1-n)z^{-n}h_{2,2}\bigr)\mathcal{Y}_i(v_{X_2},z),\\
[L_{-n},\mathcal{Y}_i(G_{-\frac{1}{2}}v_{X_2},z)]&=\bigl({z^{-n+1}}\partial_{z}+(1-n)z^{-n}(h_{2,2}+\frac{1}{2})\bigr)\mathcal{Y}_i(G_{-\frac{1}{2}}v_{X_2},z).
\end{split}
\end{equation}
We set
\begin{equation}
\label{sing0318}
S_{2,2}:=\frac{4t}{t^2-1}G^4_{-\frac{1}{2}}+\frac{t+1}{t-1}G_{-\frac{1}{2}}G_{-\frac{3}{2}}+\frac{t-1}{t+1}G_{-\frac{3}{2}}G_{-\frac{1}{2}}\ \ \ (t=-2m-1).
\end{equation}
This element of $U(\mathfrak{ns})$ gives an $\mathfrak{ns}$-singular vector of the Verma $\mathfrak{ns}$-module whose lowest weight and central charge are given by $h_{2,2}$ and $c_{1,2m+1}$ (cf. \cite{BA,CMOY}).
Since $U(\mathfrak{ns})v_{X_2}\simeq L(h_{2,2})$, $v_{X_2}$ satisfies the relation $S_{2,2}v_{X_2}=0$. 
Then, similar to the arguments in \cite{HM} (see also \cite{Barron,Belavin}), using (\ref{rel0205}), we can show that $R^{0}(z_1,z_2)$ and $R^{1}(z_1,z_2)$ satisfy the following differential equations
\begin{equation}
\label{kore20231014}
\begin{split}
&\Biggl\{\frac{m^2+2h_{2,2}(m+1)^2}{m(m+1)}\Bigl(\frac{1}{z^2_1}+\frac{1}{z^2_2}\Bigr)-\frac{2m+1}{m(m+1)}(\partial_{z_1}+\partial_{z_2})^2\\
&\ \ \ -\frac{2m^2+2m+1}{m(m+1)}\Bigl(\frac{\partial_{z_1}}{z_1}+\frac{\partial_{z_2}}{z_2}\Bigr)\Biggr\}R^{1}\\
&+\Biggl\{\frac{2m+1}{m(m+1)}\Bigl(\frac{1}{z_2}-\frac{1}{z_1}\Bigr)\partial_{z_1}\partial_{z_2}-\frac{2(2m+1)h_{2,2}}{m(m+1)}\Bigl(\frac{\partial_{z_1}}{z^2_2}-\frac{\partial_{z_2}}{z^2_1}\Bigr)\Biggr\}R^{0}=0,\\[10pt]
&\Biggl\{\frac{2mh_{2,2}}{m+1}\Bigl(\frac{1}{z^2_1}+\frac{1}{z^2_2}\Bigr)-\frac{2m+1}{m(m+1)}(\partial_{z_1}+\partial_{z_2})^2\\
&\ \ -\frac{2m^2+2m+1}{m(m+1)}\Bigl(\frac{\partial_{z_1}}{z_1}+\frac{\partial_{z_2}}{z_2}\Bigr)\Biggr\}R^{0}+\frac{2m+1}{m(m+1)}\Bigl(\frac{1}{z_1}-\frac{1}{z_2}\Bigr)R^{1}=0.
\end{split}
\end{equation}
We define two functions
$
\overline{R}^0(z)=R^0(1,z)
$
and
$
\overline{R}^{1}(z)=R^{1}(1,z).
$
From the $L_0$-conjugation formula for intertwining operators, we have
$
R^{i}(z_1,z_2)=z^{-2h_{2,2}-i}_1\overline{R}^{i}(z_2/z_1)
$, $(i=0,1)$.
Then, from (\ref{kore20231014}), we can show that $\overline{R}^{0}(z)$ and $\overline{R}^{1}(z)$ satisfy 
\begin{equation*}
\begin{split}
&\Bigl\{\frac{d^2}{dz^2}+\frac{(4m^2+2m+1)z-2m^2-2m-1}{(2m+1)z(z-1)}\frac{d}{dz}-\frac{3m^4}{(2m+1)^2z^2(z-1)^2}\Bigr\}\overline{R}^{0}(z)\\
&\ \ \ \ \ \ \ \ \ +\frac{\overline{R}^{1}(z)}{z(1-z)}=0,\\
\end{split}
\end{equation*}
and
\begin{equation*}
\begin{split}
&\Bigl\{\frac{d^2}{dz^2}+\frac{(4m^2+6m+3)z-2m^2-2m-1}{(2m+1)z(z-1)}\frac{d}{dz}\\
&\qquad \qquad\qquad\qquad+\frac{(2m+1)^2z^2-m^2(3m^2+8m+4)}{(2m+1)^2z^2(z-1)^2}\Bigr\}\overline{R}^{1}(z)\\
&\qquad+\Bigl\{\frac{1}{1-z}\frac{d^2}{dz^2}-\frac{(6m^2+2m+1)z-2m-1}{(2m+1)z(z-1)^2}\frac{d}{dz}\\
&\qquad\qquad\qquad\qquad\qquad\qquad\qquad\qquad-\frac{9m^4}{(2m+1)^2z^2(z-1)^2}\Bigr\}\overline{R}^{0}(z)=0.
\end{split}
\end{equation*}
From these equations, we can see that $\overline{R}^{0}(z)$ satisfies the following Fuchsian differential equation
\begin{equation}
\label{fuchsian}
\Bigl(\frac{d^4}{dz^4}+\frac{p_3(z)}{z(z-1)}\frac{d^3}{dz^3}+\frac{p_2(z)}{z^2(z-1)^2}\frac{d^2}{dz^2}+\frac{p_1(z)}{z^3(z-1)^3}\frac{d}{dz}+\frac{p_0(z)}{z^4(z-1)^4}\Bigr)\Phi(z)=0,
\end{equation}
where
\begin{align*}
& p_0(z)=\frac{3m^4}{(2m+1)^4}\bigl((16m^3-8m^2-16m-4)z^2\\
&\ \ \ \ \ \ \ \ \ \ \ \ \ \ \ \ \ \ \ \ +(-16m^3+8m^2+16m+4)z+(3m^4+12m^3+2m^2-4m-1)\bigr),\\
& p_1(z)=\frac{2}{(2m+1)^3}\bigl(   (16m^5+48m^4+56m^3+34m^2+12m+2)z^3\\
&\ \ \ \ \ \ \ \ \ \ \ \ \ \ \ \ \ \ \ \ +(-24m^5-72m^4-84m^3-51m^2-18m-3)z^2\\
&\ \ \ \ \ \ \ \ \ \ \ \ \ \ \ \ \ \ \ \ +(-12m^6-8m^5+12m^3+13m^2+6m+1)z\\
& \ \ \ \ \ \ \ \ \ \ \ \ \ \ \ \ \ \ \ \ +(6m^6+8m^5+12m^4+8m^3+2m^2)\bigr),\\
& p_2(z)=\frac{2}{(2m+1)^2}\bigl((8m^4+32m^3+44m^2+28m+7)z^2\\
& \ \ \ \ \ \ \ \ \ +(-8m^4-32m^3-44m^2-28m-7)z+(-m^4+2m^3+5m^2+4m+1)\bigr),\\
& p_3(z)=\frac{4(m+1)^2(2z-1)}{(2m+1)}.
\end{align*}
The Riemann scheme of the Fuchsian differential equation (\ref{fuchsian}) is given by
\begin{equation}
\label{Scheme}
\begin{split}
\begin{bmatrix}
0 & 1 & \infty \\
\rho_{1,1} & \rho_{1,1}  &0\\
\rho_{0,1} & \rho_{0,1}&\frac{1}{2m+1}\\
\rho_{1,0} & \rho_{1,0} &\frac{4m^2}{2m+1}\\
\rho_{0,0} &  \rho_{0,0} &2m+1
\end{bmatrix}
.
\end{split}
\end{equation}
where
\begin{equation}
\begin{aligned}
&\rho_{1,1}=h_{3,3}-2h_{2,2}=\frac{m^2}{2m+1},
&& \rho_{0,1}=h_{3,1}+\frac{1}{2}-2h_{2,2}=\frac{m^2+4m+1}{2m+1},\\
&\rho_{1,0}=h_{1,3}+\frac{1}{2}-2h_{2,2}=\frac{1-3m^2}{2m+1},
&&\rho_{0,0}=h_{1,1}-2h_{2,2}=-\frac{3m^2}{2m+1}
\end{aligned}
\label{charexp}
\end{equation}
(see \cite{Haraoka} for basic facts of the Fuchsian differential equations).
Note that $\rho_{1,1}-\rho_{1,0},\ \rho_{0,1}-\rho_{0,0}\in \mathbb{Z}_{\geq 1}$. 
In Subsection \ref{ConSol}, we show that (\ref{fuchsian}) has no logarithmic solutions at $z=0,1$. 

Similarly, we can see that the correlation function
\begin{equation}
\label{eq:ou2}
\langle v^*_{X_2},\mathcal{Y}^1(\mathcal{Y}^2(v_{X_2},1-z)v_{X_2},z)v_{X_2}\rangle
\end{equation}
satisfies the Fuchsian differential equation (\ref{fuchsian}), where $\mathcal{Y}^1$ and $\mathcal{Y}^2$ are even $\mathcal{SW}(m)$-intertwining operators of type {\footnotesize{$\begin{pmatrix}
\ X_2 \\
N\ \ X_2
\end{pmatrix}$}} and  {\footnotesize{$\begin{pmatrix}
\ N \\
X_2\ \ X_2
\end{pmatrix}$}}, respectively, for some $\mathcal{SW}(m)$-module $N$.

\begin{remark}
In \cite[Subsection 4.3]{CMOY}, unlike in our cases ((\ref{eq:ou1}), (\ref{eq:ou2})), a slightly different correlation function is examined using a certain embedding technique for Virasoro vertex operator algebras.
\end{remark}
\subsection{Regularization of the Dotsenko-Fateev integrals}
\label{DotsenkoFateev}
For each $l,m,n\in \mathbb{N}$ not all zero, let
\begin{align}
\label{eq:notbox}
\Box^{l,m,n}_{\bm{x}}=[-\infty,0]^l_{x_{1},\dots,x_{l}}\times [0,1]^m_{x_{l+1},\dots,x_{l+m}}\times [1,\infty]^n_{x_{l+m+1},\dots,x_{N}}\subseteqq\overline{\mathbb{R}}^N,
\end{align}
where $N=l+m+n$.
We define
\begin{align}
\label{eq:notint}
J_{l,m,n}(\bm{a},\bm{b},\bm{\gamma})=\int_{\Box^{l,m,n}_{\bm{x}}}\prod_{i=1}^Nx^{a_i}_i(x_i-1)^{b_i}\prod_{1\leq j<k\leq N}(x_k-x_j+i0)^{-2\gamma_{j,k}}{\rm d}x_1\dots{\rm d}x_N,
\end{align}
where $\bm{a}=\{a_i\}_{i=1}^N$, $\bm{b}=\{b_i\}_{i=1}^N$ and $\bm{\gamma}=\{\gamma_{j,k}=\gamma_{k,j}\}_{1\leq j<k\leq N}$.
These integrals are called $Dotsenko$-$Fateev$ $integrals$ \cite{DF1,DF2}.
Note that for $N=1$, $J_{l,m,n}(\bm{a},\bm{b},\bm{\gamma})$ is just Euler Beta integral.
The Dotsenko-Fateev integrals were introduced in \cite{DF1,DF2,Felder,FG} for the approachs to the construction of the BPZ-minimal models \cite{BPZ}. 
In these papers \cite{DF1,DF2,Felder,FG}, the parameters $\gamma_{j,k}$ are fixed as $\gamma_{j,k}=1$, a condition that follows naturally from the free-field realization of correlation functions (see also Subsection \ref{ConSol}).
However, this integer condition makes the problem for the regularization of integrals very difficult. In fact, we see that the construction of cycles in \cite{TK} is not applicable under condition $\gamma_{j,k}=1$.
Thus, the regularization of the Dotsenko-Fateev integrals has been an open problem until a recent result by \cite{S2}.

Let us introduce some notation in accordance with \cite[Section 1]{S2}.
For $0\leq r<1<R\leq \infty$, we define
\begin{equation*}
\begin{split}
\mathcal{M}_{1,0;N}&:=\{\bm{x}\in\mathbb{C}^N \ |\ (x_i\neq x_j,\ i\neq j)\land (x_i\neq 1,0)\},\\
\mathcal{M}_{1,0;N}({r,R})&:=\{(u,v)\in\mathcal{M}_{1,0;N} \ |\ r<|x_i|<{R}\}.
\end{split}
\end{equation*}
The manifold $\mathcal{M}_{1,0;N}({r,R})$ is the moduli space of two pairwise distinct elements of the punctured annulus $\{z\in \mathbb{C}\ |\ r<|x_i|<R,\ x_i\neq 1\}$.
Let $\widetilde{\mathcal{M}}_{1,0;N}(r,R)$ be the universal cover of ${\mathcal{M}}_{1,0;N}(r,R)$. Then  the integrand
\begin{align}
\label{GammaFactor0}
\mathcal{V}_N({{\bm{a},\bm{b}}},\bm{\gamma};\bm{x}):=\prod_{i=1}^Nx^{a_i}_i(x_i-1)^{b_i}\prod_{1\leq j<k\leq N}(x_k-x_j+i0)^{-2\gamma_{j,k}}
\end{align}
of $J_{l,m,n}(\bm{a},\bm{b},\bm{\gamma})$
is a $single$-$valued$ analytic function on the monodromy cover
\begin{equation*}
\begin{split}
\widehat{\mathcal{M}}_{1,0;N}(r,R)&:=\widetilde{\mathcal{M}}_{1,0;N}(r,R)/[\pi_1({\mathcal{M}}_{1,0;N}(r,R)),\pi_1({\mathcal{M}}_{1,0;N}(r,R))].
\end{split}
\end{equation*}

\begin{remark}
The above setting is simplified by considering the case of $N=1$ and $(r,R)=(0,\infty)$.
Let $g$ and $h$ be the generators of $\pi_1(\mathcal{M}_{1,0;1})=\pi_1(\mathbb{C}\setminus\{1,0\})$ corresponding to the one counterclockwise circuit around $0$, $1$, respectively.
Then the element 
\begin{equation*}
g^{-1}h^{-1}gh\in [\pi_1(\mathbb{C}\setminus\{1,0\}),\pi_1(\mathbb{C}\setminus\{1,0\})]
\end{equation*}
defines the {\rm Pochhammer} {\rm contour} arround $1$ and $0$ (Figure \ref{FigPoch}).
Circling along the Pochhammer contour $g^{-1}h^{-1}gh$, we see that the monodromy from $x^a(x-1)^b$ becomes trivial.
Thus the integrand $x^a(x-1)^b$ is well-defined on the quotient space $\widetilde{\mathcal{M}}_{1,0;1}/[\pi_1(\mathbb{C}\setminus\{1,0\}),\pi_1(\mathbb{C}\setminus\{1,0\})]$.
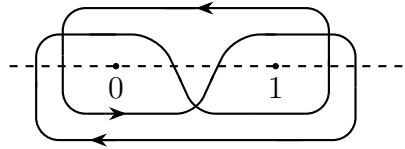
\begin{figure}[h]
  \centering
\begin{tikzpicture}[scale=0.7]
\draw[thick, rounded corners=3mm](6.9,0)--(10.0,0)--(10.0,-1.1);
\draw [->, >=Stealth,thick, ](7.6,0)--(7.5,0);
\draw[thick, rounded corners=3mm](10.0,-1.0)--(10.0,-2.0)--(8.5,-2.0);
\draw[thick, rounded corners=2.5mm](8.5,-2.0)--(7.5,-2.0)--(7.18,-1.3);
\draw[thick, rounded corners=3.5mm](7.22,-1.38)--(6.8,-0.5)--(5.5,-0.5);
\draw[thick, rounded corners=3mm](6.0,-0.5)--(4.5,-0.5)--(4.5,-1.5);
\draw[thick, rounded corners=3mm](4.5,-1.4)--(4.5,-2.5)--(9.5,-2.5);
\draw [-Stealth,thick, ](6.0,-2.5)--(5.5,-2.5);
\draw[thick, rounded corners=3mm](9.0,-2.5)--(10.5,-2.5)--(10.5,-1.5);
\draw[thick, rounded corners=3mm](10.5,-1.5)--(10.5,-0.5)--(8.8,-0.5);
\draw[thick, rounded corners=3.5mm](9.0,-0.5)--(8.2,-0.5)--(7.78,-1.38);
\draw[thick, rounded corners=2.5mm](7.82,-1.3)--(7.5,-2.0)--(5.8,-2.0);
 \draw [-Stealth,thick, ](6.0,-2.0)--(6.2,-2.0);
\draw[thick, rounded corners=3mm](5.8,-2.0)--(5.0,-2.0)--(5.0,-0.9);
\draw[thick, rounded corners=3mm](5.0,-0.9)--(5.0,0.0)--(7.5,0);
\filldraw[black] (9.0,-1.1) circle (0.5mm) node[below]{${1}$};
\filldraw[black] (6.0,-1.1) circle (0.5mm) node[below]{${0}$};
\draw [dashed,thick, ](4,-1.1)--(11.5,-1.1);
\end{tikzpicture}
\caption{The Pochhammer contour arround $1$ and $0$.}
  \label{FigPoch}
\end{figure}
\end{remark}
According to \cite{S2}, the problem of regularizing the Dotsenko–Fateev integral
$J_{l,m,n}(\bm{a},\bm{b},\bm{\gamma})$
can be formulated as finding multi-contours
\begin{align*}
{\Gamma}_{l,m,n}\in H_N(\widehat{\mathcal{M}}_{1,0;N}(r,R);\mathbb{Z})
\end{align*}
satisfying 
\begin{align}
\label{prp-c}
\int_{{\Gamma}_{l,m,n}}\mathcal{V}_N({{\bm{a},\bm{b}}},\bm{\gamma};\bm{x}){\rm d}x_1\wedge\cdots\land{\rm d}x_N\propto J_{l,m,n}(\bm{a},\bm{b},\bm{\gamma})
\end{align}
where \(``\propto\)'' denotes proportionality up to constants and trigonometric functions of \({{\bm{a},\bm{b}}},\gamma\).
In \cite{S2}, these cycles \({\Gamma}_{l,m,n}\) are constructed explicitly, and each proportional constant in (\ref{prp-c}) is determined.
For each $\mathcal{S}\subseteq \{1,\dots,N\}$, let
\begin{align*}
\begin{aligned}
a_{\mathcal{S}}&=\sum_{j\in \mathcal{S}}a_j-2\sum_{\substack{1\leq j<k\leq N\\ j,k\in \mathcal{S}}}\gamma_{j,k},\qquad \qquad b_{\mathcal{S}}=\sum_{j\in \mathcal{S}}b_j-2\sum_{\substack{1\leq j<k\leq N\\ j,k\in \mathcal{S}}}\gamma_{j,k},\\
\zeta_{\mathcal{S}}&=-\sum_{j\in \mathcal{S}}(a_j+b_j)+2\sum_{\substack{1\leq j<k\leq N\\ j\in \mathcal{S}\ {\rm or}\ k\in \mathcal{S}}}\gamma_{j,k}.
\end{aligned}
\end{align*}
\begin{theorem}[\cite{S2}]
\label{twisted00}
There exists a cycle \({\Gamma}_{l,m,n}\in H_N(\widehat{\mathcal{M}}_{1,0;N}(r,R);\mathbb{Z})\) such that
\begin{equation}
\label{0611rem}
\int_{{\Gamma}_{l,m,n}}\mathcal{V}_N({{\bm{a},\bm{b}}},\bm{\gamma};\bm{x}){\rm d}x_1\wedge\cdots\land{\rm d}x_N=c_{l,m,n}(\bm{a},\bm{b},\bm{\gamma}) J_{l,m,n}(\bm{a},\bm{b},\bm{\gamma})
\end{equation}
for all $\bm{a},\bm{b}\in \mathbb{C}^N$ and $\bm{\gamma}\in \mathbb{C}^{\frac{N(N-1)}{2}}$ for which $J_{l,m,n}(\bm{a},\bm{b},\bm{\gamma})$ is defined, 
where by denoting $\mathfrak{e}(x)=1-{\rm exp}(2\pi ix)$, $c_{l,m,n}(\bm{a},\bm{b},\bm{\gamma})$ is given by
\begin{align}
\label{eq:trigc}
\begin{aligned}
c_{l,m,n}=\prod_{\varnothing\subsetneq\mathcal{S}\subseteq \{1,\dots,l+m\}}\mathfrak{e}(a_{\mathcal{S}})
\prod_{\varnothing\subsetneq\mathcal{S}\subseteq \{l+1,\dots,N\}}\mathfrak{e}(b_{\mathcal{S}})
\prod_{\varnothing\subsetneq\mathcal{S}\subseteq \{1,\dots,l\}\cup\{l+m+1,\dots,N\}}\mathfrak{e}(\zeta_{\mathcal{S}}).
\end{aligned}
\end{align}
\end{theorem}

\begin{remark}
\begin{enumerate}
\item 
Let $\mathcal{H}_{l,m,n}$ be the following collection of hyperplanes
\begin{align*}
\mathcal{H}_{l,m,n}=\Bigl(\bigcup_{\mathcal{S}\subseteq \{1,\dots,N\}}\{a_{\mathcal{S}}\in \mathbb{Z}\}\Bigr)
\cup \Bigl(\bigcup_{\mathcal{S}}\{b_{\mathcal{S}}\in \mathbb{Z}\}\Bigr)
\cup\Bigr(\bigcup_{\mathcal{S}}\{\zeta_{\mathcal{S}}\in \mathbb{Z}\}\Bigr).
\end{align*}
Then from the above theorem, 
$J_{l,m,n}$ can be extended to an analytic function
\begin{align*}
J_{l,m,n}:\bigl(\mathbb{C}^N\times \mathbb{C}^N\times \mathbb{C}^{\frac{N(N-1)}{2}}\bigr)\setminus \mathcal{H}_{l,m,n}\rightarrow \mathbb{C}.
\end{align*}
In particular, since $\mathcal{H}_{l,m,n}$ does not contain any affine hyperplanes $\{\gamma_{j,k}=c\}$ $(c\in \mathbb{C})$,
the variables $\gamma_{j,k}$ of $J_{l,m,n}(\bm{a},\bm{b},\bm{\gamma})$ are apparent singularities at $\gamma_{j,k}=1$.
\item Each cycle \({\Gamma}_{l,m,n}\in H_N(\widehat{\mathcal{M}}_{1,0;N}(r,R);\mathbb{Z})\) is defined by a lifting of a cycle in ${\mathcal{M}}_{1,0;N}(r,R)$. These cycles and the lifting are constructed explicitly in \cite[Sections 4 and 5]{S2}.
\end{enumerate}
\end{remark}

In this paper we are mainly interested in the properties of the Dotsenko–Fateev integrals in two dimensional case.
Then in the following, we introduce our notation for $N=2$.
Fix $z_1,z_2\in\mathbb{R}_{>0}$ satisfying $z_2> z_1-z_2>0$. Given $\bm{a}=\{a_1,a_2\}$, $\bm{b}=\{b_1,b_2\}$, $\bm{c}=\{c_1,c_2\}$ $\in \mathbb{C}^2$ and $\gamma\in \mathbb{C}$, we define a new multivalued function
\begin{equation}
\label{GammaFactor}
\begin{split}
&\mathcal{U}(\bm{a},\bm{b},\bm{c},\gamma;u,v;{z_1,z_2})\\
&:=u^{a_1}(u-z_2)^{b_1}(u-z_1)^{c_1}v^{a_2}(v-z_2)^{b_2}(v-z_1)^{c_2}(u-v)^{-2\gamma}
\end{split}
\end{equation}
on 
$
\mathcal{N}_{z_1,z_2,0}=\{(u,v)\in\mathbb{C}^2 \ |\ (u\neq v)\land (u\neq z_1,z_2,0)\land (v\neq z_1,z_2,0)\}.
$
For $i,j=1,0$, let $\Box^\pm_{i,j}$ and $\overline{\Box}^\pm_{i,j}$ be the open subsets of ${\mathcal{N}}_{z_1,z_2,0}$ and ${\mathcal{M}}_{1,0;2}$ defined by
\begin{align*}
&\Box^+_{1,1}:=\{z_1<u,z_1<v\},
&&\Box^+_{0,1}:=\{0<u<z_2,z_1<v\},\\
&\Box^+_{1,0}:=\{z_1<u,0<v<z_2\},
&&\Box^+_{0,0}:=\{0<u<z_2,0<v<z_2\},\\
&\overline{\Box}^+_{1,1}:=\{1<u,1<v\},
&&\overline{\Box}^+_{0,1}:=\{0<u<1,1<v\},\\
&\overline{\Box}^+_{1,0}:=\{1<u,0<v<1\},
&&\overline{\Box}^+_{0,0}:=\{0<u<1,0<v<1\},
\end{align*}
and
\begin{align*}
&\Box^-_{1,1}:=\{u<0,v<0\},
&&\Box^-_{0,1}:=\{z_2<u<z_1,v<0\},\\
&\Box^-_{1,0}:=\{u<0,z_2<v<z_1\},
&&\Box^-_{0,0}:=\{z_2<u<z_1,z_2<v<z_1\},\\
&\overline{\Box}^-_{1,1}:=\{u<0,v<0\},
&&\overline{\Box}^-_{0,1}:=\{0<u<1,v<0\},\\
&\overline{\Box}^-_{1,0}:=\{u<0,0<v<1\},
&&\overline{\Box}^-_{0,0}:=\overline{\Box}^+_{0,0}
\end{align*}
(see Figure \ref{opensub}).
Note that by using the notation (\ref{eq:notbox}), we have
\begin{align*}
&\overline{\Box}^+_{1,1}=\Box^{0,0,2}_{u,v},
&&\overline{\Box}^+_{0,1}=\Box^{0,1,1}_{u,v},
&&&\overline{\Box}^\pm_{0,0}=\Box^{0,2,0}_{u,v},
&&&&\overline{\Box}^+_{1,0}=\Box^{0,1,1}_{v,u},\\
&\overline{\Box}^-_{1,0}=\Box^{1,1,0}_{u,v},
&&\overline{\Box}^-_{1,1}=\Box^{2,0,0}_{u,v},
&&&\overline{\Box}^-_{0,1}=\Box^{1,1,0}_{v,u}.
\end{align*}

\begin{figure}[t]
  \centering
\begin{tikzpicture}[scale=1]
\draw  (0,-1)--(0,3.0);
\draw  (-1,1.0)--(3.0,1.0);
\draw  (-1,2.0)--(3.0,2.0);
\draw  (-1,0)--(3.0,0);
\draw  (1,-1)--(1.0,3);
\draw  (2,-1)--(2.0,3);
\filldraw[black] (0.5,0.5) node{$\Box^+_{0,0}$};
\filldraw[black] (0.5,2.5) node{$\Box^+_{0,1}$};
\filldraw[black] (2.5,0.5) node{$\Box^+_{1,0}$};
\filldraw[black] (2.5,2.5) node{$\Box^+_{1,1}$};
\filldraw[black] (1.5,1.5) node{$\Box^-_{0,0}$};
\filldraw[black] (0.5+1,2.5-3) node{$\Box^-_{0,1}$};
\filldraw[black] (2.5-3,0.5+1) node{$\Box^-_{1,0}$};
\filldraw[black] (2.5-3,2.5-3) node{$\Box^-_{1,1}$};

\draw (3.5, 0) node[inner sep=0.3pt]{\scriptsize${v=0}$};
\filldraw[black] (3.5,1) node{\scriptsize$v=z_2$};
\filldraw[black] (3.5,2) node{\scriptsize$v=z_1$};
\filldraw[black] (0,3.1) node{\scriptsize$u=0$};
\filldraw[black] (1,3.1) node{\scriptsize$u=z_2$};
\filldraw[black] (2,3.1) node{\scriptsize$u=z_1$};

\draw[shift={(6,0)}]  (0,-1)--(0,3.0);
\draw[shift={(6,0)}]  (-1,2.0)--(3.0,2.0);
\draw[shift={(6,0)}]  (-1,0)--(3.0,0);
\draw[shift={(6,0)}]  (2,-1)--(2.0,3);
\filldraw[black,shift={(6,0)}] (1,1) node{$\overline{\Box}^\pm_{0,0}$};
\filldraw[black,shift={(6,0)}] (1.0,2.5) node{$\overline{\Box}^+_{0,1}$};
\filldraw[black,shift={(6,0)}] (2.5,1) node{$\overline{\Box}^+_{1,0}$};
\filldraw[black,shift={(6,0)}] (2.5,2.5) node{$\overline{\Box}^+_{1,1}$};
\filldraw[black,shift={(6,0)}] (1,2.5-3) node{$\overline{\Box}^-_{0,1}$};
\filldraw[black,shift={(6,0)}] (2.5-3,1) node{$\overline{\Box}^-_{1,0}$};
\filldraw[black,shift={(6,0)}] (2.5-3,2.5-3) node{$\overline{\Box}^-_{1,1}$};

\filldraw[black,shift={(6,0)}] (3.5, 0) node[inner sep=0.3pt]{\scriptsize${v=0}$};
\filldraw[black,shift={(6,0)}] (3.5,2) node{\scriptsize$v=1$};
\filldraw[black,shift={(6,0)}] (0,3.1) node{\scriptsize$u=0$};
\filldraw[black,shift={(6,0)}] (2,3.1) node{\scriptsize$u=1$};

\end{tikzpicture}
\caption{The open subsets $\Box^\pm_{i,j}$ and $\overline{\Box}^\pm_{i,j}$.}
  \label{opensub}
\end{figure}

We set 
\begin{equation*}
\begin{split}
{\rm P}_{x,y,z}(u,v)&:=\mathbb{C}[(u-x)^{\pm1},(v-x)^{\pm1},(u-y)^{\pm1},(v-y)^{\pm1},(u-z)^{\pm1},(v-z)^{\pm1}],\\
{\rm P}_{x,y}(u,v)&:={\rm P}_{x,y,y}(u,v),\\
{\rm P}_{x}(u,v)&:={\rm P}_{x,x}(u,v),
\end{split}
\end{equation*}
for $x,y,z\in \mathbb{C}$.
We define 
\begin{align*}
l^{\pm, s}_{i,j},l^{\pm,t}_{i,j},r^{\pm,s}_{i,j},r^{\pm,t}_{i,j}\in \{-\infty,0,z_1,z_2,\infty\},\qquad 
\bar{l}^{\pm,s}_{i,j},\bar{l}^{\pm,t}_{i,j},\bar{r}^{\pm,s}_{i,j},\bar{r}^{\pm,t}_{i,j}\in \{-\infty,0,1,\infty\}
\end{align*}
as follows
\begin{align*}
\Box^\pm_{i,j}&=\{l^{\pm,s}_{i,j}<u<l^{\pm,t}_{i,j},r^{\pm,s}_{i,j}<v<r^{\pm,t}_{i,j}\},\\
\overline{\Box}^\pm_{i,j}&=\{\bar{l}^{\pm,s}_{i,j}<u<\bar{l}^{\pm,t}_{i,j},\bar{r}^{\pm,s}_{i,j}<v<\bar{r}^{\pm,t}_{i,j}\}.\\
\end{align*}
Then we introduce the following two types of integrals.
\begin{definition}
\begin{enumerate}
\item For $i,j\in\{1,0\}$ and $E\in {\rm P}_{z_1,z_2,0}$, we define 
\begin{align*}
&\mathcal{I}^{\pm}_{i,j}[E](\bm{a},\bm{b},\bm{c},\gamma;z_1,z_2)\\
&:=\int_{\Box^\pm_{i,j}}\mathcal{U}(\bm{a},\bm{b},\bm{c},\gamma;u,v;{z_1,z_2})E(u,v){\rm d}u{\rm d}v\\
&:=\int_{l^{\pm,s}_{i,j}}^{l^{\pm,t}_{i,j}}\int_{r^{\pm,s}_{i,j}}^{r^{\pm,t}_{i,j}}u^{a_1}(u-z_2)^{b_1}(u-z_1)^{c_1}v^{a_2}(v-z_2)^{b_2}(v-z_1)^{c_2}\\
&\ \ \ \ \ \ \ \ \ \ \ \ \ \ \ \ \ \times (u-v+i0)^{-2\gamma}E(u,v){\rm d}u{\rm d}v.
\end{align*}
\item For $i,j\in\{1,0\}$ and $F\in {\rm P}_{1,0}$, we define  
\begin{align*}
&\mathcal{J}^{\pm}_{i,j}[F](\bm{a},\bm{b},\gamma)\\
&:=\int_{\overline{\Box}^\pm_{i,j}}\mathcal{V}_2({{\bm{a},\bm{b}}},\gamma;u,v)F(u,v){\rm d}u{\rm d}v\\
&:=\int_{\bar{l}^{\pm,s}_{i,j}}^{\bar{l}^{\pm,t}_{i,j}}\int_{\bar{r}^{\pm,s}_{i,j}}^{\bar{r}^{\pm,t}_{i,j}}u^{a_1}(u-1)^{b_1}v^{a_2}(v-1)^{b_2}(u-v+i0)^{-2\gamma}F(u,v){\rm d}u{\rm d}v,
\end{align*}
where we set $\mathcal{J}^{+}_{0,0}[F](\bm{a},\bm{b},\gamma)=\mathcal{J}^{-}_{0,0}[F](\bm{a},\bm{b},\gamma)$.
\end{enumerate}
\label{dfnIJ}
\end{definition}
We refer to these integrals collectively as the Dotsenko-Fateev integrals.
Applying Theorem \ref{twisted00} to our notation for $N=2$, we have the following theorem.
\begin{theorem}[\cite{S2}]
\label{twisted0}
Let $i,j\in\{1,0\}$. Then for each $(r,R)$ such that $0\leq r<1<R\leq \infty$, there exists a cycle \({\Gamma}^\pm_{i,j}\in H_2(\widehat{\mathcal{M}}_{1,0;2}(r,R);\mathbb{Z})\) such that
\begin{equation}
\label{0611rem}
\int_{{\Gamma}^\pm_{i,j}}\mathcal{V}_2({{\bm{a},\bm{b}}},\gamma;u,v){\rm d}u\land{\rm d}v= c^\pm_{i,j}(\bm{a},\bm{b},\gamma)\mathcal{J}^{\pm}_{i,j}[1](\bm{a},\bm{b},\gamma)
\end{equation}
for all $(\bm{a},\bm{b},\gamma)\in
\mathbb{C}^5_{\bm{a},\bm{b},\gamma}
$ for which $\mathcal{J}^{\pm}_{i,j}[1](\bm{a},\bm{b},\gamma)$ is defined, 
where we set ${\Gamma}^+_{0,0}={\Gamma}^-_{0,0}$ and $c^\pm_{i,j}(\bm{a},\bm{b},\gamma)$ are given by
\begin{align}
\label{eq:trigc}
\begin{aligned}
c^+_{1,1}(\bm{a},\bm{b},\gamma)=&\mathfrak{e}( b_1)\mathfrak{e}( b_2)\mathfrak{e}(b_1+b_2-2\gamma)\mathfrak{e}(a_1+b_1-2\gamma)\mathfrak{e}(a_2+b_2-2\gamma)\\
&\times \mathfrak{e}(a_1+a_2+b_1+b_2-2\gamma),\\
c^-_{1,1}(\bm{a},\bm{b},\gamma)=&c^+_{1,1}(\bm{b},\bm{a},\gamma),\\
c^\pm_{0,0}(\bm{a},\bm{b},\gamma)=&\mathfrak{e}(a_1)\mathfrak{e}( a_2)\mathfrak{e}( b_1)\mathfrak{e}( b_2)\mathfrak{e}(a_1+a_2-2\gamma)\mathfrak{e}(b_1+b_2-2\gamma),\\
c^+_{0,1}(\bm{a},\bm{b},\gamma)=&\mathfrak{e}(a_1)\mathfrak{e}( b_1)\mathfrak{e}( b_2)\mathfrak{e}(b_1+b_2-2\gamma)\mathfrak{e}(a_2+b_2-2\gamma),\\
c^+_{1,0}(\bm{a},\bm{b},\gamma)=&\mathfrak{e}(a_2)\mathfrak{e}(b_1)\mathfrak{e}(b_2)\mathfrak{e}(b_1+b_2-2\gamma)\mathfrak{e}(a_1+b_1-2\gamma),\\
c^-_{0,1}(\bm{a},\bm{b},\gamma)=&c^+_{0,1}(\bm{b},\bm{a},\gamma),\\
c^-_{1,0}(\bm{a},\bm{b},\gamma)=&c^+_{1,0}(\bm{b},\bm{a},\gamma).
\end{aligned}
\end{align}
\end{theorem}

By (\ref{0611rem}), each \(\mathcal{J}^{\pm}_{i,j}[1](\bm{a},\bm{b},\gamma)\) admits an analytic continuation for the variables $(\bm{a},\bm{b},\gamma)$. 
Since the left-hand side of (\ref{0611rem}) is holomorphic with respect to $(\bm{a},\bm{b},\gamma)$, from the explicit form of $c^\pm_{i,j}(\bm{a},\bm{b},\gamma)$, we have the following theorem.
\begin{theorem}[\cite{S,S2}]
Let $i,j=1,0$. 
Then $\mathcal{J}^{\pm}_{i,j}[1](\bm{a},\bm{b},\gamma)$ admits an analytic extension to
$
\mathbb{C}^5_{\bm{a},\bm{b},\gamma}\setminus\mathcal{H}_{\bm{a},\bm{b},\gamma},
$
where 
\begin{equation}
\label{not:hyp}
\begin{split}
\mathcal{H}_{\bm{a},\bm{b},\gamma}=&\{a_1\in \mathbb{Z}_{}\}\cup\{a_2\in \mathbb{Z}_{}\}\cup\{a_1+a_2-2\gamma\in \mathbb{Z}_{}\}\\
&\cup \{b_1\in \mathbb{Z}_{}\}\cup\{b_2\in \mathbb{Z}_{}\}\cup\{b_1+b_2-2\gamma\in \mathbb{Z}_{}\}\\
&\cup \{a_1+b_1\in \mathbb{Z}_{}\}\cup\{a_2+b_2\in \mathbb{Z}_{}\}\\
&\cup\{a_1+a_2+b_1+b_2-2\gamma\in \mathbb{Z}_{}\}.
\end{split}
\end{equation}
\label{thmDF00}
\end{theorem}

Let $z$ and $w$ be real numbers satisfying $z>w\geq 0$. 
Then we set
\begin{equation*}
\mathcal{M}_{z,w;N}:=\{\bm{x}\in\mathbb{C}^N \ |\ (x_i\neq x_j,\ i\neq j)\land (x_i\neq z,w)\}.
\end{equation*}
For $0\leq r<z-w<R\leq \infty$
we define
\begin{equation*}
\begin{split}
\mathcal{M}_{z,w;N}({r,R})&:=\{\bm{x}\in \mathbb{C}^N\ |\ r<|x_i-w|<{R}\}\cap \mathcal{M}_{z,w;N}.
\end{split}
\end{equation*}
Let $\widetilde{\mathcal{M}}_{z,w;N}(r,R)$ be the universal cover of ${\mathcal{M}}_{z,w;N}(r,R)$. Then 
\begin{align*}
\mathcal{V}_N({{\bm{a},\bm{b}}},\bm{\gamma};(z-w)^{-1}(x_1-w),\dots,(z-w)^{-1}(x_N-w))
\end{align*}
is a single-valued analytic function on the monodromy cover
\begin{equation*}
\begin{split}
\widehat{\mathcal{M}}_{z,w;N}(r,R)&:=\widetilde{\mathcal{M}}_{z,w;N}(r,R)/[\pi_1({\mathcal{M}}_{z,w;N}(r,R)),\pi_1({\mathcal{M}}_{z,w;N}(r,R))].
\end{split}
\end{equation*}
Note that there exists a diffeomorphism 
\begin{align*}
t_{z,w;N;r,R}:\mathcal{M}_{z,w;N}({r,R})\rightarrow \mathcal{M}_{1,0;N}({(z-w)^{-1}r,(z-w)^{-1}R})
\end{align*}
defined by
\begin{align*}
t_{z,w;N;r,R}((x_1,\dots,x_N))=((z-w)^{-1}(x_1-w),\dots,(z-w)^{-1}(x_N-w)).
\end{align*}
This diffeomorphism leads to a natural isomorphism
\begin{align*}
 s_{z,w;N;r,R}:H_N(\widehat{\mathcal{M}}_{1,0;N}({(z-w)^{-1}r,(z-w)^{-1}R});\mathbb{Z})\xrightarrow{\simeq} H_N(\widehat{\mathcal{M}}_{z,w;N}(r,R);\mathbb{Z}).
\end{align*}
Then, by Theorem \ref{twisted00}, we have 
\begin{align}
\label{eq:rest-int}
\begin{aligned}
&\int_{s_{z,w;N;r,R}({\Gamma}_{l,m,n})}\mathcal{V}_N({{\bm{a},\bm{b}}},\bm{\gamma};(z-w)^{-1}(x_1-w),\dots,(z-w)^{-1}(x_N-w)){\rm d}\bm{x}\\
&\qquad=(z-w)^{N+\sum_{i=1}^N(a_i+b_i)-2\sum_{j<k}\gamma_{j,k}}c_{l,m,n}(\bm{a},\bm{b},\bm{\gamma}) J_{l,m,n}(\bm{a},\bm{b},\bm{\gamma}).
\end{aligned}
\end{align}
where ${\rm d}\bm{x}={\rm d}x_1\wedge\cdots\land{\rm d}x_N$.
Using these settings, let us define the regularization of the integrals $\mathcal{I}^{\pm}_{i,j}[E](\bm{a},\bm{b},\bm{c},\gamma;z_1,z_2)$.
Since it is sufficient to consider the case where $E\in P_{z_1,z_2,0}$ is a monomial, we can set $E=1$.
First we consider the case $i=j$.
Note that for each $(z,w)\in \{(z_1,z_2),(z_1,0),(z_2,0)\}$,
by choosing $(r,R)$ appropriately, 
the integrand $\mathcal{U}(\bm{a},\bm{b},\bm{c},\gamma;u,v;{z_1,z_2})$ is
a single-valued analytic function on
$
\widehat{\mathcal{M}}_{z,w;2}(r,R).
$
Then for each 
$(i,\pm)$,
by choosing $(r,R,z,w)$ appropriately, 
we can pair the two form $\mathcal{U}(\bm{a},\bm{b},\bm{c},\gamma;u,v;{z_1,z_2}){\rm d}u\wedge{\rm d}v$ with the cycles $s_{z,w;2;r,R}({\Gamma}^\pm_{i,i})$.
Let $\delta>0$ be a small real number satisfying $z_2>z_1-z_2+\delta$ and let
\begin{align*}
\begin{aligned}
{[{\Box}^+_{1,1}]}&:=(c^+_{1,1}(\bm{a}+\bm{b},\bm{c},\gamma))^{-1}s_{z_1,0;2;r,R}({\Gamma}^+_{1,1}),\qquad &&(r,R)=(z_1-\delta,\infty),\\
[{\Box}^+_{0,0}]&:=(c^+_{0,0}(\bm{a},\bm{b},\gamma))^{-1}s_{z_2,0;2;r,R}({\Gamma}^+_{0,0}),\qquad &&(r,R)=(0,z_2+\delta),\\
[{\Box}^-_{1,1}]&:=(c^-_{1,1}(\bm{a},\bm{c}+\bm{b},\gamma))^{-1}s_{z_1,0;2;r,R}({\Gamma}^-_{1,1}),\qquad &&(r,R)=(z_1-\delta,\infty),\\
[{\Box}^-_{0,0}]&:=(c^-_{0,0}(\bm{b},\bm{c},\gamma))^{-1}s_{z_1,z_2;2;r,R}({\Gamma}^-_{0,0}),\qquad &&(r,R)=(0,z_1-z_2+\delta).
\end{aligned}
\end{align*}
Similar to the argument in \cite[Section 6]{S2}, we can show that each pairing
\begin{align}
\label{eq:pairingreg}
\int_{[\Box^\kappa_{i,i}]}\mathcal{U}(\bm{a},\bm{b},\bm{c},\gamma;u,v;{z_1,z_2}){\rm d}u\wedge{\rm d}v,\qquad \kappa=\pm,\ i=1,0
\end{align}
is well-defined and defines the regularization of $\mathcal{I}^{\kappa}_{i,i}[1](\bm{a},\bm{b},\bm{c},\gamma;z_1,z_2)$.
This fact can also be seen from the following argument.
By the Taylor expansion of the integrand \(\mathcal{U}(\bm{a},\bm{b},\bm{c},\gamma;u,v;{z_1,z_2})\)
and by (\ref{eq:rest-int}), we see that
the above pairing (\ref{eq:pairingreg}) admits the following expansion:
\begin{enumerate}
\item For $\kappa=+$, we have
\begin{align}
\label{eq:exp1}
\begin{aligned}
&\int_{[{\Box}^+_{1,1}]}\mathcal{U}(\bm{a},\bm{b},\bm{c},\gamma;u,v;{z_1,z_2}){\rm d}u\wedge{\rm d}v\\
&=z_1^{\sum_{i=1,2}(a_i+b_i+c_i)-2\gamma}\int_{[{\Box}^+_{1,1}]}\mathcal{V}_2(\bm{a}+\bm{b},\bm{c},\gamma;z^{-1}_1u,z^{-1}_1v)\Bigl(\sum_{k\geq 0}z^k_2F^+_{1;k}\Bigr){\rm d}u\wedge{\rm d}v\\
&=z_1^{\sum_{i=1,2}(a_i+b_i+c_i)-2\gamma+2}\sum_{k\geq 0}\Bigl(\frac{z_2}{z_1}\Bigr)^k\mathcal{J}^+_{1,1}[F^+_{1;k}](\bm{a}+\bm{b},\bm{c},\gamma),\\
\vspace{2mm}\\
&\int_{[{\Box}^+_{0,0}]}\mathcal{U}(\bm{a},\bm{b},\bm{c},\gamma;u,v;{z_1,z_2}){\rm d}u\wedge{\rm d}v\\
&=z_2^{\sum_{i=1,2}(a_i+b_i)-2\gamma}z^{c_1+c_2}_1\int_{[{\Box}^+_{0,0}]}\mathcal{V}_2(\bm{a},\bm{b},\gamma;z^{-1}_2u,z^{-1}_2v)\Bigl(\sum_{k\geq 0}z^{-k}_1F^+_{0;k}\Bigr){\rm d}u\wedge{\rm d}v\\
&=z_2^{\sum_{i=1,2}(a_i+b_i)-2\gamma+2}z^{c_1+c_2}_1\sum_{k\geq 0}\Bigl(\frac{z_2}{z_1}\Bigr)^k\mathcal{J}^+_{0,0}[F^+_{0;k}](\bm{a},\bm{b},\gamma),\\
\end{aligned}
\end{align}
where
\begin{align*}
\begin{aligned}
F^+_{1;k}(u,v)&=(k!)^{-1}\left.\partial^k_{z}(1-zu^{-1})^{b_1}(1-zv^{-1})^{b_2}\right|_{z=0},\\
F^+_{0;k}(u,v)&=(k!)^{-1}\left.\partial^k_{z}(1-zu)^{c_1}(1-zv)^{c_2}\right|_{z=0},\\
\end{aligned}
\end{align*}
up to phases, and
we use notation $\bm{x}+\bm{y}=\{x_1+y_1,x_2+y_2\}$ for $\bm{x}=\{x_1,x_2\}$, $\bm{y}=\{y_1,y_2\}$.
\item For $\kappa=-$, we have
\begin{align}
\label{eq:exp2}
\begin{aligned}
&\int_{[{\Box}^-_{1,1}]}\mathcal{U}(\bm{a},\bm{b},\bm{c},\gamma;u,v;{z_1,z_2}){\rm d}u\wedge{\rm d}v\\
&=z_1^{\sum_{i=1,2}(a_i+b_i+c_i)-2\gamma}\int_{[{\Box}^-_{1,1}]}\mathcal{V}_2(\bm{a},\bm{c}+\bm{b},\gamma;z^{-1}_1u,z^{-1}_1v)\\
&\qquad \qquad \qquad \qquad \cdot\Bigl(\sum_{k\geq 0}(z_1-z_2)^kF^-_{1;k}(u+1-z_1,v+1-z_1)\Bigr){\rm d}u\wedge{\rm d}v\\
&=z_1^{\sum_{i=1,2}(a_i+b_i+c_i)-2\gamma+2}\sum_{k\geq 0}\Bigl(\frac{z_1-z_2}{z_1}\Bigr)^k\mathcal{J}^-_{1,1}[F^-_{1;k}](\bm{a},\bm{c}+\bm{b},\gamma),\\
\vspace{2mm}\\
&(z_1-z_2)^{-\sum_{i=1,2}(c_i+b_i)+2\gamma}\int_{[{\Box}^-_{0,0}]}\mathcal{U}(\bm{a},\bm{b},\bm{c},\gamma;u,v;{z_1,z_2}){\rm d}u\wedge{\rm d}v\\
&=z_1^{a_1+a_2}\int_{[{\Box}^-_{0,0}]}\mathcal{V}_2(\bm{b},\bm{c},\gamma;(z_1-z_2)^{-1}(u-z_2),(z_1-z_2)^{-1}(v-z_2))\\
&\qquad \qquad \qquad \qquad \cdot\Bigl(\sum_{k\geq 0}z^{-k}_1F^-_{0;k}(u+1-z_1,v+1-z_1)\Bigr){\rm d}u\wedge{\rm d}v\\
&=(z_1-z_2)^2z_1^{a_1+a_2}\sum_{k\geq 0}\Bigl(\frac{z_1-z_2}{z_1}\Bigr)^k\mathcal{J}^-_{0,0}[F^-_{0;k}](\bm{b},\bm{c},\gamma),\\
\end{aligned}
\end{align}
\end{enumerate}
where
\begin{align*}
\begin{aligned}
F^-_{1;k}(u,v)&=(k!)^{-1}\left.\partial^k_{z}\bigl(1-z(1-u)^{-1}\bigr)^{b_1}\bigl(1-z(1-v)^{-1}\bigr)^{b_2}\right|_{z=0},\\
F^-_{0;k}(u,v)&=(k!)^{-1}\left.\partial^k_{z}\bigl(1-z(1-u)\bigr)^{a_1}\bigl(1-z(1-v)\bigr)^{a_2}\right|_{z=0}.
\end{aligned}
\end{align*}
Since $\mathcal{I}^{\kappa}_{i,i}[1](\bm{a},\bm{b},\bm{c},\gamma;z_1,z_2)$ has the same expansion for an appropriate range of $(\bm{a},\bm{b},\bm{c},\gamma)$, 
the pairing (\ref{eq:pairingreg}) gives the regularization.

Next, let us define the regularization of $\mathcal{I}^{\pm}_{i,j}[1](\bm{a},\bm{b},\bm{c},\gamma;z_1,z_2)$
for the case $i\neq j$.
Note that in this case, the factor $(u-v)^{-2\gamma}$ of the integrand $\mathcal{U}(\bm{a},\bm{b},\bm{c},\gamma;u,v;{z_1,z_2})$ is holomorphic on an appropriate neighborhood of $\Box^\pm_{i,j}$.
Then we define the products of one dimensional cycles
\begin{align*}
\begin{aligned}
{[{\Box}^+_{0,1}]}&:=(d^+_{0,1}(\bm{a},\bm{b},\bm{c},\gamma))^{-1}s_{z_2,0;1;0,z_2+\delta}({\Gamma}_{0,1,0})\times s_{z_1,0;1;z_1-\delta,\infty}({\Gamma}_{0,0,1}),\\
[{\Box}^+_{1,0}]&:=(d^+_{1,0}(\bm{a},\bm{b},\bm{c},\gamma))^{-1}s_{z_1,0;1;z_1-\delta,\infty}({\Gamma}_{0,0,1})\times s_{z_2,0;1;0,z_2+\delta}({\Gamma}_{0,1,0}),\\
[{\Box}^-_{0,1}]&:=(d^-_{0,1}(\bm{a},\bm{b},\bm{c},\gamma))^{-1}s_{z_1,z_2;1;0,z_1-z_2+\delta}({\Gamma}_{0,1,0})\times s_{z_1,0;1;z_1-\delta,\infty}({\Gamma}_{1,0,0}),\\
[{\Box}^-_{1,0}]&:=(d^-_{1,0}(\bm{a},\bm{b},\bm{c},\gamma))^{-1}s_{z_1,0;1;z_1-\delta,\infty}({\Gamma}_{1,0,0})\times s_{z_1,z_2;1;0,z_1-z_2+\delta}({\Gamma}_{0,1,0}),
\end{aligned}
\end{align*}
where
\begin{align}
\label{Int202406112}
\begin{aligned}
&d^+_{0,1}(\bm{a},\bm{b},\bm{c},\gamma)=\mathfrak{e}(a_1)\mathfrak{e}(b_1)\mathfrak{e}(c_2)\mathfrak{e}(a_2+b_2+c_2-2\gamma),\\
&d^+_{1,0}(\bm{a},\bm{b},\bm{c},\gamma)=\mathfrak{e}(a_2)\mathfrak{e}(b_2)\mathfrak{e}(c_1)\mathfrak{e}(a_1+b_1+c_1-2\gamma),\\
&d^-_{0,1}(\bm{a},\bm{b},\bm{c},\gamma)=\mathfrak{e}(b_1)\mathfrak{e}(c_1)\mathfrak{e}(a_2)\mathfrak{e}(a_2+b_2+c_2-2\gamma),\\
&d^-_{1,0}(\bm{a},\bm{b},\bm{c},\gamma)=\mathfrak{e}(b_2)\mathfrak{e}(c_2)\mathfrak{e}(a_1)\mathfrak{e}(a_1+b_1+c_1-2\gamma).
\end{aligned}
\end{align}
Similar to the arguments for the case $i=j$, we can show 
\begin{equation*}
\begin{split}
\int_{[\Box^\pm_{i,j}]}\mathcal{U}(\bm{a},\bm{b},\bm{c},\gamma;u,v;{z_1,z_2}){\rm d}u\wedge{\rm d}v=\mathcal{I}^{\pm}_{i,j}[1](\bm{a},\bm{b},\bm{c},\gamma;z_1,z_2),\quad i\neq j,
\end{split}
\end{equation*}
for all $\bm{a},\bm{b}\in \mathbb{C}^2$ and ${\gamma}\in \mathbb{C}$ for which $\mathcal{I}^{\pm}_{i,j}[1](\bm{a},\bm{b},\bm{c},\gamma;z_1,z_2)$ is defined.

Let $\mathcal{H}_{\bm{a},\bm{b},\bm{c},\gamma}$ be a collection of hyperplanes
\begin{equation}
\label{not:hyp2}
\mathcal{H}_{\bm{a},\bm{b},\bm{c},\gamma}:=\mathcal{H}_{\bm{a},\bm{b},\gamma}\cup\mathcal{H}_{\bm{b},\bm{c},\gamma}\cup\mathcal{H}_{\bm{a},\bm{c},\gamma}\cup \mathcal{H}_{\bm{a}+\bm{b},\bm{c},\gamma}\cup\mathcal{H}_{\bm{b}+\bm{c},\bm{a},\gamma}\cup\mathcal{H}_{\bm{a}+\bm{c},\bm{b},\gamma},
\end{equation}
where we use the notation (\ref{not:hyp}).
From the explicit forms of the trigonometric functions (\ref{eq:trigc}), (\ref{Int202406112}), we obtain the following proposition.
\begin{proposition}
Let $i,j\in\{1,0\}$ and let $z_1,z_2\in\mathbb{R}_{>0}$ be real numbers satisfying $z_2> z_1-z_2>0$. Then for $E(u,v)\in {\rm P}_{z_1,z_2,0}$,
$\mathcal{I}^{\pm}_{i,j}[E](\bm{a},\bm{b},\bm{c},\gamma;z_1,z_2)$ admits an analytic extension to
$
\mathbb{C}^7_{\bm{a},\bm{b},\bm{c},\gamma}\setminus \mathcal{H}_{\bm{a},\bm{b},\bm{c},\gamma}
$.
\label{thmDF0}
\end{proposition}

\subsection{Properties of the Dotsenko-Fateev integrals}
\label{DotsenkoFateev2}
From this subsection, we mainly use the following notation
\begin{equation}
\label{UIJ}
\begin{split}
{U}(a,\rho;u,v;{z_1,z_2})&:=\mathcal{U}(\{a,a'\},\{a,a'\},\{a,a'\},1;u,v;{z_1,z_2}),\\
{I}^{\pm}_{i,j}[E](a,\rho;z_1,z_2)&:=\mathcal{I}^{\pm}_{i,j}[E](\{a,a'\},\{a,a'\},\{a,a'\},1;z_1,z_2),\\
J^{\pm}_{i,j}[F]({a,b},{\rho})&:=\mathcal{J}^{\pm}_{i,j}[F](\{a,a'\},\{b,b'\},1),
\end{split}
\end{equation}
where 
\begin{equation}
\label{rel rho}
a'=-a\rho',\ \ \ \ b'=-b\rho',\ \ \ \ \rho'=1/\rho
\end{equation}
(see (\ref{GammaFactor}) for the notation $\mathcal{U}$).
The condition (\ref{rel rho}) appears naturally in a certain free-field realization of correlation functions (see \cite{DF1,DF2} and Subsection \ref{ConSol}).

\begin{definition}[\cite{S}]
\label{DFSymdef}
Fix $x\in \mathbb{C}$ and $\rho\in \mathbb{C}\setminus\{0,1\}$. We define a $\mathbb{C}$-subalgebra ${\rm P}^{{\rm DF},\rho}_{x}(u,v)$ of ${\rm P}_{x}(u,v)=\mathbb{C}[(u-x)^{\pm1},(v-x)^{\pm1}]$ as follows
\begin{align*}
{\rm P}^{{\rm DF},\rho}_{x}(u,v):=\left.\Bigl\{F(u,v)\in{\rm P}_{x}(u,v) \ \right|\ \bigl(1-\frac{1}{\rho}\bigr)\left.\Bigl(\frac{\partial}{\partial u}F\Bigr)\right|_{u=v}=\frac{\partial}{\partial v}\left.\bigl(F\right|_{u=v}\bigr)\Bigr\}.
\end{align*}
\end{definition}
Let us call the Laurent polynomials in Definition \ref{DFSymdef} ``DF-symmetric polynomials'' according to \cite{S}.
For example, we have
\begin{equation*}
{(u-x)^n}+{(-\rho^{-1})}{(v-x)^n}\in {\rm P}^{{\rm DF},\rho}_{x}(u,v)\ \ \ (n\in \mathbb{Z}).
\end{equation*}
Note that the relation of ${\rm P}^{{\rm DF},\rho}_{x}(u,v)$ is equivalent to
\begin{align*}
\bigl(1-{\rho}\bigr)\left.\Bigl(\frac{\partial}{\partial v}F\Bigr)\right|_{v=u}=\frac{\partial}{\partial u}\left.\bigl(F\right|_{v=u}\bigr).
\end{align*}
Given a DF-symmetric polynomial $F\in {\rm P}^{{\rm DF},\rho}_{x}$ $(x=1,0)$, $J^+_{0,0}[F]({{a,b,\rho}})$ satisfies the following very important property.
\begin{theorem}[\cite{S}]
\label{thmDF}
Fix $x\in \{1,0\}$, $\rho\in \mathbb{C}\setminus\{0,1\}$ and
$F(u,v)\in {\rm P}^{{\rm DF},\rho}_{x}$. Then, for $J^+_{0,0}[F]({{a,b,\rho}})$, there exists an entire function 
$
\tilde{J}[F]:\mathbb{C}^2\rightarrow \mathbb{C}
$
such that
\begin{align*}
J^+_{0,0}[F]({{a,b,\rho}})=\frac{{\rm sin}(\pi (a+b)){\rm sin}(\pi (a'+b'))}{{\rm sin}(\pi a){\rm sin}(\pi b){\rm sin}(\pi a'){\rm sin}(\pi b')}\tilde{J}[F]({{a,b}}),
\end{align*}
where $a'=-\rho^{-1}a$ and $b'=-\rho^{-1}b$ (see the notation (\ref{rel rho})).
\end{theorem}
\vspace{3mm}

The following transformation formulas hold among $J^\pm_{i,j}[F](a,b,\rho)$ $(i,j\in\{1,0\})$.
\begin{theorem}[\cite{DF2,Forrester,S}]
\label{DFint0204}
Fix $x\in \{1,0\}$, $\rho\in \mathbb{C}\setminus\{0,1\}$ and
$F(u,v)\in {\rm P}^{{\rm DF},\rho}_{x}$.
Let $J^\pm_{i,j}=J^\pm_{i,j}[F](a,b,\rho)$ $(i,j\in\{1,0\})$ and $s(x)={\rm sin}(\pi x)$.
Then we have
\begin{align*}
\begin{aligned}
&J^+_{1,0}=\frac{s(a)}{s(a+b)}J^+_{0,0},&&J^-_{1,0}=\frac{s(b)}{s(a+b)}J^+_{0,0},\\
&J^+_{0,1}=\frac{s(a')}{s(a'+b')}J^+_{0,0},&&J^-_{0,1}=\frac{s(b')}{s(a'+b')}J^+_{0,0},\\
&J^+_{1,1}=\frac{s(a)s(a')}{s(a+b)s(a'+b')}J^+_{0,0},&&J^-_{1,1}=\frac{s(b)s(b')}{s(a+b)s(a'+b')}J^+_{0,0},
\end{aligned}
\end{align*}
up to phase factors.
\end{theorem}
The following formula was first given explicitly in \cite{Forrester0}. 
\begin{theorem}[\cite{DF2,Forrester0,Forrester}]
\label{DFint0206}
We have
\begin{equation*}
\begin{split}
J^+_{0,0}[1](a,b,\rho)=&
(\rho')^2\frac{{\rm sin}\pi(a+b)}{{\rm sin}(\pi a)}\frac{\Gamma(\rho'-1)}{\Gamma(\rho')}\frac{\Gamma(1+b)\Gamma(1-a-b)}{\Gamma(-a)}\frac{\Gamma(a')\Gamma(b')}{\Gamma(1+a'+b')},
\end{split}
\end{equation*}
up to a phase factor.
\end{theorem}
These three theorems will be important in constructing the solutions of the Fuchsian differential equation (\ref{fuchsian}). 

\vspace{3mm}
In the following, we introduce transformation and expansion formulas for $I^\pm_{i,j}[1]({{a,\rho}};1,z)$. 
We use the shorthand notation $I^\pm_{i,j}(a,\rho;z)=I^\pm_{i,j}[1]({{a,\rho}};1,z)$ for $i,j=1,0$, and
set 
\begin{equation*}
\begin{split}
\mathcal{H}^{\rm DF}_{a,\rho}:=\mathcal{H}_{\{a,-\rho^{-1}a\},\{a,-\rho^{-1}a\},\{a,-\rho^{-1}a\},1}
\end{split}
\end{equation*}
(see (\ref{not:hyp}) and (\ref{not:hyp2})).  
\begin{lemma}
\label{fourrel}
Fix $(a,\rho)\in \mathbb{C}^2_{a,\rho}\setminus\mathcal{H}^{\rm DF}_{a,\rho}$. We define an involutory matrix
\begin{equation*}
\mathcal{C}
=
\begin{pmatrix}
 (cc')^{-1}& -(cc')^{-1}& -(cc')^{-1}& (cc')^{-1}\\
 (-c+c^{-1})(c')^{-1}& -(cc')^{-1}&(c-c^{-1})(c')^{-1}&(cc')^{-1}\\
 -c^{-1}(c'-(c')^{-1})&c^{-1}(c'-(c')^{-1})&-(cc')^{-1}&(cc')^{-1}\\
 (c-c^{-1})(c'-(c')^{-1})&c^{-1}(c'-(c')^{-1})&(c-c^{-1})(c')^{-1}&(cc')^{-1}
\end{pmatrix}
,
\end{equation*}
where $c=2{\rm cos}(\pi a)$ and $c'=2{\rm cos}(\pi a')$. Then we have
\begin{equation*}
\begin{pmatrix}
I^+_{1,1}(a,\rho;z)\\
I^+_{0,1}(a,\rho;z)\\
I^+_{1,0}(a,\rho;z)\\
I^+_{0,0}(a,\rho;z) 
\end{pmatrix}
=\mathcal{C}
\begin{pmatrix}
I^-_{1,1}(a,\rho;z) \\
I^-_{0,1}(a,\rho;z)\\
I^-_{1,0}(a,\rho;z)\\
I^-_{0,0}(a,\rho;z) 
\end{pmatrix}
.
\end{equation*}

\begin{proof} 
We prove this lemma by similar methods in \cite{DF1}, \cite[Chapter 9]{Haraoka} and \cite[Proposition 3.10]{S}. Let us show the identity
\begin{equation}
\label{idenF}
I^+_{1,1}(a,\rho;z)=(cc')^{-1}\bigl(I^-_{1,1}(a,\rho;z)-I^-_{0,1}(a,\rho;z)-I^-_{1,0}(a,\rho;z)+I^-_{0,0}(a,\rho;z)\bigr).
\end{equation}
The other identities can be proved in the same way. 
Let us derive (\ref{idenF}) using 
$
\mathcal{I}^\pm_{i,j}[1](\{a,a'\},\{a,a'\},\{a,a'\},\gamma;1,z).
$
We use the shorthand notation 
\begin{equation*}
\begin{split}
\mathcal{I}^\pm_{i,j}&=\mathcal{I}^\pm_{i,j}[1](\{a,a'\},\{a,a'\},\{a,a'\},\gamma;1,z),\\
\mathcal{U}(u,v)&=\mathcal{U}(\{a,a'\},\{a,a'\},\{a,a'\},\gamma;u,v;1,z). 
\end{split}
\end{equation*}
Suppose that $a$, $\rho$ and $\gamma$ satisfy $(a,\rho,\gamma)\in \mathbb{C}^3_{a,\rho,\gamma}\setminus \mathcal{H}_{\{a,\rho'a\},\{a,\rho'a\},\{a,\rho'a\},\gamma}$ and
\begin{align}
\label{0215eq}
-\frac{1}{2}&<{\rm Re}(a)\ll0,\ \ \ \ -\frac{1}{2}<{\rm Re}(a')\ll0,\ \ \ \ 0<{\rm Re}(\gamma)\ll1.
\end{align}
Let $C^{+}$ be a counterclockwise simple closed contour in the upper half plane $\{ u\in\mathbb{C}\ |\ {\rm Im}(u)>0\}$, and let $C^{-}$ be a clockwise simple closed contour in the lower half plane $\{ u\in\mathbb{C}\ |\ {\rm Im}(u)<0\}$. 
For $\mathcal{U}(u,v)$ with $v\in \mathbb{R}_{>1}$ fixed, we take the branch cut along the real axis $\{u\in \mathbb{R}\ |\ u\leq v\}$.
Then by noting the condition (\ref{0215eq}), from the Cauchy's theorem, we obtain
\begin{equation*}
\begin{split}
0&=e^{2i\pi \gamma}\oint_{C^+}\Bigl(\int_{1}^{\infty}\mathcal{U}(u+i0,v){\rm d}v\Bigr){\rm d}u\\
&=e^{2i\pi \gamma}\mathcal{I}^+_{1,1}+e^{i\pi a}\int_{z}^{1}\int_{1}^{\infty}\mathcal{U}{\rm d}v{\rm d}u+e^{2i\pi a}\int_{0}^z\int_{1}^{\infty}\mathcal{U}{\rm d}v{\rm d}u+e^{3i\pi a}\int_{-\infty}^{0}\int_{1}^{\infty}\mathcal{U}{\rm d}v{\rm d}u
\end{split}
\end{equation*}
and
\begin{equation*}
\begin{split}
0&=e^{-2i\pi \gamma}\oint_{C^-}\Bigl(\int_{1}^{\infty}\mathcal{U}(u-i0,v){\rm d}v\Bigr){\rm d}u\\
&=e^{-6i\pi \gamma}\mathcal{I}^{+}_{1,1}{+}e^{-i\pi a}\int_{z}^{1}\int_{1}^{\infty}\mathcal{U}{\rm d}v{\rm d}u{+}e^{-2i\pi a}\int_{0}^z\int_{1}^{\infty}\mathcal{U}{\rm d}v{\rm d}u{+}e^{-3i\pi a}\int_{-\infty}^{0}\int_{1}^{\infty}\mathcal{U}{\rm d}v{\rm d}u.
\end{split}
\end{equation*}
From these identities, we have
\begin{equation}
\label{e0216}
\begin{split}
e^{-2i\pi \gamma}s(2a-4\gamma)\mathcal{I}^+_{1,1}=s(a)\int_{1}^{z}\int_{1}^{\infty}\mathcal{U}{\rm d}v{\rm d}u-s(a)\int_{0}^{-\infty}\int_{1}^{\infty}\mathcal{U}{\rm d}v{\rm d}u,
\end{split}
\end{equation}
 where $s(x)={\rm sin}(\pi x)$. 
Similarly, 
we can show the identities
\begin{equation*}
\begin{split}
s(2(a'-\gamma))\int_{1}^{z}\int_{1}^{\infty}\mathcal{U}{\rm d}v{\rm d}u&=e^{-2i\pi \gamma}s(a'-4\gamma)\mathcal{I}^-_{0,0}-s(a')\mathcal{I}^-_{0,1}\\[10pt]
s(2a')\int_{0}^{-\infty}\int_{1}^{\infty}\mathcal{U}{\rm d}v{\rm d}u&=s(a')\mathcal{I}^-_{1,0}-e^{-2i\pi \gamma}s(a'-2\gamma)\mathcal{I}^-_{1,1}.
\end{split}
\end{equation*}
Thus, from these identities and (\ref{e0216}), we obtain the desired identity (\ref{idenF}) by analytically continuing $\mathcal{I}^\pm_{i,j}$ to $\gamma=1$.
\end{proof}
\end{lemma}

From the expansions (\ref{eq:exp1})-(\ref{eq:exp2}), we obtain the following lemma.
\begin{lemma}
\label{conv}
Fix $(a,\rho)\in \mathbb{C}^2_{a,\rho}\setminus\mathcal{H}^{\rm DF}_{a,\rho}$. Then $I^\pm_{i,j}(a,\rho;z)$ admits the following expansion.
\begin{enumerate}
\item The functions $z^{2(i-1)(a+a')}I^+_{i,i}(a,\rho;z)$ and $(1-z)^{2(i-1)(a+a')}I^-_{i,i}(a,\rho;z)$ $(i=1,0)$ admit analytic continuations to the complex domains $\{|z|<1\}$ and $\{|z-1|<1\}$, respectively, and on these domains, satisfy convergent series
\begin{equation*}
\begin{split}
z^{2(i-1)(a+a')}I^+_{i,i}(a,\rho;z)&=\sum_{k\geq 0}J^+_{i,i}[F^+_{i;k}]((i+1)a,a,\rho)z^k,\\
(1-z)^{2(i-1)(a+a')}I^-_{i,i}(a,\rho;z)&=\sum_{k\geq 0}J^-_{i,i}[F^-_{i;k}](a,(i+1)a,\rho)(1-z)^k,
\end{split}
\end{equation*}
where 
$F^+_{1;k}(u,v)$ is 
defined by
\begin{align*}
F^+_{1;k}(u,v)=(k!)^{-1}\left.\partial^k_{z}\bigl((1-zu^{-1})^{a}(1-zv^{-1})^{a'}\bigr)\right|_{z=0},
\end{align*}
and
\begin{align*}
\begin{aligned}
F^+_{0;k}(u,v)&=F^+_{1;k}(u^{-1},v^{-1}),
&F^-_{1;k}(u,v)&=F^+_{1;k}(u-1,v-1),\\
F^-_{0;k}(u,v)&=F^+_{1;k}((1-u)^{-1},(1-v)^{-1}).
\end{aligned}
\end{align*}
up to phases.
\item For $i\neq j$, $z^{2(i-1)a+2(j-1)a'}I^+_{i,j}(a,\rho;z)$ and $(1-z)^{2(i-1)a+2(j-1)a'}I^-_{i,j}(a,\rho;z)$ admit analytic continuations to the complex domains $\{|z|<1\}$ and $\{|z-1|<1\}$, respectively, and on these domains, satisfy convergent series
\begin{equation*}
\begin{split}
&\sum_{k\geq 0}\mathcal{J}^+_{i,j}[F^+_{\{i,j\};k}](\{(i+1)a-2i,(j+1)a'-2j\},\{a,a'\},0)z^{1+k},\\
&\sum_{k\geq 0}\mathcal{J}^-_{i,j}[F^-_{\{i,j\};k}](\{a,a'\},\{(i+1)a-2i,(j+1)a'-2j\},0)(1-z)^{1+k},
\end{split}
\end{equation*}
respectively, where 
by setting
\begin{align*}
G_k(x_1,x_2)=(k!)^{-1}\left.\partial^k_{z}\bigl((1-zx_1)^{a}(1-zx_2)^{a'}(1-zx_1x_2)^{-2}\bigr)\right|_{z=0},
\end{align*}
$F^\pm_{\{i,j\};k}$ are defined by 
\begin{align*}
\begin{aligned}
&F^+_{\{0,1\};k}(u,v)=G_k(u,v^{-1}),&F^+_{\{1,0\};k}(u,v)&=G_k(u^{-1},v),\\
&F^-_{\{0,1\};k}(u,v)=G_k(1-u,(1-v)^{-1}),&F^-_{\{1,0\};k}(u,v)&=G_k((1-u)^{-1},1-v)
\end{aligned}
\end{align*}
up to phases.
\end{enumerate}
\end{lemma}
\begin{remark}
\label{rem:dfsymex}
\begin{enumerate}
\item We can see that
$
\left.\partial^n_{z}\bigl((1-zu)^a(1-zv)^{a'}\bigr)\right|_{z=0}\in {\rm P}^{{\rm DF},\rho}_{0}.
$
Then we have $F^+_{1;k},F^+_{0;k}\in {\rm P}^{{\rm DF},\rho}_{0}$ and $F^-_{1;k},F^-_{0;k}\in {\rm P}^{{\rm DF},\rho}_{1}$.
\item From Definition 
\ref{dfnIJ} and the formulas of the beta integrals, we have
\begin{equation*}
\begin{split}
\hspace{-8mm}\mathcal{J}^\pm_{0,1}[1](\bm{a},\bm{b},0)&=\frac{\Gamma(a_1+1)\Gamma(b_1+1)\Gamma(-a_2-b_2-1)\Gamma(b_2+1)^{\frac{1\pm1}{2}}\Gamma(a_2+1)^{\frac{1\mp 1}{2}}}{\Gamma(a_1+b_1+2)\Gamma(-a_2)^{\frac{1\pm1}{2}}\Gamma(-b_2)^{\frac{1\mp 1}{2}}},\\
\hspace{-8mm}\mathcal{J}^\pm_{1,0}[1](\bm{a},\bm{b},0)&=\frac{\Gamma(a_2+1)\Gamma(b_2+1)\Gamma(-a_1-b_1-1)\Gamma(b_1+1)^{\frac{1\pm1}{2}}\Gamma(a_1+1)^{\frac{1\mp 1}{2}}}{\Gamma(a_2+b_2+2)\Gamma(-a_1)^{\frac{1\pm1}{2}}\Gamma(-b_1)^{\frac{1\mp 1}{2}}},
\end{split}
\end{equation*}
up to phase factors (see \cite[Appendix A]{S2}).
\end{enumerate}
\end{remark}

\subsection{Construction of solutions}
\label{ConSol}
We want to use the Dotsenko-Fateev integrals $I^\pm_{i,j}[1]({{a,\rho}};1,z)$ $(i,j=1,0)$ with parameters $a=\alpha_-\beta_{2,2}$, $\rho=\alpha_-^2$  (for the definition of $I^\pm_{i,j}[1]({{a,\rho}};1,z)$, see (\ref{UIJ})). However, in this case, we have $a'=\alpha_+\beta_{2,2}=-m\in \mathbb{Z}$, and then the integral cannot be defined. 
So let us add a small complex parameter to the variables $(a,\rho)$ to make the integrals well-defined.
Before that, let us define some notation. We set
\begin{equation*}
\mathbb{D}(x):=\{x\in \mathbb{C}\ |\ |x|<1\},\ \ \ \ 
\mathbb{D}^\times(x):=\mathbb{D}(x)\setminus\{x=0\}.
\end{equation*}
We also use the shorthand notation $\mathbb{D}=\mathbb{D}(x)$ and $\mathbb{D}^\times=\mathbb{D}^\times(x)$.
For any non-empty open set $U\in \mathbb{C}^n$, we define
\begin{equation*}
\begin{split}
\mathcal{O}(U)&:=\{f\in F(U)\ |\ f\ {\rm is}\ {\rm holomorphic}\ {\rm on}\ U\},
\end{split}
\end{equation*}
where $F(U)$ is the set of complex functions on $U$.
\begin{remark}
In the above definition, we say that $f(z_1,\dots,z_n)\in F(U)$ is holomorphic on $U\in \mathbb{C}^n$, when $f(z_1,\dots,z_n)$ is holomorphic for each variable $z_i$ $((z_1,\dots,z_n)\in U)$.
\end{remark}
Following \cite{TW}, we introduce $\epsilon$-deformations of $\alpha_\pm$, $\alpha_0$ and ${\beta}_{r,s}$ $(r,s\in \mathbb{Z})$ as follows (see (\ref{eq:lattip})-(\ref{betars}))
\begin{align*}
\begin{aligned}
&\widetilde{\alpha}_{+}(\epsilon)=\alpha_++\theta{\epsilon},&&\widetilde{\alpha}_-(\epsilon)=-\frac{1}{\widetilde{\alpha}_{+}(\epsilon)},\\
&\widetilde{\alpha}_0(\epsilon)=\widetilde{\alpha}_+(\epsilon)+\widetilde{\alpha}_-(\epsilon),&&\widetilde{\beta}_{r,s}=\frac{(1-r)}{2}\widetilde{\alpha}_++\frac{(1-s)}{2}\widetilde{\alpha}_-,
\end{aligned}
\end{align*}
where we fix a sufficiently small $\theta\in \mathbb{R}_{>0}$ to satisfy the following condition
\begin{equation}
\label{eq:condition}
(a,\rho)=(\widetilde{\alpha}_-(\epsilon)\widetilde{\beta}_{2,2}(\epsilon),\widetilde{\alpha}_-(\epsilon)^2)\notin \mathcal{H}^{\rm DF}_{a,\rho}\ \ \ \ \ (\epsilon\in \mathbb{D}^\times(\epsilon)).
\end{equation}
Note that $\widetilde{\alpha}_\pm,\widetilde{\alpha}_0,\widetilde{\beta}_{r,s}\in \mathcal{O}(\mathbb{D})$, $\widetilde{\alpha}_\pm(0)={\alpha}_\pm$, $\widetilde{\alpha}_0(0)={\alpha}_0$ and $\widetilde{\beta}_{r,s}(0)={\beta}_{r,s}$.
Consider the integrals $I^\pm_{i,j}[1](\widetilde{\alpha}_-(\epsilon)\widetilde{\beta}_{2,2}(\epsilon),\widetilde{\alpha}_-(\epsilon)^2;1,z)$ for $i,j=1,0$. 
Since (\ref{eq:condition}), from Proposition \ref{thmDF0}, we see that these integrals are well-defined for $\epsilon\in \mathbb{D}^\times$.
We choose branches of $I^\pm_{i,j}[1](\widetilde{\alpha}_-(\epsilon)\widetilde{\beta}_{2,2}(\epsilon),\widetilde{\alpha}_-(\epsilon)^2;1,z)$ at $z=\frac{1}{2}$ as
$
{\rm arg}z={\rm arg}(1-z)=0.
$
Then by the analytic continuations, $I^\pm_{i,j}[1](\widetilde{\alpha}_-(\epsilon)\widetilde{\beta}_{2,2}(\epsilon),\widetilde{\alpha}_-(\epsilon)^2;1,z)$ are single-valued functions on the simply-connected domain $\mathbb{D}(z)\cap \mathbb{D}(1-z)=\{|z|<1\}\cap \{|z-1|<1\}$.
\begin{proposition}
\begin{enumerate}
\item For $i,j=1,0$, we have
\begin{equation*}
\epsilon I^\pm_{i,j}[1](\widetilde{\alpha}_-(\epsilon)\widetilde{\beta}_{2,2}(\epsilon),\widetilde{\alpha}_-(\epsilon)^2;1,z)\in \mathcal{O}((\mathbb{D}(z)\cap \mathbb{D}(1-z))\times\mathbb{D}(\epsilon)).
\end{equation*}

\item Let $C_0\subset \mathbb{D}^\times(\epsilon)$ be a simple closed curve arround $\epsilon=0$.
Then for $z\in \mathbb{D}(z)\cap\mathbb{D}(1-z)$ and $i,j=1,0$, we have
\begin{equation*}
\begin{split}
\int_{C_{0}}I^\pm_{i,j}[1](\widetilde{\alpha}_-(\epsilon)\widetilde{\beta}_{2,2}(\epsilon),\widetilde{\alpha}_-(\epsilon)^2;1,z){\rm d}\epsilon\neq 0
\end{split}
,
\end{equation*}
and asymptotic behaviors
\begin{equation}
\label{asympt}
\int_{C_{0}}I^\kappa_{i,j}[1](\widetilde{\alpha}_-(\epsilon)\widetilde{\beta}_{2,2}(\epsilon),\widetilde{\alpha}_-(\epsilon)^2;1,z){\rm d}\epsilon\sim
\begin{cases}
C^+_{i,j}z^{\lambda_{i,j}}&(\kappa=+,\ z\rightarrow 0),\\
C^-_{i,j}(1-z)^{\lambda_{i,j}}&(\kappa=-,\ z\rightarrow 1),
\end{cases}
\end{equation}
where $\lambda_{i,j}=2(1-i)\alpha_-\beta_{2,2}+2(1-j)\alpha_+\beta_{2,2}+(1-i)j+(1-j)i$ and $C^\pm_{i,j}$ are some non-zero constants.
\end{enumerate}
\begin{proof}
We only prove the case $i=j$. The cases $i\neq j$ can be proved in the same way.
We use the shorthand notation
\begin{align*}
\begin{aligned}
J^+_{i,i}[F^+_{i;k}](\epsilon)&=J^+_{i,i}[F^+_{i;k}]((i+1)a,a,\rho),
&&J^-_{i,i}[F^-_{i;k}](\epsilon)=J^-_{i,i}[F^-_{i;k}](a,(i+1)a,\rho),\\
I^\pm_{i,j}(z,\epsilon)&=I^\pm_{i,j}[1](a,\rho;1,z),
\end{aligned}
\end{align*}
as $a=\widetilde{\alpha}_-(\epsilon)\widetilde{\beta}_{2,2}(\epsilon)$, $a'=\widetilde{\alpha}_+(\epsilon)\widetilde{\beta}_{2,2}(\epsilon)$ and $\rho=\widetilde{\alpha}_-(\epsilon)^2$, where $F^+_{1;k},F^+_{0;k}\in {\rm P}^{\rm DF,\widetilde{\alpha}_-(\epsilon)^2}_0$ and $F^-_{1;k},F^-_{0;k}\in {\rm P}^{\rm DF,\widetilde{\alpha}_-(\epsilon)^2}_1$ are defined in Lemma \ref{conv} (see also Remark \ref{rem:dfsymex}).

By Proposition \ref{thmDF0} and by Lemma \ref{conv}, we see that
\begin{equation}
\label{HolE-1}
I^\pm_{i,i}(z,\epsilon)\in \mathcal{O}((\mathbb{D}(z)\cap \mathbb{D}(1-z))\times\mathbb{D}^\times(\epsilon))\ \ \ (i=1,0).
\end{equation}
By Theorems \ref{thmDF00}, \ref{thmDF} and \ref{DFint0204}, we see that 
\begin{equation}
\label{DFJint}
\epsilon J^+_{i,i}[F^+_{i;k}](\epsilon),\ \epsilon J^-_{i,i}[F^-_{i;k}](\epsilon)\in \mathcal{O}(\mathbb{D}(\epsilon))\ \ \ \ (k\geq 0,\ i=1,0).
\end{equation}
Furthermore, by Theorems \ref{DFint0204}-\ref{DFint0206}, we can see that
\begin{equation}
\label{DFJint2}
{\rm Res}_{\epsilon=0} J^\pm_{i,j}[1](\epsilon)\neq 0\ \ \ (i,j=1,0).
\end{equation}
Then by Lemma \ref{conv} and by (\ref{HolE-1})-(\ref{DFJint}), we see that $\epsilon I^\pm_{i,i}(z,\epsilon)$ $(i=1,0)$ are holomorphic for $\epsilon\in \mathbb{D}(\epsilon)$.

Let $C_0\subset \mathbb{D}^\times(\epsilon)$ be a simple closed curve arround $\epsilon=0$. By Lemma \ref{conv} and by (\ref{DFJint})-(\ref{DFJint2}), we see that $\int_{C_0}I^\pm_{i,i}(z,\epsilon){\rm d}\epsilon\neq 0$ $(i=1,0)$ and $\int_{C_0}I^\pm_{i,i}(z,\epsilon){\rm d}\epsilon$ satisfy the asymptotic behavior (\ref{asympt}). 
\end{proof}
\label{lemDf}
\end{proposition}

We introduce $\epsilon$-deformations of the screening currents (\ref{screencurrent}) as follows
\begin{align*}
Q^{(\epsilon)}_\pm(z)=b(z)V_{\widetilde{\alpha}_\pm(\epsilon)}(z).
\end{align*}
We set $T^{(\epsilon)}(z)=T^{(\widetilde{\alpha}_0(\epsilon))}(z)$, $G^{(\epsilon)}(z)=G^{(\widetilde{\alpha}_0(\epsilon))}(z)$, 
$L^{(\epsilon)}_n=L^{(\widetilde{\alpha}_0(\epsilon))}_n$ and $G^{(\epsilon)}_r=G^{(\widetilde{\alpha}_0(\epsilon))}_r$
(see (\ref{eq:TG})-(\ref{eq:TGex}) for the right-hand notation).
From the definition of these operators, we have the operator product expansions
\begin{align}
&T^{(\epsilon)}(z)Q^{(\epsilon)}_{\pm}(w)=\partial_{w}\frac{Q^{(\epsilon)}_{\pm}(w)}{z-w}+\cdots,
&G^{(\epsilon)}(z)Q^{(\epsilon)}_{\pm}(w)=\frac{1}{\widetilde{\alpha}_{\pm}(\epsilon)}\partial_{w}\frac{V_{\widetilde{\alpha}_{\pm}(\epsilon)}(w)}{z-w}+\cdots.
\label{sc0}
\end{align}
Then the operators \(Q^{(\epsilon)}_\pm(z)\) define the screening currents for \(U^{(\epsilon)}(\mathfrak{ns})\), where $U^{(\epsilon)}(\mathfrak{ns})$ is the universal enveloping algebra of $\mathfrak{ns}$ with the central charge fixed to $c_{\widetilde{\alpha}_0(\epsilon)}$. 

For $\beta,\beta'$, let $Y^{(\epsilon)}_{\beta,\beta'}$ be a non-zero even free field $\mathcal{F}_{\widetilde{\alpha}_0(\epsilon)}$-intertwining operator of type {\footnotesize{$\begin{pmatrix}
\ F_{{\beta+\beta'}} \\
F_{{\beta}}\ \ F_{{\beta}'}
\end{pmatrix}$}}.
Note that by the restriction, $Y^{(\epsilon)}_{\beta,\beta'}$ gives
an even $\mathfrak{ns}$-intertwining operator of type
\({\scriptsize{
\begin{pmatrix}
\ U^{(\epsilon)}(\mathfrak{ns})\ket{{\beta+\beta'}} \\
U^{(\epsilon)}(\mathfrak{ns})\ket{{\beta}}\ \ U^{(\epsilon)}(\mathfrak{ns})\ket{{\beta}'}
\end{pmatrix}
}}.
\)
From the last argument in Subsection \ref{FreeNS}, we can fix the normalization of $Y^{(\epsilon)}_{\beta,\beta'}$ as $Y^{(\epsilon)}_{\beta,\beta'}(\ket{{\beta}},z)=V_{{{\beta}}}(z)$.
Then, based on the free-field realization techniques in \cite{DF1,DF2,Felder}, we define
\begin{align}
\label{def:Rep}
\begin{aligned}
&R^{0,\pm}_{i,j}(z_1,z_2,\epsilon)\\
=&\int_{[\Box^\pm_{i,j}]}\langle{\widetilde{\beta}}_{2,2}(\epsilon){\mid}Q^{(\epsilon)}_{-}(u)Q^{(\epsilon)}_+(v)V_{{\widetilde{\beta}}_{2,2}(\epsilon)}(z_1)V_{{\widetilde{\beta}}_{2,2}(\epsilon)}(z_2){\mid}{\widetilde{\beta}}_{2,2}(\epsilon)\rangle{\rm d}u{\rm d}v\\
=&z^{\widetilde{\beta}_{2,2}(\epsilon)\widetilde{\beta}_{2,2}(\epsilon)}_1z^{\widetilde{\beta}_{2,2}(\epsilon)\widetilde{\beta}_{2,2}(\epsilon)}_2(z_1-z_2)^{\widetilde{\beta}_{2,2}(\epsilon)\widetilde{\beta}_{2,2}(\epsilon)}I^\pm_{i,j}[1](\widetilde{\alpha}_-(\epsilon)\widetilde{\beta}_{2,2}(\epsilon),\widetilde{\alpha}_-(\epsilon)^2;z_1,z_2),\\[10pt]
&R^{1,\pm}_{i,j}(z_1,z_2,\epsilon)\\
=&\int_{[\Box^\pm_{i,j}]}\langle\widetilde{\beta}_{2,2}(\epsilon){\mid}Q^{(\epsilon)}_{-}(u)Q^{(\epsilon)}_+(v)[G^{(\epsilon)}_{-\frac{1}{2}},V_{{\widetilde{\beta}}_{2,2}(\epsilon)}(z_1)]
[G^{(\epsilon)}_{-\frac{1}{2}},V_{{\widetilde{\beta}}_{2,2}(\epsilon)}(z_2)]{\mid}\widetilde{\beta}_{2,2}(\epsilon)\rangle{\rm d}u{\rm d}v\\
=&\widetilde{\beta}_{2,2}(\epsilon)\widetilde{\beta}_{2,2}(\epsilon)z^{\widetilde{\beta}_{2,2}(\epsilon)\widetilde{\beta}_{2,2}(\epsilon)}_1z^{\widetilde{\beta}_{2,2}(\epsilon)\widetilde{\beta}_{2,2}(\epsilon)}_2(z_1-z_2)^{\widetilde{\beta}_{2,2}(\epsilon)\widetilde{\beta}_{2,2}(\epsilon)-1}\\
&\ \ \ \times I^\pm_{i,j}[E](\widetilde{\alpha}_-(\epsilon)\widetilde{\beta}_{2,2}(\epsilon),\widetilde{\alpha}_-(\epsilon)^2;z_1,z_2),
\end{aligned}
\end{align}
for $i,j=1,0$ and $\epsilon\in \mathbb{D}^\times$, where $E(u,v)=\bigl((u-z_1)(u-z_2)(v-z_1)(v-z_2)\bigr)^{-1}$
and we use the operator product expansions (\ref{eq:opeVV}) and (\ref{eq:opebb}).
By Proposition \ref{thmDF0}, these functions are well-defined for $\epsilon \in \mathbb{D}^\times$. 
Then, from the following lemma, we see that each $R^{0,\pm}_{i,j}(z_1,z_2,\epsilon)$ realizes a four point correlation function with respect to $U^{(\epsilon)}(\mathfrak{ns})$.

\begin{lemma}
\label{Lem1101}
Let $i,j\in \{1,0\}$, $\epsilon\in \mathbb{D}^\times$ and $e_1,e_2\in F_{\widetilde{\beta}_{2,2}(\epsilon)}$. Then for any $A\in U^{(\epsilon)}(\mathfrak{ns})$, the function
\begin{align*}
&\int_{[\Box^\pm_{i,j}]}\langle\widetilde{\beta}_{2,2}(\epsilon){\mid}[A,Q^{(\epsilon)}_{-}(u)Q^{(\epsilon)}_+(v)]Y^{(\epsilon)}(e_1,z_1)Y^{(\epsilon)}(e_2,z_2){\mid}{\widetilde{\beta}}_{2,2}(\epsilon)\rangle{\rm d}u{\rm d}v
\end{align*}
is identically zero, where we omit the subscripts of the intertwining operators $Y^{(\epsilon)}_{\beta,\beta'}$.
\begin{proof}
By the operator expansions (\ref{eq:opeTG}) and (\ref{sc0}), we see that the function becomes a total derivative form
\begin{equation*}
\begin{split}
&z^{\widetilde{\beta}_{2,2}(\epsilon)\widetilde{\beta}_{2,2}(\epsilon)}_1z^{\widetilde{\beta}_{2,2}(\epsilon)\widetilde{\beta}_{2,2}(\epsilon)}_2(z_1-z_2)^{\widetilde{\beta}_{2,2}(\epsilon)\widetilde{\beta}_{2,2}(\epsilon)}\\
&\times\Bigl[f\int_{[\Box^\pm_{i,j}]}d_{u,v}\Bigl(U(\widetilde{\alpha}_-(\epsilon)\widetilde{\beta}_{2,2}(\epsilon),\widetilde{\alpha}_-(\epsilon)^2;u,v;z_1,z_2)\bigl(E{\rm d}u+F{\rm d}v\bigr)\Bigr)\\
&\ \ +g\int_{[\Box^\pm_{i,j}]}d_{u,v}\Bigl(U^{\frac{1}{2}}(\widetilde{\alpha}_-(\epsilon)\widetilde{\beta}_{2,2}(\epsilon),\widetilde{\alpha}_-(\epsilon)^2;u,v;z_1,z_2)\bigl(E'{\rm d}u+F'{\rm d}v\bigr)\Bigr)\Bigr]
\end{split}
\end{equation*}
for some $f,g \in \mathbb{C}[z^{\pm1}_1,z^{\pm1}_2,(z_1-z_2)^{-1}]$ and Laurent polynomials $E,F,E',F'$ in 
\(
\mathbb{C}[u^{\pm 1},v^{\pm 1},(u-z_1)^{-1},(v-z_1)^{-1},(u-z_2)^{-1},(v-z_2)^{-1},(u-v)^{-1}],
\)
where $d_{u,v}$ is the total derivative with respect to $u,v$ and we use the notation
\begin{equation*}
{U}^{\frac{1}{2}}(a,\rho;u,v;{z_1,z_2})=\mathcal{U}(\{a,a'\},\{a,a'\},\{a,a'\},2^{-1};u,v;{z_1,z_2})
\end{equation*}
as $a=\widetilde{\alpha}_-(\epsilon)\widetilde{\beta}_{2,2}(\epsilon)$ and $\rho=\widetilde{\alpha}_-(\epsilon)^2$. 
Thus, by Theorem \ref{twisted0} and by the Stokes theorem, this function must vanish.  
\end{proof}
\end{lemma}
We define
\(
\overline{R}^{0,\pm}_{i,j}(z,\epsilon)=R^{0,\pm}_{i,j}(1,z,\epsilon),
\)
\(
\overline{R}^{1,\pm}_{i,j}(z,\epsilon)=R^{1,\pm}_{i,j}(1,z,\epsilon).
\)
From Lemma \ref{Lem1101} and the $L^{(\epsilon)}_0$-conjugation formula for intertwining operators, we have
\begin{equation}
\label{Lo}
\begin{split}
R^{0,\pm}_{i,j}(z_1,z_2,\epsilon)&=z^{-2h_{\widetilde{\beta}_{2,2}(\epsilon)}}_1\overline{R}^{0,\pm}_{i,j}(z_2/z_1,\epsilon),\\
R^{1,\pm}_{i,j}(z_1,z_2,\epsilon)&=z^{-2h_{\widetilde{\beta}_{2,2}(\epsilon)}-1}_1\overline{R}^{1,\pm}_{i,j}(z_2/z_1,\epsilon),
\end{split}
\end{equation}
where $h_{\widetilde{\beta}_{2,2}(\epsilon)}=2^{-1}\widetilde{\beta}_{2,2}(\epsilon)(\widetilde{\beta}_{2,2}(\epsilon)-\widetilde{\alpha}_0(\epsilon))$ (see (\ref{hbeta})). 
Then by Proposition \ref{lemDf} and by the definition (\ref{def:Rep}), we see that $\int_{C_0}\overline{R}^{0,\pm}_{i,j}(z,\epsilon){\rm d}\epsilon$ is nonzero for any simple closed curve $C_{0}\subset \mathbb{D}^\times(\epsilon)$ arround $\epsilon=0$.
Then
we define
\begin{align}
\label{Def1011}
\Psi^\pm_{i,j}(z):=\int_{C_0}\overline{R}^{0,\pm}_{i,j}(z,\epsilon){\rm d}\epsilon
\end{align}
for $i,j=1,0$. 
Note that these functions are invariant under the homotopy deformation of $C_0$ in $\mathbb{D}^\times(\epsilon)$.

\begin{proposition}
\label{solutionF}
The functions $\Psi^\pm_{i,j}(z)$ $(i,j=1,0)$ satisfy the Fuchsian differential equation (\ref{fuchsian}).
\begin{proof}
From the structure of the Fock module \(F_{2,2}\) (see Proposition \ref{socleFockmodule}), we see that the lowest weight vector $\ket{\widetilde{\beta}_{2,2}(0)}=\ket{{\beta}_{2,2}}\in F_{2,2}$ satisfies the relation $S^0_{2,2}\ket{{\beta}_{2,2}}=0$, where $S^0_{2,2}$ is defined by (\ref{sing0318}) with the central charge fixed to $c_{1,2m+1}$. We define an $\epsilon$-deformation of \(S^0_{2,2}\) as follows
\begin{equation*}
\begin{split}
S_{2,2}(\epsilon)=s_1(\epsilon)(L^{(\epsilon)}_{-1})^2+s_2(\epsilon)L^{(\epsilon)}_{-2}-s_1(\epsilon)G^{(\epsilon)}_{-\frac{3}{2}}G^{(\epsilon)}_{-\frac{1}{2}}\in U^{(\epsilon)}(\mathfrak{ns}),
\end{split}
\end{equation*}
where
\begin{equation*}
s_1(\epsilon)=\frac{4\widetilde{\alpha}_+(\epsilon)^2}{1-\widetilde{\alpha}_+(\epsilon)^4},\ \ \ \ \ \ \ \ \ \ s_2(\epsilon)=-\frac{2-2\widetilde{\alpha}_+(\epsilon)^2}{\widetilde{\alpha}_+(\epsilon)^2+1}.
\end{equation*}
From the definition, we see that 
\begin{equation}
\label{S220423}
S_{2,2}(0)=S^0_{2,2},\ \ \ \ s_1(\epsilon),s_2(\epsilon)\in \mathcal{O}(\mathbb{D}(\epsilon))\setminus\epsilon \mathcal{O}(\mathbb{D}(\epsilon)),
\end{equation}
where $\epsilon \mathcal{O}(\mathbb{D}(\epsilon))$ is the subring of $\mathcal{O}(\mathbb{D}(\epsilon))$ defined by
\begin{equation*}
\Bigl\{\left.f(\epsilon)\in \mathcal{O}(\mathbb{D}(\epsilon))\ \right|\ \frac{f(\epsilon)}{\epsilon}\in \mathcal{O}(\mathbb{D}(\epsilon))\Bigr\}.
\end{equation*}
By a straightforward calculation, we see that the lowest weight vector $\ket{\widetilde{\beta}_{2,2}(\epsilon)}$ satisfies the relation $S_{2,2}(\epsilon)\ket{\widetilde{\beta}_{2,2}(\epsilon)}=0$ (cf. \cite{BA,IK2}).
Then, by Lemma \ref{Lem1101}, similar to (\ref{kore20231014}), we can show that $R^{0,\pm}_{i,j}(z_1,z_2,\epsilon)$ and $R^{1,\pm}_{i,j}(z_1,z_2,\epsilon)$ satisfy the following differential equations
\begin{equation}
\label{kore20231015}
\begin{split}
&\Biggl\{\bigl(-2s_1(\epsilon)h_{\widetilde{\beta}_{2,2}(\epsilon)}+s_2(\epsilon)h_{\widetilde{\beta}_{2,2}(\epsilon)}+\frac{s_2(\epsilon)}{2}\bigr)\Bigl(\frac{1}{z^2_1}+\frac{1}{z^2_2}\Bigr)+s_1(\epsilon)(\partial_{z_1}+\partial_{z_2})^2\\
&\ \ \ \ \ \ \ \ \ \ \ \ \ \ \ \ \ \ \ \ \ \ \ \ \ \ \ \ \ \ \ \ \ \ \ \ +(s_1(\epsilon)-s_2(\epsilon))\Bigl(\frac{\partial_{z_1}}{z_1}+\frac{\partial_{z_2}}{z_2}\Bigr)\Biggr\}R^{1,\pm}_{i,j}(z_1,z_2,\epsilon)\\
&+\Biggl\{-s_1(\epsilon)\Bigl(\frac{1}{z_2}-\frac{1}{z_1}\Bigr)\partial_{z_1}\partial_{z_2}+2s_1(\epsilon)h_{\widetilde{\beta}_{2,2}(\epsilon)}\Bigl(\frac{\partial_{z_1}}{z^2_2}-\frac{\partial_{z_2}}{z^2_1}\Bigr)\Biggr\}R^{0,\pm}_{i,j}(z_1,z_2,\epsilon)=0,\\[10pt]
&\Biggl\{s_2(\epsilon)h_{\widetilde{\beta}_{2,2}(\epsilon)}\Bigl(\frac{1}{z^2_1}+\frac{1}{z^2_2}\Bigr)+s_1(\epsilon)(\partial_{z_1}+\partial_{z_2})^2\Biggr\}R^{0,\pm}_{i,j}(z_1,z_2,\epsilon)\\
&+(s_1(\epsilon)-s_2(\epsilon))\Bigl(\frac{\partial_{z_1}}{z_1}+\frac{\partial_{z_2}}{z_2}\Bigr)R^{0,\pm}_{i,j}(z_1,z_2,\epsilon)-s_1(\epsilon)\Bigl(\frac{1}{z_1}-\frac{1}{z_2}\Bigr)R^{1,\pm}_{i,j}(z_1,z_2,\epsilon)=0.
\end{split}
\end{equation}
Thus by (\ref{Lo}) and (\ref{kore20231015})$, \overline{R}^{0,\pm}_{i,j}(z,\epsilon)$ satisfy a fourth-order differential equation. Let $\mathcal{D}^{(\epsilon)}_z$ be the fourth-order differential operator defined by 
\begin{equation*}
\mathcal{D}^{(\epsilon)}_z\overline{R}^{0,\pm}_{i,j}(z,\epsilon)=0,\ \ \ \ \mathcal{D}^{(\epsilon)}_z=\partial^4_{z}+({\rm lower}\ {\rm order}\ {\rm terms}\ {\rm for}\ \partial^k_{z}),
\end{equation*} 
and 
let $\mathcal{D}_z$ be the fourth-order differential operator of (\ref{fuchsian}).
Then, by (\ref{S220423}), we can see that
\begin{equation}
\label{eq:0-edif}
\mathcal{D}^{(\epsilon)}_z-\mathcal{D}_z\in \epsilon\mathcal{O}(\mathbb{D}(\epsilon))[z,z^{-1},(z-1)^{-1}]\langle\partial_{z}\rangle,
\end{equation}
where the right hand side is the ring of linear differential operators with coefficients in the ring $\epsilon\mathcal{O}(\mathbb{D}(\epsilon))[z,z^{-1},(z-1)^{-1}]:=\epsilon\mathcal{O}(\mathbb{D}(\epsilon))\otimes\mathbb{C}[z,z^{-1},(z-1)^{-1}]$.
Thus by Proposition \ref{lemDf}, we have
\begin{align*}
0&=\int_{C_0}\mathcal{D}^{(\epsilon)}_z\overline{R}^{0,\pm}_{i,j}(z,\epsilon){\rm d}\epsilon\\
&=\int_{C_0}(\mathcal{D}^{(\epsilon)}_z-\mathcal{D}_z)\overline{R}^{0,\pm}_{i,j}(z,\epsilon){\rm d}\epsilon
+\int_{C_0}\mathcal{D}_z\overline{R}^{0,\pm}_{i,j}(z,\epsilon){\rm d}\epsilon\\
&=\mathcal{D}_z\int_{C_0}\overline{R}^{0,\pm}_{i,j}(z,\epsilon){\rm d}\epsilon.
\end{align*}
Therefore, $\Psi^\pm_{i,j}$ satisfy the Fuchsian differential equation (\ref{fuchsian}).
\end{proof}
\end{proposition}
By noting the Riemann scheme (\ref{Scheme}) of (\ref{fuchsian}), from Propositions \ref{lemDf} and \ref{solutionF}, we obtain the following proposition.
\begin{proposition}
\label{propfund}
The sets $\{\Psi^+_{i,j}(z)\ |\ i,j=1,0\}$ and $\{\Psi^-_{i,j}(z)\ |\ i,j=1,0\}$ are fumdamental systems of solutions of (\ref{fuchsian}) at $z=0$ and $z=1$, respectively.
Let $\rho^+_{i,j}$ and $\rho^-_{i,j}$ be the characteristic exponents of $\Psi^+_{i,j}(z)$, $\Psi^-_{i,j}(z)$ at $z=0$ and $z=1$, respectively. Then we have $\rho^{\pm}_{i,j}=\rho_{i,j}$ for all $i,j=1,0$, where $\rho_{i,j}$ are defined by (\ref{charexp}).
\end{proposition}
\begin{remark}
The important parts of the above construction of solutions are Lemma \ref{Lem1101}, Proposition \ref{lemDf} and (\ref{eq:0-edif}). The explicit forms of $\mathcal{D}_z$ and $\mathcal{D}^{(\epsilon)}_z$ are not important.
Hence, we believe that our construction is applicable to the correlation functions of other logarithmic minimal models.
\end{remark}
From Lemma \ref{fourrel}, we obtain the following proposition.
\begin{proposition}
\label{connection}
We have the the following connection formulas:
\begin{equation*}
\begin{pmatrix}
\Psi^+_{1,1}(z)\\
\Psi^+_{0,1}(z)\\
\Psi^+_{1,0}(z)\\
\Psi^+_{0,0}(z) 
\end{pmatrix}
=\frac{(-1)^m}{2}
\begin{pmatrix}
 c^{-1}& -c^{-1}& -c^{-1}& c^{-1}\\
 -c+c^{-1}& -c^{-1}&c-c^{-1}&c^{-1}\\
 -3c^{-1}&3c^{-1}&-c^{-1}&c^{-1}\\
 3(c-c^{-1})&3c^{-1}&c-c^{-1}&c^{-1}
\end{pmatrix}
\begin{pmatrix}
\Psi^-_{1,1}(z) \\
\Psi^-_{0,1}(z)\\
\Psi^-_{1,0}(z)\\
\Psi^-_{0,0}(z) 
\end{pmatrix}
,
\end{equation*}
where $c=2{\rm cos}(\pi \alpha_-\beta_{2,2})$. 
\end{proposition}
\begin{remark}
From Proposition \ref{connection}, we obtain
\begin{equation}
\label{mono1}
\begin{pmatrix}
\Psi^+_{1,1}(z)+\Psi^+_{1,0}(z)  \\
\Psi^+_{0,1}(z)+\Psi^+_{0,0}(z) \\
\end{pmatrix}
=(-1)^m\begin{pmatrix}
 -c^{-1} & c^{-1} \\
 c-c^{-1}& c^{-1}
\end{pmatrix}
\begin{pmatrix}
\Psi^-_{1,1}(z)+\Psi^-_{1,0}(z)  \\
\Psi^-_{0,1}(z)+\Psi^-_{0,0}(z) \\
\end{pmatrix}
.
\end{equation}
Then, from (\ref{mono1}), it seems that the monodoromy of the Fuchsian differential equation (\ref{fuchsian}) is reducible and $\{\Psi^\pm_{1,1}(z)+\Psi^\pm_{1,0}(z),
\Psi^\pm_{0,1}(z)+\Psi^\pm_{0,0}(z)\}$ gives a two dimensional subspace of the monodromy representation of (\ref{fuchsian}).
\end{remark}

\section{Tensor structure on $\mathcal{SW}(m)$-{\rm mod}}
\label{NSfusion}
Since the super triplet $W$-algebra $\mathcal{SW}(m)$ is $C_2$-cofinite, Proposition 2.1 in \cite{CGNS} (see also \cite[Theorem 2.25]{CMOY} and \cite[Theorem 4.13]{H}) show that $\mathcal{SW}(m)$ has braided tensor supercategory structure developed in the papers \cite{CKM} and \cite{HLZ1}-\cite{HLZ8}. We denote by ($\mathcal{SW}(m)\mathchar`-{\rm mod},\boxtimes$) the tensor supercategory on $\mathcal{SW}(m)\mathchar`-{\rm mod}$, where the unit object is given by $X_1$ and the symbol $\boxtimes$ denotes the tensor product. 
It is known that the tensor product $\boxtimes$ of ($\mathcal{SW}(m)\mathchar`-{\rm mod},\boxtimes$) is right exact \cite[Proposition 2.1]{CGNS}.
In this section, we study the fusion structure of $\mathcal{SW}(m)$ and determine the structure of the projective covers of all simple $\mathcal{SW}(m)$-modules. See \cite{CRR,CMOY} for the detailed structure of fusion rules and the tensor category for the $N=1$ super Virasoro minimal models.

\subsection{Tensor product $\boxtimes$ and $P(w)$-intertwining operators}
In this subsection, we review the definition of the tensor product $\boxtimes$ and $P(w)$-intertwining operators in accordance with \cite{CKM,HLZ3,Kanade} and derive some identities known as the Nahm-Gaberdiel-Kausch fusion algorithm(cf. \cite{CRR,GK1,Na}). 
\begin{definition}
Let $V$ be a $\frac{1}{2}\mathbb{Z}_{\geq 0}$-graded vertex operator superalgebra and let $\mathcal{C}$ be a category of grading-restricted generalized $V$-modules. A {\rm tensor product} $($or {\rm fusion product}$)$ of $M_1$ and $M_2$ in $\mathcal{C}$ is a pair $(M_1\boxtimes M_2,\mathcal{Y}_{\boxtimes})$, with $M_1\boxtimes M_2$ and $\mathcal{Y}_{\boxtimes}$ an intertwining operator of type 
{\scriptsize{$\begin{pmatrix}
   M_1\boxtimes M_2  \\
   M_1\ M_2
\end{pmatrix}
$}}, 
which satisfies the following universal property: For any $M_3\in\mathcal{C}$ and intertwining operator $\mathcal{Y}$ of type 
{\scriptsize{$\begin{pmatrix}
   M_3  \\
   M_1\ M_2
\end{pmatrix}
$}}, there is a unique $V$-module homomorphism $f:M_1\boxtimes M_2\rightarrow M_3$ such that $\mathcal{Y}=f\circ \mathcal{Y}_{\boxtimes}$.
\end{definition}

In the paper \cite{HLZ3}, the notion of $P(w)$-intertwining operators and the $P(w)$-tensor product are introduced. The definitions are as follows. 
\begin{definition}
\label{P(w)}
Fix $w\in \mathbb{C}^\times$. Let $V$ be a $\frac{1}{2}\mathbb{Z}_{\geq 0}$-graded vertex operator superalgebra and let $\mathcal{C}$ be a category of grading-restricted generalized $V$-modules. Given $M_1$, $M_2$ and $M_3$ in $\mathcal{C}$, a {\rm parity}-{\rm homogeneous} $P(w)$-{\rm intertwining operator} $I$ {\rm of type}
{\scriptsize{$\begin{pmatrix}
   \ M_3  \\
   M_1\ M_2
\end{pmatrix}
$}}
is a parity-homogeneous bilinear map $I$ $:$ $M_1\otimes M_2\rightarrow \overline{M}_3$ that satisfies the following properties:
\begin{enumerate}
\item For any $\psi_1\in M_1$ and $\psi_2\in M_2$, $\pi_h(I[\psi_1\otimes \psi_2])=0$ for all $h\ll0$, where $\pi_h$ denotes the projection onto the generalised eigenspace $M_3[h]$ of $L_0$-eigenvalue $h$.
\item For any $\psi_1\in M_1$, $\psi_2\in M_2$, $\psi^*_3\in M^*_3$ and $v\in V$, the three point functions
\begin{align*}
&\langle \psi^*_3, Y_3(v,z)I[\psi_1\otimes \psi_2]\rangle,&
\langle \psi^*_3, I[Y_1(v,z-w)\psi_1\otimes \psi_2]\rangle&,\\  
&\langle \psi^*_3, I[\psi_1\otimes Y_2(v,z)\psi_2]\rangle
\end{align*}
are absolutely convergent in the regions ${\mid}z{\mid}>{\mid}w{\mid}>0$, ${\mid}w{\mid}>{\mid}z-w{\mid}>0$, ${\mid}w{\mid}>{\mid}z{\mid}>0$, respectively, where $Y_i$ is the action of $V$-module on $M_i$.
\item Given any $f(t)\in R_{P(w)}:=\mathbb{C}[t,t^{-1},(t-w)^{-1}]$ and parity-homogeneous vectors $v\in V$, $\psi_1\in M_1$, $\psi_2\in M_2$, $\psi^*_3\in M^*_3$, we have the following identity
\begin{equation}
\label{P-comp}
\begin{aligned}
&(-1)^{|v||I|}\oint_{0,w}f(z)\langle \psi^*_3, Y_3(v,z)I[\psi_1\otimes \psi_2]\rangle\frac{{\rm d}z}{2\pi i}\\
&=\oint_{w}f(z)\langle\psi^*_3, I[Y_1(v,z-w)\psi_1\otimes \psi_2]\rangle\frac{{\rm d}z}{2\pi i}\\
&\ \ \ \ \ +(-1)^{|v||\psi_1|}\oint_{0}f(z)\langle \psi^*_3, I[\psi_1\otimes Y_2(v,z)\psi_2]\rangle\frac{{\rm d}z}{2\pi i}.
\end{aligned}
\end{equation}
\end{enumerate}
A general $P(w)$-intertwining operator of type {\scriptsize{$\begin{pmatrix}
   M_3  \\
   M_1\ M_2
\end{pmatrix}
$}} is a sum of parity-homogeneous ones.
\end{definition}

\begin{definition}
Let $V$ be a $\frac{1}{2}\mathbb{Z}_{\geq 0}$-graded vertex operator superalgebra and let $\mathcal{C}$ be a category of grading-restricted generalized $V$-modules. A $P(w)$-{\rm tensor product} of $M_1$ and $M_2$ in $\mathcal{C}$ is a pair $(M_1\boxtimes_{P(w)} M_2,{\boxtimes_{P(w)}})$, with $M_1\boxtimes_{P(w)} M_2$ and ${\boxtimes_{P(w)}}$ a $P(w)$-intertwining operator of type 
{\scriptsize{$\begin{pmatrix}
   M_1\boxtimes_{P(w)} M_2  \\
   M_1\ \ M_2
\end{pmatrix}
$}}, 
which satisfies the following universal property: For any $M_3\in\mathcal{C}$ and $P(w)$-intertwining operator $I$ of type 
{\scriptsize{$\begin{pmatrix}
   M_3  \\
   M_1\ M_2
\end{pmatrix}
$}}, there is a unique $V$-module homomorphism $\eta:M_1\boxtimes_{P(w)} M_2\rightarrow M_3$ such that 
\begin{align*}
\overline{\eta}\circ {\boxtimes_{P(w)}}[\psi_1\otimes \psi_2]=I[\psi_1\otimes \psi_2]
\end{align*}
for all $\psi_1\in M_1$ and $\psi_2\in M_2$, where $\overline{\eta}$ denotes the extension of $\eta$ to a map between the completions of $M_1\boxtimes_{P(w)} M_2$ and $M_3$.
\end{definition}
It is known that the definition $P(w)$-tensor product $\boxtimes_{P(w)}$ does not depend on the choice of $w\in \mathbb{C}^\times$. Precisely, the following proposition holds \cite[Corollary 3.36]{CKM} (see also \cite[Remark 4.22]{HLZ3}).
\begin{proposition}[\cite{CKM}]
\label{notw}
Let $V$ be a $\frac{1}{2}\mathbb{Z}_{\geq 0}$-graded vertex operator superalgebra and let $M_1,M_2$ be $V$-modules. Suppose that for some $w_0\in \mathbb{C}^\times$, $P(w_0)$-tensor product $M_1\boxtimes_{P(w_0)} M_2$ exists. Then, for any $w,w'\in \mathbb{C}^\times$, $M_1\boxtimes_{P(w)} M_2$ and $M_1\boxtimes_{P(w')} M_2$ are isomorphic.
\end{proposition}
It is known that certain specializations for intertwining operators yield $P(w)$-intertwining operators and by these specializations give linear isomorphisms from the the spaces of intertwining operators to the spaces of $P(w)$-intertwining operators of the same types \cite[Proposition 3.15]{CKM}(see also \cite{HLZ3}). Thus, from Proposition \ref{notw}, $\mathcal{SW}(m)$ has a braided tensor supercategory structure with $\boxtimes=\boxtimes_{P(1)}$.
\vspace{2mm}

In the following, we will introduce some useful formulas derived from the $P(w)$-compatibility conditions.
We define a translation map
\begin{eqnarray*}
T_1:\mathbb{C}(t)\rightarrow \mathbb{C}(t),\ \ \ \ \ \ \ {\rm by}\ \ \ \ \ \ \ f(t)\mapsto f(t+1),
\end{eqnarray*}
and a expansion map 
\begin{align*}
\iota_+:\mathbb{C}(t)\hookrightarrow \mathbb{C}((t))
\end{align*}
that expands a given rational function in $t$ as a power series around $t=0$.
Given a $\frac{1}{2}\mathbb{Z}_{\geq 0}$-graded vertex operator superalgebra $V$ and a homogeneous vector $v\in V[h]$, we use notation $v_{n-h+1}=vt^n=v\otimes t^n\in V\otimes \mathbb{C}[t,t^{-1}]$ for $n\in \mathbb{Z}$.
Given $V$-modules $M_1$, $M_2$, $M_3$ and a parity-homogeneous $P(1)$-intertwining operator $I$ of type
{\scriptsize{$\begin{pmatrix}
   \ M_3  \\
   M_1\ M_2
\end{pmatrix}
$}}, as detailed in \cite{Kanade}, by using the action (\ref{cont-act}) and the identities (\ref{P-comp}), we can define the action of $V\otimes \mathbb{C}[t,t^{-1},(t-1)^{-1}]$ or $V\otimes \mathbb{C}((t))$ on $M^*_3$ as 
\begin{equation}
\begin{aligned}
(-1)^{|v||I|}\langle vf(t)\psi^*_3, I[\psi_1\otimes\psi_1]\rangle&=(-1)^{|v||I|}\langle v\iota_+(f(t))\psi^*_3, I[\psi_1\otimes\psi_1]\rangle\nonumber\\
&=\langle \psi^*_3,I[\iota_+\circ T_1\bigl(v^{{\rm opp}}f(t^{-1})\bigl)\psi_1\otimes\psi_2]\rangle\\
&\ \ \ +(-1)^{|v||\psi_1|}\langle\psi^*_3,I[\psi_1\otimes\iota_+\bigl(v^{{\rm opp}}f(t^{-1})\bigl)\psi_2]\rangle  
\end{aligned}
\label{NGK0}
\end{equation}
where $\psi^*_3\in M^*_3$, $\psi_i\in M_i$, and
\begin{align*}
v^{{\rm opp}}:=e^{t^{-1}L_{1}}(-t^2)^{L_0}vt^{-2}.
\end{align*}

\begin{lemma}
\label{NGK01}
Let $V$ be a $\frac{1}{2}\mathbb{Z}_{\geq 0}$-graded vertex operator superalgebra and let $v\in V$ be a non-zero Virasoro primary vector with the conformal weight $h$:
\begin{align*} 
&L_0v=hv,
&L_nv=0\ \ \ (n\geq 1).
\end{align*}
Given $V$-modules $M_1$, $M_2$, $M_3$ and a parity-homogeneous $P(1)$-intertwining operator $I$ {\rm of type}
{\scriptsize{$\begin{pmatrix}
   \ M_3  \\
   M_1\ M_2
\end{pmatrix}
$}}, we have the following identities:
\begin{align*}
&(-1)^{|v||I|}\langle v_n\psi^*_3, I[\psi_1\otimes\psi_2]\rangle\\
&=\sum_{i=0}^\infty\binom{h-n-1}{i}\langle\psi^*_3,I[\bigl(v_{i-h+1}\psi_1\bigr)\otimes\psi_2]\rangle+(-1)^{|v||\psi_1|}\langle\psi^*_3,I[\psi_1\otimes\bigl(v_{-n}\psi_2\bigr)]\rangle,\\[10pt]
&(-1)^{|v||I|}\sum_{i=0}^\infty\binom{n+h-1}{i}(-1)^i\langle v_{i-n}\psi^*_3,I[\psi_1\otimes\psi_2]\rangle\\
&=\langle\psi^*_3,I[\bigl(v_{n}\psi_1\bigr)\otimes\psi_2]\rangle+(-1)^{|v||\psi_1|}\sum_{i=0}^\infty\binom{n+h-1}{i}(-1)^{i}\langle\psi^*_3,I[\psi_1\otimes\bigl(v_{n-i}\psi_2\bigr)]\rangle,\\[10pt]
&(-1)^{|v||I|}\sum_{i=0}^\infty\binom{n+h-2}{i}(-1)^i\langle v_{i-n}\psi^*_3,I[\psi_1\otimes\psi_2]\rangle\\
&=\langle\psi^*_3,I[\bigl((v_{n-1}+v_{n})\psi_1\bigr)\otimes\psi_2]\rangle\\
&\ \ \ \ \ \ \ \ \ +(-1)^{|v||\psi_1|}\sum_{i=0}^\infty\binom{n+h-2}{i}(-1)^{n-i+h-2}\langle\psi^*_3,I[\psi_1\otimes\bigl(v_{i-h+2}\psi_2\bigr)]\rangle,
\end{align*}
where $\psi^*_3\in M^*_3$, $\psi_i\in M_i$ and $n\in \mathbb{Z}$.
\begin{proof}
For the first identity, let $f(t)=t^{n+h-1}$ in (\ref{NGK0}), for the second, let $f(t)=t^{2h-2}(t^{-1}-1)^{n+h-1}$ in (\ref{NGK0}), and for the third identity, let $f(t)=t^{2h-3}(t^{-1}-1)^{n+h-2}$ in (\ref{NGK0}).
\end{proof}
\end{lemma}

From the next subsection,  we use the shorthand notation
\begin{eqnarray}
\label{abr0204}
 \langle \psi^*_3,\psi_1\otimes\psi_2\rangle= \langle \psi^*_3,I[\psi_1\otimes\psi_2]\rangle
\end{eqnarray}
for $P(1)$-intertwining operators $I$.


\subsection{Self duality of the simple module $X_2$}
\label{X''dual}
Note that all simple $\mathcal{V}_L$-modules can be written as forms
\begin{eqnarray*}
\mathcal{V}_{L+\beta}:=\bigoplus_{n\in\mathbb{Z}}F_{\beta+n\alpha_+},\ \ \beta\in \{\ \beta_{r,s;n}\ {\mid}\ r,s,n\in\mathbb{Z} \}
\end{eqnarray*}
(for the definition of the simple $\mathcal{V}_L$-modules, see (\ref{simpleLattice})).
Given simple $\mathcal{V}_L$-modules $\mathcal{V}_{L+\beta}$, $\mathcal{V}_{L+\beta'}$ and $\mathcal{V}_{L+\beta''}$, it can be proved easily that there are no $\mathcal{V}_{L}$-intertwining operators of type 
{\footnotesize{$\begin{pmatrix}
\ \mathcal{V}_{L+\beta''} \\
\mathcal{V}_{L+\beta'}\ \ \mathcal{V}_{L+\beta}
\end{pmatrix}$}}
unless $\beta''\equiv\beta'+\beta\ {\rm mod}\ L$, and ${\rm dim}_{\mathbb{C}}I^{0}_{\mathcal{V}_{L}}${\footnotesize{$\begin{pmatrix}
\ \mathcal{V}_{L+\beta'+\beta} \\
\mathcal{V}_{L+\beta'}\ \ \mathcal{V}_{L+\beta}
\end{pmatrix}$}}$=1$.
Let $\mathcal{Y}$ be the even $\mathcal{V}_{L}$-intertwining operator of type {\footnotesize{$\begin{pmatrix}
\ \mathcal{V}_{L+\beta'+\beta} \\
\mathcal{V}_{L+\beta'}\ \ \mathcal{V}_{L+\beta}
\end{pmatrix}$}}. Then, by restricting the action of $\mathcal{V}_{L}$ to $\mathcal{SW}(m)$, $\mathcal{Y}$ defines an $\mathcal{SW}(m)$-intertwining operator of type {\footnotesize{$\begin{pmatrix}
\ \mathcal{V}_{L+\beta'+\beta} \\
\mathcal{V}_{L+\beta'}\ \ \mathcal{V}_{L+\beta}
\end{pmatrix}$}}.
We denote by $\mathcal{Y}_{\beta',\beta}$ this even $\mathcal{SW}(m)$-intertwining operator.
\begin{lemma}
\label{Intertwining}
For $2\leq s\leq 2m$, we have
\begin{equation*}
I_{\mathcal{SW}(m)}
\begin{pmatrix}
\ X_{s-1} \\
X_2\ \ X_s
\end{pmatrix}
\neq 0
,\ \ \ \ \ \ \ \ 
I_{\mathcal{SW}(m)}
\begin{pmatrix}
\ X_{s+1} \\
X_2\ \ X_s
\end{pmatrix}
\neq 0
.
\end{equation*}
\begin{proof}
We only prove 
\begin{equation*}
I_{\mathcal{SW}(m)}
\begin{pmatrix}
\ X_{2i-1} \\
X_2\ \ X_{2i}
\end{pmatrix}
\neq 0
.
\end{equation*}
The other cases can be shown in similar ways.

Let us consider the $\mathcal{SW}(m)$-intertwining operator $\mathcal{Y}=\mathcal{Y}_{\beta_1,\beta_2}$, where $\beta_1=\beta_{2,2}$ and $\beta_2=\beta_{1,2(m-i)+1;2}$. 
By (\ref{0202V}), we have
\begin{align}
\label{eqde}
\langle\beta_{1,2(m-i+1);1}{\mid}\mathcal{Y}({\mid}\beta_{2,2}\rangle,z){\mid}\beta_{1,2(m-i)+1;2}\rangle\neq 0.
\end{align}
From Proposition \ref{socleV}, we have
\begin{align*}
&X_2\simeq \mathcal{SW}(m) {\mid}\beta_{2,2}\rangle,&
&\mathcal{V}_{L+\beta_2}\simeq \mathcal{SW}(m) {\mid}\beta_{1,2(m-i)+1;2}\rangle,\\
&\mathcal{V}_{L+\beta_2+\beta_1}\simeq \mathcal{SW}(m) {\mid}\beta_{1,2(m-i+1);1}\rangle,
\end{align*}
and the exact sequence
\begin{align*}
0\rightarrow X_{2(m-i+1)}\rightarrow \mathcal{V}_{L+\beta_2+\beta_1}\rightarrow X_{2i-1}\rightarrow 0.
\end{align*}
Thus, by (\ref{eqde}), we have
\begin{equation}
I_{\mathcal{SW}(m)}
\begin{pmatrix}
\ X_{2i-1} \\
X_2\ \mathcal{V}_{L+\beta_2}
\end{pmatrix}
\neq 0
.
\label{202208130}
\end{equation}
From Proposition \ref{socleV}, we have the exact sequence
\begin{align*}
0\rightarrow X_{2(m-i)+1}\rightarrow \mathcal{V}_{L+\beta_2}\rightarrow X_{2i}\rightarrow 0.
\end{align*}
Then we have the following exact sequence
\begin{align*}
X_2\boxtimes X_{2(m-i)+1}\rightarrow X_2\boxtimes \mathcal{V}_{L+\beta_2}\rightarrow X_2\boxtimes X_{2i}\rightarrow 0.
\end{align*}
From this exact sequence, we have the following exact sequence
\begin{align}
0&\rightarrow{\rm Hom}_{\mathcal{SW}(m)}(X_2\boxtimes X_{2i},X_{2i-1})\rightarrow {\rm Hom}_{\mathcal{SW}(m)}(X_2\boxtimes \mathcal{V}_{L+\beta_2},X_{2i-1})\nonumber\\
&\rightarrow {\rm Hom}_{\mathcal{SW}(m)}(X_2\boxtimes X_{2(m-i)+1},X_{2i-1}).
\label{12140}
\end{align}
By Lemma \ref{NGKlem}, we have
$
{\rm Hom}_{\mathcal{SW}(m)}(X_2\boxtimes X_{2(m-i)+1},X_{2i-1})=0.
$
Therefore by (\ref{202208130}) and (\ref{12140}), we obatin 
$
{\rm Hom}_{\mathcal{SW}(m)}(X_2\boxtimes X_{2i},X_{2i-1})\neq 0.
$
\end{proof}
\end{lemma}

From this subsection, we will use the following notation.
\begin{definition}
For any $M\in \mathcal{SW}(m)\mathchar`-{\rm mod}$, we define the following vector space 
\begin{align*}
A_0(M)=\bigl\{\psi\in M\setminus\{0\}\ {\mid}\ v_n\psi=0\ {\rm for}\ v\in \mathcal{SW}(m)\ {\rm and}\ n\in \frac{1}{2}\mathbb{Z}_{>0}\bigr\}.
\end{align*}
\end{definition}

\begin{lemma}
\label{NGKlem}
For $i=1,\dots,m$, the vector space $A_0((X_2\boxtimes X_{2i})^*)$ is at most two dimensional. $L_0$ acts semisimply on $A_0((X_2\boxtimes X_{2i})^*)$ and any $L_0$ eigenvalue of this space is contained in $\{h_{1,2i-1},h_{1,2i+1}\}$, where $h_{1,2i-1}$ and $h_{1,2i+1}$ are the minimal conformal weights of $X_{2i-1}$ and $X_{2i+1}$, respectively. 
\end{lemma}
Similar results are obtained for the $N=1$ super Virasoro algebra by using Zhu bimodules \cite[Subsection 4.1]{CMOY}.
\begin{proof}
By Lemma \ref{Intertwining}, we see that the tensor product $X_2\boxtimes X_{2i}$ is non-zero. Let $\psi^*$ be an arbitrary non zero $L_0$ homogeneous vector of $A_0((X_2\boxtimes X_{2i})^*)$.
Let $\phi_1$ and $\phi_2$ be arbitrary non zero $L_0$ homogeneous vectors $X_2$ and $X_{2i}$ such that $\langle\psi^*,\phi_1\otimes \phi_2\rangle\neq 0$, where we use the shorthand notation (\ref{abr0204}) as $I={\boxtimes_{P(1)}}$. For $1\leq j\leq m$, let $\{v^+_{j},v^-_{j}\}$ be a basis of the minimal conformal weight space of $X_{2j}$ such that
\begin{align*}
&\widehat{W}^\pm_0v^\pm_{j}=0,\ \ \ \ \ \ \ \ \ 
&\widehat{W}^\pm_0v^\mp_{j}\in\mathbb{C}^\times v^\pm_{j},
\end{align*}
where $\widehat{W}^\pm_0$ are the zero-mode of the fields $Y(\widehat{W}^\pm,z)$.
For $n\geq 1$, let 
\begin{align*}
v^{(n)}_{\frac{2k-1}{2}}\ \ \ \ \ (k=-n,\dots,n+1)
\end{align*}
be the minimal conformal weight vectors of the subspace $(2n+2)L(h_{2n+2,2i})\subset X_{2i}$ defined in Proposition \ref{dfn0830}. 
First let us show
\begin{align}
\label{kieru}
\langle\psi^*, U(\mathfrak{ns})v^\pm_1\otimes v^{(n)}_{\frac{2k-1}{2}}\rangle=0\ \ (n\geq 1,\ k\in \mathbb{Z}).
\end{align}
Note that $v^\pm_1$ satisfy the following relations (cf.\cite{BA,CMOY}) 
\begin{equation}
\label{sing1104}
\Bigl\{  \frac{4t}{t^2-1}G^4_{-\frac{1}{2}}+\frac{t+1}{t-1}G_{-\frac{1}{2}}G_{-\frac{3}{2}}+\frac{t-1}{t+1}G_{-\frac{3}{2}}G_{-\frac{1}{2}} \Bigr\}v^\pm_1=0\ \ \ (t=-2m-1).
\end{equation}
By using Lemma \ref{NGK01} and the relation (\ref{sing1104}),
we see that, depending on whether $\psi^*$ is even or odd, the values $\langle\psi^*,U(\mathfrak{ns})v^\pm_1\otimes \otimes v^{(n)}_{\frac{2k-1}{2}}\rangle$ is determined by the numbers
\begin{align*}
\langle\psi^*,v^\delta_1\otimes v^{(2n)}_{\frac{2l-1}{2}}\rangle,&&\langle\psi^*,L_{-1}v^\delta_1\otimes v^{(2n)}_{\frac{2l-1}{2}}\rangle
\end{align*}
or 
\begin{align*}
\langle\psi^*,G_{-\frac{1}{2}}v^\delta_1\otimes v^{(2n+1)}_{\frac{2l-1}{2}}\rangle,&&\langle\psi^*,G^3_{-\frac{1}{2}}v^\delta_1\otimes v^{(2n+1)}_{\frac{2l-1}{2}}\rangle,
\end{align*}
for $\delta=\pm$ and some finite $n$ and $l$ (cf.\cite[Section 7]{Kanade}). 
By using Lemma \ref{NGK01} and (\ref{sing1104}), we have
\begin{equation*}
\begin{split}
&
\begin{pmatrix}
\langle L_0\psi^*,v^\delta_1\otimes v^{(n)}_{\frac{2k-1}{2}}\rangle  \\
\langle L_0\psi^*,(L_{-1}v^\delta_1)\otimes v^{(n)}_{\frac{2k-1}{2}}\rangle \\
\end{pmatrix}
\\
&=\begin{pmatrix}
 h_{2,2}+h_{2n+2,2i} & 1 \\
 \frac{2m^2}{2m+1}h_{2n+2,2i}& h_{2,2}+h_{2n+2,2i}+1-\frac{(2m+1)^2+1}{2(2m+1)}
\end{pmatrix}
\begin{pmatrix}
\langle \psi^*, v^\delta_1\otimes v^{(n)}_{\frac{2k-1}{2}}\rangle  \\
\langle \psi^*,(L_{-1}v^\delta_1)\otimes v^{(n)}_{\frac{2k-1}{2}}\rangle \\
\end{pmatrix}
.
\end{split}
\end{equation*}
We see that the eigenvalues of this matrix do not correspond to the minimal conformal weights of all simple $\mathcal{SW}(m)$-modules. Thus we have
\begin{align*}
&\langle \psi^*, v^\delta_1\otimes v^{(n)}_{\frac{2k-1}{2}}\rangle=0,
&\langle \psi^*,(L_{-1}v^\delta_1)\otimes v^{(n)}_{\frac{2k-1}{2}}\rangle=0
\end{align*}
for any $n\in \mathbb{Z}_{\geq 1}$ and $k\in \mathbb{Z}$. Similary we can show that
\begin{align*}
&\langle \psi^*, G_{-\frac{1}{2}}v^\delta_1\otimes v^{(n)}_{\frac{2k-1}{2}}\rangle=0,
&\langle \psi^*,G^3_{-\frac{1}{2}}v^\delta_1\otimes v^{(n)}_{\frac{2k-1}{2}}\rangle=0
\end{align*}
for any $n\in \mathbb{Z}_{\geq 1}$ and $k\in \mathbb{Z}$. Therefore we obtain (\ref{kieru}).

By Proposition \ref{sl2action2}, Lemma \ref{NGK01}, (\ref{sing1104}) and (\ref{kieru}), we see that, depending on whether $\psi^*$ is even or odd, $\langle\psi^*,\phi_1\otimes \phi_2\rangle$ is determined by the numbers
\begin{align}
\label{even0718}
&\langle\psi^*,G^l_{-\frac{1}{2}}v^\delta_1\otimes v^{\delta'}_i\rangle\ \ \ (l=0,2,\ \ \delta=\pm,\ \delta'=\pm)
\end{align}
or
\begin{align}
\label{odd0718}
&\langle\psi^*,G^l_{-\frac{1}{2}}v^\delta_1\otimes v^{\delta'}_i\rangle\ \ \ (l=1,3,\ \ \delta=\pm,\ \delta'=\pm).
\end{align}

Let us assume that the parity of $\psi^*$ is odd. Then $\langle\psi^*,\phi_1\otimes \phi_2\rangle$ is determined by the numbers (\ref{odd0718}).
By using Lemma \ref{NGK01} and the relation (\ref{sing1104}), we have
\begin{equation*}
\begin{pmatrix}
\langle L_0\psi^*,G_{-\frac{1}{2}}v^\delta_1\otimes v^{\delta'}_i\rangle  \\
\langle L_0\psi^*,G^3_{-\frac{1}{2}}v^\delta_1\otimes v^{\delta'}_i\rangle \\
\end{pmatrix}
=M_1
\begin{pmatrix}
\langle \psi^*, G_{-\frac{1}{2}}v^\delta_1\otimes v^{\delta'}_i\rangle  \\
\langle \psi^*,G^3_{-\frac{1}{2}}v^\delta_1\otimes v^{\delta'}_i\rangle \\
\end{pmatrix}
,
\end{equation*}
where
\begin{equation*}
M_1=
\begin{pmatrix}
 h_{2,2}+h_{2,2i}+\frac{1}{2} & 1 \\
 \frac{m^2}{2m+1}+\frac{2m^2}{2m+1}h_{2,2i}& h_{2,2}+h_{2,2i}+\frac{3}{2}-\frac{(2m+1)^2+1}{2(2m+1)}.
\end{pmatrix}
\end{equation*}
We see that this matrix $M_1$ is diagonalizable and the eigenvalues are given by $h_{1,2i+1}$ and $h_{3,2i-1}$. Note that the eigenvalue $h_{3,2i-1}$ does not correspond to any minimal conformal weight of the simple $\mathcal{SW}(m)$-modules. Thus $L_0$ acts semisimply on $\psi^*$ and the $L_0$-weight of $\psi^*$ is $h_{1,2i+1}$ which is the minimal conformal weight of $X_{2i+1}$.

Next let us assume that the parity of $\psi^*$ is even. Then $\langle\psi^*,\phi_1\otimes \phi_2\rangle$ is determined by the numbers (\ref{even0718}).
By using Lemma \ref{NGK01} and (\ref{sing1104}), we have
\begin{equation*}
\begin{pmatrix}
\langle L_0\psi^*,v^{\delta}_1\otimes v^{\delta'}_i\rangle  \\
\langle L_0\psi^*,(L_{-1}v^\delta_1)\otimes v^{\delta'}_i\rangle \\
\end{pmatrix}
=M_2
\begin{pmatrix}
\langle \psi^*, v^\delta_1\otimes v^{\delta'}_i\rangle  \\
\langle \psi^*,(L_{-1}v^\delta_1)\otimes v^{\delta'}_i\rangle \\
\end{pmatrix}
,
\end{equation*}
where
\begin{equation*}
M_2=\begin{pmatrix}
 h_{2,2}+h_{2,2i} & 1 \\
 \frac{2m^2}{2m+1}h_{2,2i}& h_{2,2}+h_{2,2i}+1-\frac{(2m+1)^2+1}{2(2m+1)}
\end{pmatrix}
.
\end{equation*}
We see that this matrix $M_2$ is diagonalizable and eigenvalues are given by $h_{1,2i-1}$ and $h_{3,2i+1}$. Note that the eigenvalue $h_{3,2i+1}$ does not correspond to any minimal conformal weight of the simple $\mathcal{SW}(m)$-modules. Thus $L_0$ acts semisimply on $\psi^*$ and the $L_0$-weight of $\psi^*$ is $h_{1,2i-1}$ which is the minimal conformal weight of $X_{2i-1}$.  

Hence the $L_0$-weight of $\psi^*$ is given by $h_{1,2i-1}$ or $h_{1,2i+1}$.
Note that from Proposition \ref{sl2action2}, $\widehat{W}^\pm_0$ act trivially on the minimal conformal weight spaces of $X_{2i-1}$ and $X_{2i+1}$. Then we have
$
\widehat{W}^\pm_0\psi^*=0.
$
Thus, by Lemma \ref{NGK01}, we see that 
$\langle\psi^*,\phi_1\otimes \phi_2\rangle$ is determined by the numbers
\begin{align*}
&\langle\psi^*,v^+_1\otimes v^{-}_i\rangle,\ \ \ \ \ \ \ \ \ 
&\langle\psi^*,L_{-1}v^+_1\otimes v^{-}_i\rangle
\end{align*}
in the case of $L_0\psi^*=h_{1,2i-1}\psi^*$, and 
\begin{align*}
&\langle\psi^*,G_{-\frac{1}{2}}v^+_1\otimes v^{-}_i\rangle,\ \ \ \ \ 
&\langle\psi^*,G^{3}_{-\frac{1}{2}}v^+_1\otimes v^{-}_i\rangle
\end{align*}
in the case of $L_0\psi^*=h_{1,2i+1}\psi^*$.
Let $\binom{\mu_1}{\mu_2}$ be an eigenvector of ${}^tM_1$ with the eigenvalue $h_{3,2i-1}$ and $\binom{\nu_1}{\nu_2}$ be an eigenvector of ${}^tM_2$ with the eigenvalue $h_{3,2i+1}$.
Assume $L_0\psi^*=h_{1,2i+1}\psi^*$. Then we have
\begin{equation*}
\begin{split}
h_{1,2i+1}(\mu_1,\mu_2)
\begin{pmatrix}
\langle \psi^*,G_{-\frac{1}{2}}v^+_1\otimes v^{-}_i\rangle  \\
\langle \psi^*,G^3_{-\frac{1}{2}}v^+_1\otimes v^{-}_i\rangle \\
\end{pmatrix}
&=
(\mu_1,\mu_2)
\begin{pmatrix}
\langle L_0\psi^*,G_{-\frac{1}{2}}v^+_1\otimes v^{-}_i\rangle  \\
\langle L_0\psi^*,G^3_{-\frac{1}{2}}v^+_1\otimes v^{-}_i\rangle \\
\end{pmatrix}
\\
&=(\mu_1,\mu_2)M_1
\begin{pmatrix}
\langle \psi^*, G_{-\frac{1}{2}}v^+_1\otimes v^{-}_i\rangle  \\
\langle \psi^*,G^3_{-\frac{1}{2}}v^+_1\otimes v^{-}_i\rangle \\
\end{pmatrix}
\\
&=h_{3,2i-1}(\mu_1,\mu_2)
\begin{pmatrix}
\langle \psi^*, G_{-\frac{1}{2}}v^+_1\otimes v^{-}_i\rangle  \\
\langle \psi^*,G^3_{-\frac{1}{2}}v^+_1\otimes v^{-}_i\rangle \\
\end{pmatrix}
.
\end{split}
\end{equation*}
Thus we obtain
\begin{equation*}
\mu_1\langle \psi^*,G_{-\frac{1}{2}}v^+_1\otimes v^{-}_i\rangle+\mu_2\langle \psi^*,G^3_{-\frac{1}{2}}v^+_1\otimes v^{-}_i\rangle=0.
\end{equation*}
Similarly, assuming $L_0\psi^*=h_{1,2i-1}\psi^*$, we obtain
\begin{equation*}
\nu_1\langle \psi^*,v^+_1\otimes v^{-}_i\rangle+\nu_2\langle \psi^*,L_{-1}v^+_1\otimes v^{-}_i\rangle=0.
\end{equation*}
Therefore the vector space $A_0((X_2\boxtimes X_{2i})^*)$ is at most two dimensional.
\end{proof}

\begin{lemma}
\label{X1X1}
We have
\begin{eqnarray}
\label{2022110400}
X_2\boxtimes X_2=X_1\oplus \Gamma(X_3),
\end{eqnarray}
where $\Gamma(X_3)$ is a lowest weight module whose top composition factor is $X_3$.
\begin{proof}
From Lemmas \ref{Intertwining} and \ref{NGKlem}, we have
\begin{equation*}
X_2\boxtimes X_2=\Gamma(X_1)\oplus \Gamma(X_3),
\end{equation*} 
where $\Gamma(X_{2i+1})$ $(i=0,1)$ are lowest weight modules whose top composition factors are given by $X_{2i+1}$. Let us show $\Gamma(X_1)=X_1$. Assume $\Gamma(X_1)\ncong X_1$. 
Then by Proposition \ref{AMsplit}, $\Gamma(X_1)$ must have a composition factor $X_{2m}$.
In particular, we have $\psi^*\in A_0((X_2\boxtimes X_{2})^*)$ such that
\begin{equation*}
L_0\psi^*=h_{1,1}\psi^*=0,\ \ \ \ \ \ \ \ \ \ \ G_{-\frac{1}{2}}\psi^*\neq 0.
\end{equation*}
Note that $G_{-\frac{1}{2}}\psi^*\in A_0((X_2\boxtimes X_{2})^*)$. Then, from the proof of Lemma \ref{NGKlem}, we see that the $L_0$-weight of $G_{-\frac{1}{2}}\psi^*$ must be $h_{1,3}$. But, since $L_0G_{-\frac{1}{2}}\psi^*=h_{2,2m}G_{-\frac{1}{2}}\psi^*$, we have a contradiction.
\end{proof}
\end{lemma}

Recall that in Subsection \ref{ConSol}, we construct the solutions $\{\Psi^+_{i,j}(z)\ |\ i,j=1,0\}$ and $\{\Psi^-_{i,j}(z)\ |\ i,j=1,0\}$ of the Fuchsian differential equation (\ref{fuchsian}) (see (\ref{Def1011}) and Proposition \ref{propfund}).
By the connection formulas in Proposition \ref{connection}, we can show the following theorem (cf. \cite[Theorem 4.7]{CMOY}).
\begin{theorem}
\label{rigid12}
$X_2$ is rigid and self-dual.
\begin{proof}
We show the rigidity of $X_2$ following the methods in \cite{CMY,CMY2,McRae,TWFusion}. By Lemma \ref{X1X1}, we have parity-homogeneous homomorphisms
\begin{align*}
&i_{1}:X_1\rightarrow X_2\boxtimes X_2,
&&p_1:X_2\boxtimes X_2\rightarrow X_1,\\
&i_{3}:\Gamma(X_3)\rightarrow X_2\boxtimes X_2,
&&p_3:X_2\boxtimes X_2\rightarrow \Gamma(X_3)
\end{align*}
such that
\begin{equation}
\label{1017NH}
\begin{split}
&p_1\circ i_1={\rm id}_{X_1},\ \ \ \ \ p_3\circ i_3={\rm id}_{\Gamma(X_3)},\\
&i_1\circ p_1+i_3\circ p_3={\rm id}_{X_2\boxtimes X_2},
\end{split}
\end{equation}
where $\Gamma(X_3)$ is the lowest weight module defined by (\ref{2022110400}). The maps $i_1$ and $p_1$ are candidates of the coevaluation and evaluation, respectively.

We define two homomorphisms $f,g:X_2\rightarrow X_2$ as the compositions
\begin{equation*}
\begin{split}
&f:X_2  \xrightarrow{r^{-1}} X_2{\boxtimes} X_1 \xrightarrow{{\rm id}\boxtimes i_1} X_2{\boxtimes}(X_2{\boxtimes} X_2)\xrightarrow{\mathcal{A}}(X_2{\boxtimes} X_2){\boxtimes} X_2\xrightarrow{p_1\boxtimes{\rm id}}X_1{\boxtimes} X_2\xrightarrow{l} X_2, \\
&g:X_2  \xrightarrow{l^{-1}} X_1{\boxtimes} X_2 \xrightarrow{i_1\boxtimes{\rm id} } (X_2{\boxtimes} X_2){\boxtimes} X_2\xrightarrow{\mathcal{A}^{-1}}X_2{\boxtimes} (X_2{\boxtimes} X_2)\xrightarrow{{\rm id}\boxtimes p_1}X_2{\boxtimes} X_1\xrightarrow{r} X_2,
\end{split}
\end{equation*}
where $\mathcal{A}$ is the associativity isomorphism, and $l$, $r$ are the left and right unit isomorphisms. The left and right unit isomorphisms $l,r$ are characterized by
\begin{equation*}
\begin{split}
\bar{l}_{X_2}(u_{X_1}\boxtimes u_{X_2})&=Y_{X_2}(u_{X_1},1)u_{X_2},\\
\bar{r}_{X_2}(u_{X_2}\boxtimes u_{X_1})&=(-1)^{|u_{X_1}||u_{X_2}|}e^{L_{-1}}Y_{X_2}(u_{X_1},-1)u_{X_2},
\end{split}
\end{equation*}
for parity-homogeneous $u_{X_1}\in X_1$ and $u_{X_2}\in X_2$.
By Lemma 4.2.1 and Corollary 4.2.2 of \cite{CMY2},
it is enough to show that one of $f$ and $g$ is non-zero in order to show that $X_2$ is rigid and self-dual.
Let us show $f\neq 0$. 

Let $\mathcal{Y}_{2\boxtimes {2}}$ and $\mathcal{Y}_{2\boxtimes (2\boxtimes 2)}$ be non-zero even intertwining operators of type
\begin{equation*}
\begin{pmatrix}
\ X_2\boxtimes X_2 \\
X_2\ \ X_2
\end{pmatrix}
,
\ \ \ \ \ \ 
\begin{pmatrix}
\ X_2\boxtimes(X_2\boxtimes X_2) \\
X_2\ \ X_2\boxtimes X_2
\end{pmatrix}
,
\end{equation*}
respectively. 
We introduce even intertwining operators
\begin{align*}
&\mathcal{Y}^2_{21}=l_{X_2}\circ (p_1\boxtimes {\rm id}_{X_2})\circ \mathcal{A}_{X_2,X_2,X_2}\circ \mathcal{Y}_{2\boxtimes (2\boxtimes 2)}\circ ({\rm id}_{X_2}\otimes i_1),\\
&\mathcal{Y}^2_{23}=l_{X_2}\circ (p_1\boxtimes {\rm id}_{X_2})\circ \mathcal{A}_{X_2,X_2,X_2}\circ \mathcal{Y}_{2\boxtimes (2\boxtimes 2)}\circ ({\rm id}_{X_2}\otimes i_3).
\end{align*} 
The first intertwining operator corresponds to $f$.
Fix a minimal conformal weight vector $v\in X_2[h_{2,2}]$ and let $v^*$ be a minimal conformal weight vector of $X^*_2(\simeq X_2)$ such that
$
\langle v^*, v\rangle\neq 0. 
$ 
For some $x\in\mathbb{R}$ such that $1>x>1-x>0$, we set
\begin{align*}
&\phi_1(x)=\langle v^*,\mathcal{Y}^2_{2{1}}(v,1)(p_1\circ\mathcal{Y}_{2\boxtimes{2}})(v,x)v\rangle,\\
&\phi_3(x)=\langle v^*,\mathcal{Y}^2_{2{3}}(v,1)(p_3\circ\mathcal{Y}_{2\boxtimes{2}})(v,x)v\rangle.
\end{align*}
To prove $f\neq 0$, it is enough to show that
$
\phi_1(x)\neq 0.
$
Note that $\phi_1$ and $\phi_3$ satisfy the Fuchsian differential equation (\ref{fuchsian}), and admit series
\begin{align*}
&\phi_1(x)\in \mathbb{C}x^{h_{1,1}-2h_{2,2}}\bigl(1+x\mathbb{C}[[x]]\bigr),
&\phi_3(x)\in \mathbb{C}x^{h_{1,3}+\frac{1}{2}-2h_{2,2}}\bigl(1+x\mathbb{C}[[x]]\bigr).
\end{align*}
Then, by noting the characteristic exponents of (\ref{fuchsian}), we have
\begin{equation}
\label{1102eq}
\phi_1(x)\in \mathbb{C}\Psi^+_{0,1}(x)+\mathbb{C}\Psi^+_{0,0}(x),\ \ \ \ \ \ \ \phi_3(x)\in \mathbb{C}\Psi^+_{1,1}(x)+\mathbb{C}\Psi^+_{1,0}(x).
\end{equation}
By (\ref{1017NH}), we have
\begin{equation}
\label{F}
\begin{split}
&\phi_1(x)+\phi_3(x)\\
&=\langle v^*, \overline{l_{X_2}\circ(p_1\boxtimes {\rm id}_{X_2})\circ \mathcal{A}_{X_2,X_2,X_2} }\bigl(\mathcal{Y}_{2\boxtimes(2\boxtimes {2})}(v,1)\mathcal{Y}_{2\boxtimes {2}}(v,x)v\bigr)\rangle\\
&=\langle v^*,  \overline{l_{X_2}\circ(p_1\boxtimes {\rm id}_{X_2}) } \bigl(\mathcal{Y}_{(2\boxtimes 2)\boxtimes{2}}(\mathcal{Y}_{2\boxtimes 2}(v,1-x)v,x)v\bigr)\rangle\\
&=\langle v^*, \overline{l_{X_2}}\bigl(\mathcal{Y}_{1\boxtimes {2}}((p_1\circ\mathcal{Y}_{2\boxtimes 2})(v,1-x)v,x)v\bigr) \rangle,\\
&=\langle v^*, Y_{X_2}\bigl((p_1\circ\mathcal{Y}_{2\boxtimes 2})(v,1-x)v,x)\bigr)v\rangle,
\end{split}
\end{equation}
where $\mathcal{Y}_{1\boxtimes {2}}$ is a non-zero even intertwining operator of type
{\scriptsize{$\begin{pmatrix}
   \ X_2  \\
   X_1\ X_2
\end{pmatrix}
$}}.
Note that $p_1\circ\mathcal{Y}_{2\boxtimes 2}$ is a non-zero even intertwining operator of type
{\scriptsize{$\begin{pmatrix}
   \ X_1  \\
   X_2\ X_2
\end{pmatrix}
$}}. Then we have
\begin{equation}
\label{fourpoint}
\begin{split}
\langle v^*, Y_{X_2}\bigl((p_1\circ\mathcal{Y}_{2\boxtimes 2})(v,1-x)v,x)\bigr)v\rangle\in \mathbb{C}^\times\Psi^-_{0,0}(x)+\mathbb{C}\Psi^-_{0,1}(x).
\end{split}
\end{equation}
Assume that $\phi_1(x)=0$. Then, from (\ref{F}) and (\ref{fourpoint}), we have $\phi_3(x)\neq 0$ and 
\begin{equation*}
\phi_3(x)\in \mathbb{C}^\times\Psi^-_{0,0}(x)+\mathbb{C}\Psi^-_{0,1}(x).
\end{equation*}
Then by (\ref{1102eq}), we see that there exist $(k,l)\in \mathbb{C}^2\setminus\{(0,0)\}$ such that
\begin{equation}
\label{cont0203}
k\Psi^+_{1,1}(x)+l\Psi^+_{1,0}(x)\in \mathbb{C}^\times\Psi^-_{0,0}(x)+\mathbb{C}\Psi^-_{0,1}(x).
\end{equation}
On the other hand, from Proposition \ref{connection}, we have
\begin{equation}
\label{rel0203}
\begin{split}
\Psi^+_{1,1}(x)&=\frac{(-1)^m}{4{\rm cos}(\pi\alpha_-\beta_{2,2})}(\Psi^-_{1,1}(x)-\Psi^-_{0,1}(x)-\Psi^-_{1,0}(x)+\Psi^-_{0,0}(x)),\\
\Psi^+_{1,0}(x)&=\frac{(-1)^m}{4{\rm cos}(\pi\alpha_-\beta_{2,2})}(-3\Psi^-_{1,1}(x)+3\Psi^-_{0,1}(x)-\Psi^-_{1,0}(x)+\Psi^-_{0,0}(x)).
\end{split}
\end{equation}
By (\ref{rel0203}), we see that (\ref{cont0203}) contradicts the linear independence of $\{\Psi^-_{i,j}(x)\}$.
Therefore, we obtain $\phi_1(x)\neq 0$.
\end{proof}
\end{theorem}

\subsection{Non-semisimple fusion rules}
\label{SecFusion}

\begin{lemma}
\label{NGKlem2}
\mbox{}
For $i=1,\dots,m-1$, the vector space $A_0((X_2\boxtimes X_{2i+1})^*)$ is at most four dimensional. $L_0$ acts semisimply on $A_0((X_2\boxtimes X_{2i+1})^*)$ and any $L_0$ eigenvalue of this space is contained in $\{h_{2,2i},h_{2,2i+2}\}$, where $h_{2,2i}$ and $h_{2,2i+2}$ are the minimal conformal weights of $X_{2i}$ and $X_{2i+2}$, respectively. 
\begin{proof}
By Lemma \ref{Intertwining}, the tensor product $X_2\boxtimes X_{2i+1}$ is non zero.
Let $\psi^*$ be an arbitrary non zero $L_0$ homogeneous vector of $A_0((X_2\boxtimes X_{2i+1})^*)$.
Let $\phi_1$ and $\phi_2$ be arbitrary non zero $L_0$ homogeneous vectors of $X_{2}$ and $X_{2i+1}$ such that $\langle\psi^*,\phi_1\otimes \phi_2\rangle\neq 0$, where we use the shorthand notation (\ref{abr0204}) as $I={\boxtimes_{P(1)}}$. Let $\{v^+,v^-\}$ be a basis of the minimal conformal weight space of $X_{2}$ such that 
$
\widehat{W}^\pm_0v^\pm=0
$
and
$
\widehat{W}^\pm_0v^\mp\in\mathbb{C}^\times v^\pm.
$
For $n\geq 1$, let $w^{(n)}_k(k=-n,\dots,n)$ be the minimal conformal weight vectors of the subspace $(2n+1)L(h_{2n+1,2i+1})\subset X_{2i+1}$ defined in Proposition \ref{dfn0830}. Similar to the arguments in Lemma \ref{NGKlem}, we have
\begin{align}
\label{kieru2}
\langle\psi^*, U(\mathfrak{ns})v^\pm\otimes w^{(n)}_k\rangle=0\ \ \ (n\geq 2,\ -n\leq k\leq n).
\end{align}
Let $u(=w^{(0)}_0)$ be the minimal conformal weight vector of $X_{2i+1}$. Then, by Lemma \ref{NGK01} and by the relations (\ref{sing1104}) and (\ref{kieru2}), we see that $\langle\psi^*,\phi_1\otimes \phi_2\rangle$ is determined by the numbers
\begin{align*}
&\langle\psi^*, G^l_{-\frac{1}{2}}v^\pm\otimes u\rangle, 
&\langle\psi^*, G^l_{-\frac{1}{2}}v^\pm\otimes w^{(1)}_k\rangle,
\end{align*}
for $l=0,1,2,3$ and $k=-1,0,1$.
From Proposition \ref{sl2action2}, we have
\begin{align}
\label{04170}
\begin{aligned}
&\widehat{W}^\pm_{-h}v^\delta\in U(\mathfrak{ns})v^\delta+U(\mathfrak{ns})v^{-\delta}\\
&{W}^\pm_{-h}v^\delta\in U(\mathfrak{ns})v^\delta+U(\mathfrak{ns})v^{-\delta}
\end{aligned}
\ \ \ \ \ h< h_{3,1}=\frac{1}{2}+2m,\ \delta=\pm,
\end{align}
and
\begin{align}
\label{041700}
\begin{aligned}
&w^{(1)}_{\pm 1}\in \mathbb{C}^\times W^\pm [h_{1,2i+1}-h_{3,2i+1}]u+U(\mathfrak{ns})u,\\
&w^{(1)}_{0}\in \mathbb{C}^\times W^0 [h_{1,2i+1}-h_{3,2i+1}]u+U(\mathfrak{ns})u.
\end{aligned}
\end{align}
From the definition (\ref{widehatW}), we have
\begin{align*}
&[G_{-\frac{1}{2}},Y({W}^a,z)]\in \mathbb{C}^\times Y(\widehat{W}^a,z),
&[G_{-\frac{1}{2}},Y(\widehat{W}^a,z)]\in \mathbb{C}^\times Y(L_{-1}{W}^a,z) 
\end{align*}
for $a=\pm,0$.
Thus, by using (\ref{04170})-(\ref{041700}) and the identities in Lemma \ref{NGK01}, we see that each $\langle\psi^*, G^s_{-\frac{1}{2}}.v^\pm\otimes w^{(1)}_k\rangle$ $(s=0,1,2,3)$ is determined by the numbers
\begin{align*}
&\langle\psi^*, G^{l}_{-\frac{1}{2}}v^\pm\otimes u\rangle\ \ \ \ (l=0,1,2,3).
\end{align*}
Therefore $\langle\psi^*,\phi_1\otimes \phi_2\rangle$ is determined by the numbers
\begin{align}
\label{ev7}
\langle\psi^*,G^l_{-\frac{1}{2}}v^\delta\otimes u\rangle\ \ \ (l=0,2,\ \ \delta=\pm)
\end{align}
or
\begin{align}
\label{od7}
\langle\psi^*,G^l_{-\frac{1}{2}}v^\delta\otimes u\rangle\ \ \ (l=1,3,\ \ \delta=\pm).
\end{align}

(a) Assume that $\langle\psi^*,\phi_1\otimes \phi_2\rangle$ is determined by the numbers (\ref{od7}).
Then, by using Lemma \ref{NGK01} and (\ref{sing1104}), we have
\begin{equation*}
\begin{pmatrix}
\langle L_0\psi^*,G_{-\frac{1}{2}}v^\delta\otimes u\rangle  \\
\langle L_0\psi^*,G^3_{-\frac{1}{2}}v^\delta_1\otimes u\rangle \\
\end{pmatrix}
=N_1
\begin{pmatrix}
\langle \psi^*, G_{-\frac{1}{2}}v^\delta\otimes u\rangle  \\
\langle \psi^*,G^3_{-\frac{1}{2}}v^\delta\otimes u\rangle \\
\end{pmatrix}
,
\end{equation*}
where
\begin{equation*}
N_1=
\begin{pmatrix}
 h_{2,2}+h_{1,2i+1}+\frac{1}{2} & 1 \\
 \frac{m^2}{2m+1}+\frac{2m^2}{2m+1}h_{1,2i+1}& h_{2,2}+h_{1,2i+1}+\frac{3}{2}-\frac{(2m+1)^2+1}{2(2m+1)}
\end{pmatrix}
.
\end{equation*}
We see that this matrix $N_1$ is diagonalizable and the eigenvalues are given by $h_{2,2i}$ and $h_{1,2(m-i)-1}$ which are the minimal conformal weights of $X_{2i}$ and $X_{2(m-i)-1}$ respectively. 

(b) Let us assume that $\langle\psi^*,\phi_1\otimes \phi_2\rangle$ is determined by the numbers (\ref{ev7}).
By using Lemma \ref{NGK01} and (\ref{sing1104}), we have
\begin{equation*}
\begin{pmatrix}
\langle L_0\psi^*,v^{\delta}\otimes u\rangle  \\
\langle L_0\psi^*,(L_{-1}v^\delta)\otimes u\rangle \\
\end{pmatrix}
=N_2
\begin{pmatrix}
\langle \psi^*, v^\delta\otimes u\rangle  \\
\langle \psi^*,(L_{-1}v^\delta)\otimes u\rangle \\
\end{pmatrix}
,
\end{equation*}
where
\begin{equation*}
N_2=\begin{pmatrix}
 h_{2,2}+h_{1,2i+1} & 1 \\
 \frac{2m^2}{2m+1}h_{1,2i+1}& h_{2,2}+h_{1,2i+1}+1-\frac{(2m+1)^2+1}{2(2m+1)}
\end{pmatrix}
.
\end{equation*}
We see that this matrix $N_2$ is diagonalizable and eigenvalues are given by $h_{2,2i+2}$ and $h_{1,2(m-i)+1}$ which are the minimal conformal weights of $X_{2i+2}$ and $X_{2(m-i)+1}$ respectively.

(c) Let us assume that the $L_0$-weight of $\psi^*$ is $h_{1,2(m-i)+1}$ or $h_{1,2(m-i)-1}$. 
Then, similar to the arguments in Lemma \ref{NGKlem}, we see that
\begin{align}
\label{SW1224}
\langle\psi^*, U(\mathfrak{ns})v^\pm\otimes w^{(1)}_k\rangle=0\ \ \ \ (k=1,0,-1).
\end{align}
From the structure of the simple modules $X_{2(m-i)+1}$ and $X_{2(m-i)-1}$, we have
\begin{align}
\label{uvan0}
\widehat{W}^\pm_{0}\psi^*=0.
\end{align}
By (\ref{prop1010}) in Proposition \ref{sl2action2}), we have
\begin{equation}
\label{uvan}
\widehat{W}^\pm_{-h}u=0\ \ \ \ (h<\frac{1}{2}-i+2m).
\end{equation}
Then, by using the second formula in Lemma \ref{NGK01} as $v=\widehat{W}^\pm$ and $n=0$, and by using the relations (\ref{rel20240509}), (\ref{kieru2}), (\ref{SW1224}), (\ref{uvan0}) and (\ref{uvan}), we have
$
\langle\psi^*, v^\pm\otimes u\rangle=0.
$
Thus we have a contradiction. 

Hence, from (a), (b) and (c), the $L_0$-weight of $\psi^*$ is $h_{2,2i}$ or $h_{2,2i+2}$. Note that by Proposition \ref{sl2action2},
\begin{align}
\label{SW12242}
&W^\pm_{-h}v^\pm=0,
&\widehat{W}^\pm_{-h}v^\pm=0
\end{align}
for $h< h_{3,1}$. 
Then, noting the argument just before (a) and using Lemma \ref{NGK01} and (\ref{uvan})-(\ref{SW12242}), we see that 
\begin{align*}
\langle\widehat{W}^\pm_{0}\psi^*, v^\pm\otimes U(\mathfrak{ns})u\rangle=0.
\end{align*}
Thus, noting that the minimal conformal weight spaces of $X_{2i}$ and $X_{2i+2}$ are two dimensional, we see that $\langle\psi^*,\phi_1\otimes \phi_2\rangle$ is determined by the numbers
\begin{align*}
&\langle\widehat{W}^\pm_{0}\psi^*+\psi^*,v^\mp\otimes u\rangle,\ \ \ \ \ \ \ \ \ 
&\langle\widehat{W}^\pm_{0}\psi^*+\psi^*,L_{-1}v^\mp\otimes u\rangle
\end{align*}
or
\begin{align*}
&\langle\widehat{W}^\pm_{0}\psi^*+\psi^*,G_{-\frac{1}{2}}v^\mp\otimes u\rangle,\ \ \ \ \ \ \ \ \ 
&\langle\widehat{W}^\pm_{0}\psi^*+\psi^*,G^{3}_{-\frac{1}{2}}v^\mp\otimes u\rangle.
\end{align*}
Let $\binom{\kappa_1}{\kappa_2}$ be an eigenvector of ${}^tN_1$ with the eigenvalue $h_{1,2(m-i)-1}$ and $\binom{\lambda_1}{\lambda_2}$ be an eigenvector of ${}^tN_2$ with the eigenvalue $h_{1,2(m-i)+1}$. Similar to the argument in the proof of Lemma \ref{NGKlem}, we have
\begin{equation*}
\begin{split}
&\kappa_1\langle\widehat{W}^\pm_{0}\psi^*+\psi^*,G_{-\frac{1}{2}}v^\mp\otimes u\rangle+\kappa_2\langle\widehat{W}^\pm_{0}\psi^*+\psi^*,G^{3}_{-\frac{1}{2}}v^\mp\otimes u\rangle=0,\\
&\lambda_1\langle\widehat{W}^\pm_{0}\psi^*+\psi^*,v^\mp\otimes u\rangle+\lambda_2\langle\widehat{W}^\pm_{0}\psi^*+\psi^*,L_{-1}v^\mp\otimes u\rangle=0.
\end{split}
\end{equation*}
Therefore $A_0((X_2\boxtimes X_{2i+1})^*)$ is at most four dimensional.
\end{proof}
\end{lemma}

Let us recall some properties of rigid and dual objects in tensor categories. For the following proposition, see, for example, \cite{Etingof,JS,KL}.
\begin{proposition}
\label{Etingof}
Let $(\mathcal{C},\otimes)$ be a tensor category, then we have:
\begin{enumerate}
\item Let $V$ be a rigid object in $\mathcal{C}$. Then there is a natural adjunction isomorphism
\begin{eqnarray*}
{\rm Hom}_{\mathcal{C}}(U\otimes V,W)\simeq {\rm Hom}_{\mathcal{C}}(U,W\otimes V^\vee),
\end{eqnarray*}
where $U$ and $W$ are any objects in $\mathcal{C}$, and $V^\vee$ is the dual object of $V$. 
\item Let $V_1$ and $V_2$ be rigid objects in $\mathcal{C}$. Then $V_1\otimes V_2$ is also rigid and 
\begin{align*}
(V_1\otimes V_2)^\vee=V^\vee_2\otimes V^\vee_1.
\end{align*}
\item Let $V$ be a rigid object in $\mathcal{C}$ and let $P$ be a projective object in $\mathcal{C}$. Then $V^\vee\otimes P$ is projective.
\item Assume that
\begin{itemize}
\item $\mathcal{C}$ has enough projective and injective objects.
\item All projective objects are injective and all injective objects are projective.
\item All projective objects are rigid.
\end{itemize}
Then if
\begin{align*}
0\rightarrow V_1\rightarrow V_2\rightarrow V_3\rightarrow 0
\end{align*}
is an exact sequence in $\mathcal{C}$ such that two of $V_1$, $V_2$, $V_3$ are rigid, then the third object is also rigid.
\end{enumerate}
\end{proposition}

\begin{proposition}
\label{fusionXX}
For $s=2,\dots,2m$, we have
\begin{eqnarray*}
X_2\boxtimes X_s=X_{s-1}\oplus X_{s+1}.
\end{eqnarray*}
\begin{proof}
First, let us show
\begin{eqnarray}
\label{eq0507}
X_2\boxtimes X_{2i}=X_{2i-1}\oplus X_{2i+1}\ \ \ \ (i=1,\dots,m).
\end{eqnarray}
By Lemma \ref{Intertwining}, we have
\begin{equation}
I_{\mathcal{SW}(m)}
\begin{pmatrix}
X_{2i-1} \\
X_2\ \ X_{2i}
\end{pmatrix}
\neq 0
,\ \ \ \\ \ \ \ \ \ \ \ 
I_{\mathcal{SW}(m)}
\begin{pmatrix}
X_{2i+1} \\
X_2\ \ X_{2i}
\end{pmatrix}
\neq 0.
\label{Int12140}
\end{equation}
By Lemma \ref{NGKlem}, Proposition \ref{Etingof} and the self-duality of $X_2$, we see that
\begin{eqnarray*}
{\rm Hom}_{\mathcal{SW}(m)}(X_{2(m-i+1)}\oplus X_{2(m-i)},X_2\boxtimes X_{2i})=0.
\end{eqnarray*}
Thus, by (\ref{Int12140}), we obtain (\ref{eq0507}).

Next, let us show
\begin{eqnarray}
\label{eq05072}
X_2\boxtimes X_{2i+1}=X_{2i}\oplus X_{2i+2}\ \ \ \ (i=0,\dots,m-1).
\end{eqnarray}
By Lemma \ref{Intertwining}, we have
\begin{equation}
I_{\mathcal{SW}(m)}
\begin{pmatrix}
X_{2i} \\
X_2\ \ X_{2i+1}
\end{pmatrix}
\neq 0
,\ \ \ \\ \ \ \ \ \ \ \ 
I_{\mathcal{SW}(m)}
\begin{pmatrix}
X_{2i+2} \\
X_2\ \ X_{2i+1}
\end{pmatrix}
\neq 0.
\label{Int12240}
\end{equation}
By Lemma \ref{NGKlem2}, Proposition \ref{Etingof} and the self-duality of $X_2$, we see that
\begin{eqnarray*}
{\rm Hom}_{\mathcal{SW}(m)}(X_{2(m-i)-1}\oplus X_{2(m-i)+1},X_2\boxtimes X_{2i+1})=0.
\end{eqnarray*}
Thus, by (\ref{Int12240}), we obtain (\ref{eq05072}).
\end{proof}
\end{proposition}

Since $\mathcal{SW}(m)$ is $C_2$-cofinite, every simple module has a projective cover \cite{H}. 
In the following, we introduce an indecomposable $\mathcal{SW}(m)$-module and some lemmas to determine the structure of these projective covers.
Let us introduce an operator $\Delta_{Q_-}(-,z)-:\mathcal{SW}(m)\times \mathcal{V}_{L}\rightarrow \overline{\mathcal{V}_{L+\gamma_{2m}}}[[z,z^{-1}]]$ as follows
\begin{align*}
\Delta_{Q_-}(A,z)=\sum_{n\geq 1}\frac{(-1)^{n+1}}{nz^n}\oint_{z'=z}(z'-z)^nQ_-(z')Y(A,z){\rm d}z'\ \ (A\in \mathcal{SW}(m)).
\end{align*}
This operator $\Delta_{Q_-}$ is called logarithmic deformation \cite{FJ} or Li's operator \cite{Li}.
Let
$
W=\mathcal{V}_{L}\oplus \mathcal{V}_{L+\gamma_{2m}}(\in \mathcal{B}_{1})
$
and let $Y_{W}$ be the ordinary $\mathcal{SW}(m)$-action on $W$.
We define an operator $J(-,z)-:\mathcal{SW}(m)\times W\rightarrow \overline{W}[[z,z^{-1}]]$ as follows
\begin{equation*}
J(A,z)
=
\begin{cases}
Y_{W}(A,z)+\Delta_{Q_-}(A,z)& on\ \mathcal{V}_{L},\\
Y_{W}(A,z) & on\ \mathcal{V}_{L+\gamma_{2m}},
\end{cases}
\end{equation*}
where $A\in \mathcal{SW}(m)$.
By the results in \cite{FJ,Li}, the operator $J$ defines an $\mathcal{SW}(m)$-action on $W$.
We set $\mathcal{P}=(W,J)$.
For the conformal vector $T$, the action $J(T,z)$ on the subspace $\mathcal{V}_{L}\subset \mathcal{P}$ is given by
\begin{align*}
J(T,z)=T(z)+\frac{Q_-(z)}{z}.
\end{align*}
Thus, by Propositions \ref{BRST}, $\mathcal{P}$ is indecomposable and has $L_0$-nilpotent rank two.
By the definition, $\mathcal{P}$ has $X_{2m}$ as a submodule and a quotient, and the total composition factors of $\mathcal{P}$ are given by $X_1\oplus X_1\oplus X_{2m}\oplus X_{2m}$.

\begin{lemma}
\label{indecomposable}
$\mathcal{P}/X_{2m}$ or $\mathcal{P}^*/X_{2m}$ is indecomposable.
\begin{proof}
Let $\overline{X}_{2m}$ be the minimal conformal weight space of $X_{2m}$. By the results in \cite{KW,Zhu}, $\overline{X}_{2m}$ has the structure of an $A(\mathcal{SW}(m))$-module. Then, from Theorem \ref{simpleclass}, we see that
\begin{equation}
\label{induceX}
{\rm Ind}^{\mathcal{SW}(m)}_{A(\mathcal{SW}(m))}\overline{X}_{2m}\simeq X_{2m}.
\end{equation}
If neither $\mathcal{P}/X_{2m}$ nor $\mathcal{P}^*/X_{2m}$ is indecomposable, then, as the quotient of $\mathcal{P}$ or $\mathcal{P}^*$, we have a non trivial extension in
\begin{align*}
{\rm Ext}^1_{\mathcal{SW}(m)}(X_{2m},X_{2m}).
\end{align*}
But from (\ref{induceX}), this contradicts Theorem \ref{AM}.
\end{proof}
\end{lemma}

The following lemma can be proved in the same way as Lemmas \ref{NGKlem} and \ref{NGKlem2}.
\begin{lemma}
\label{NGKlem3}
\mbox{}
The vector space $A_0((X_2\boxtimes X_{2m+1})^*)$ is at most six dimensional.
Any $L_0$ eigenvalue of $A_0((X_2\boxtimes X_{2m+1})^*)$ is contained in $\{h_{1,1},h_{2,2m}\}=\{0,\frac{1}{2}\}$, where $h_{1,1}$ and $h_{2,2m}$ are the minimal conformal weights of $X_{1}$ and $X_{2m}$, respectively. 
\end{lemma}

\begin{lemma}
\label{Hom12140}
For any simple $\mathcal{SW}(m)$-module $X$, we have
\begin{equation*}
\begin{split}
&{\rm Hom}_{\mathcal{SW}(m)}(X_2\boxtimes X_{2m+1},X)
=
\begin{cases}
\mathbb{C}&X=X_{2m}, \\
0&\ {\rm otherwise}
\end{cases}
\\
&{\rm Hom}_{\mathcal{SW}(m)}(X,X_2\boxtimes X_{2m+1})
=
\begin{cases}
\mathbb{C}&X=X_{2m},\\
0&\ {\rm otherwise}
\end{cases}
\end{split}
.
\end{equation*}
\begin{proof}
By Propositions \ref{Etingof}-\ref{fusionXX} and the self-duality of $X_2$, we obtain the above eaqualities.
\end{proof}
\end{lemma}

Let us use the following notation:
\begin{enumerate}
\item For $1\leq i\leq m$, let ${P}_{2i}$ be the projective cover of $X_{2(m-i)+1}$.
\item For $0\leq j\leq m-1$, let ${P}_{2j+1}$ be the projective cover of $X_{2(m-j)}$.
\end{enumerate}
By Proposition \ref{Blockdecomp}, we have ${P}_{2i}\in \mathcal{B}_{m-i+1}$ and ${P}_{2j+1}\in \mathcal{B}_{j+1}$.

\begin{proposition}
\label{Proj0}
\begin{enumerate}
\item We have
\begin{align*}
X_2\boxtimes X_{2m+1}=P_{1}.
\end{align*}
\item ${P}_1$ has the socle series
\(
{\rm Soc}_1({P}_1)\subsetneq {\rm Soc}_2({P}_1)\subsetneq {\rm Soc}_3({P}_1)={P}_1
\)
such that
\begin{equation}
\label{SoclePP}
\begin{aligned}
&{\rm Soc}_1({P}_1)\simeq X_{2m},\ {\rm Soc}_2({P}_1)/{\rm Soc}_1({P}_1)\simeq X_1\oplus X_1,\\
&{\rm Soc}_3({P}_1)/{\rm Soc}_2({P}_1)\simeq X_{2m}.
\end{aligned}
\end{equation}
\end{enumerate}
\begin{proof}
Note that from Proposition \ref{Blockdecomp}, $X_{2m+1}$ is projective.
Then, by the self-duality of $X_2$, $X_2\boxtimes X_{2m+1}$ must be projective.
Thus, by Lemma \ref{Hom12140}, we obtain
\begin{align*}
X_2\boxtimes X_{2m+1}={P}_{1}.
\end{align*}
By Lemma \ref{indecomposable} and the projectivity of ${P}_1$, we see that ${P}_1$ has $X_1\oplus X_1$ as composition factors. 
Thus, by Lemmas \ref{NGKlem3}-\ref{Hom12140}, we see that ${P}_1$ satisfies the socle series (\ref{SoclePP}).
\end{proof}
\end{proposition}

\begin{proposition}
\label{projstr}
\begin{enumerate}
\item For each $1\leq s\leq 2m$, the tensor product $X_2\boxtimes {P}_{s}$ is given by:
\begin{itemize}
\item For $s=1$, we have
\begin{align*}
X_2\boxtimes P_1=X_{2m+1}\oplus X_{2m+1}\oplus {P}_2.
\end{align*}
\item For $2\leq s\leq 2m-1$, we have
\begin{align*}
X_2\boxtimes {P}_{s}={P}_{s-1}\oplus {P}_{s+1}.
\end{align*}
\item For $s=2m$, we have
\begin{align*}
X_2\boxtimes {P}_{2m}=X_{2m+1}\oplus X_{2m+1}\oplus {P}_{2m-1}.
\end{align*}
\end{itemize}
\item 
The socle series of the projective covers of the simple modules are given by:
\begin{itemize}
\item For $1\leq i\leq m$, we have
\(
{\rm Soc}_1({P}_{2i})\subsetneq {\rm Soc}_2({P}_{2i})\subsetneq {\rm Soc}_3({P}_{2i})={P}_{2i}
\)
such that
\begin{align*}
&{\rm Soc}_1({P}_{2i})\simeq X_{2(m-i)+1},\ {\rm Soc}_2({P}_{2i})/{\rm Soc}_1({P}_{2i})\simeq X_{2i}\oplus X_{2i},\\
&{\rm Soc}_3({P}_{2i})/{\rm Soc}_2({P}_{2i})\simeq X_{2(m-i)+1}.
\end{align*}
\item For $0\leq i\leq m-1$, we have
\(
{\rm Soc}_1({P}_{2i+1})\subsetneq {\rm Soc}_2({P}_{2i+1})\subsetneq {\rm Soc}_3({P}_{2i+1})={P}_{2i+1}
\)
such that
\begin{align*}
&{\rm Soc}_1({P}_{2i+1})\simeq X_{2(m-i)},\ {\rm Soc}_2({P}_{2i+1})/{\rm Soc}_1({P}_{2i+1})\simeq X_{2i+1}\oplus X_{2i+1},\\
&{\rm Soc}_3({P}_{2i+1})/{\rm Soc}_2({P}_{2i+1})\simeq X_{2(m-i)}.
\end{align*}
\end{itemize}
\end{enumerate}
\begin{proof}
We only prove that $X_2\boxtimes P_1=X_{2m+1}\oplus X_{2m+1}\oplus {P}_2$ and 
\begin{equation}
\label{socle1107}
\begin{split}
&{\rm Soc}_1({P}_{2})\simeq X_{2m-1},\ {\rm Soc}_2({P}_{2})/{\rm Soc}_1({P}_{2})\simeq X_{2}\oplus X_{2},\\
&{P}_{2}/{\rm Soc}_2({P}_{2})\simeq X_{2m-1}.
\end{split}
\end{equation}
The other cases can be proved in the same way.

From Propositions \ref{Etingof}, \ref{fusionXX}, \ref{Proj0} and the self-duality of $X_2$, we can see that for any simple $\mathcal{SW}(m)$-module $X$
\begin{equation}
\begin{split}
&{\rm Hom}_{\mathcal{SW}(m)}(X_2\boxtimes P_1,X)
=
\begin{cases}
\mathbb{C}^2&X=X_{2m+1}, \\
\mathbb{C}&X=X_{2m-1},\\
0&\ {\rm otherwise}
\end{cases}
\\
&{\rm Hom}_{\mathcal{SW}(m)}(X,X_2\boxtimes P_1)
=
\begin{cases}
\mathbb{C}^2&X=X_{2m+1}\\
\mathbb{C}&X=X_{2m-1},\\
0&\ {\rm otherwise}
\end{cases}
\end{split}
.
\label{0618xx}
\end{equation}
Since $X_2$ is rigid and $P_1$ is projective, $X_2\boxtimes P_1$ is also projective. Thus, from (\ref{0618xx}), we obtain $X_2\boxtimes P_1=X_{2m+1}\oplus X_{2m+1}\oplus {P}_2$. Therefore, by the rigidity of $X_2$ and Proposition \ref{fusionXX}, we see that $P_2$ satisfies the socle series (\ref{socle1107}).
\end{proof}
\end{proposition}
In \cite{AM2}, the equivalence between $\mathcal{SW}(m)\mathchar`-{\rm mod}$ and $U^{small}_q(sl_2)\mathchar`-{\rm mod}$ is conjectured,
where $U^{small}_q(sl_2)\mathchar`-{\rm mod}$ is the abelian category of finite dimensional modules over the small quantum group $U^{small}_q(sl_2)$ at $q=e^{\frac{2\pi i}{2m+1}}$. At the level of the arbelian category, we see that this conjecture is true. That is, the following corollary holds.
\begin{corollary}
Two categories $\mathcal{SW}(m)\mathchar`-{\rm mod}$ and $U^{small}_{q}(sl_2)\mathchar`-{\rm mod}$ are equivalent as abelian categories.
\begin{proof}
Similar to the arguments in \cite[Section 6]{NT}, using the structure of the projective $\mathcal{SW}(m)$-modules given in Propositions \ref{Proj0}-\ref{projstr} and the projective $U^{small}_q(sl_2)$-modules classified by \cite{K,Suter,Xia}, we can prove the equivalence of the two abelian categories. We omit the details.
\end{proof}
\end{corollary}
\begin{remark}
Let $q=e^{\frac{2\pi i}{2m+1}}$. The small quantum group $U^{small}_q(sl_2)$ is an associative $\mathbb{C}$-algebra which is generated by $E,F,K,K^{-1}$ satisfying the following fundamental relations
\begin{align*}
&KK^{-1}=K^{-1}K=1,\ KEK^{-1}=q^2E,\ KFK^{-1}=q^{-2}F,\\
&EF-FE=\frac{K-K^{-1}}{q-q^{-1}},\\
&E^{2m+1}=F^{2m+1}=0,\ \ K^{2m+1}=1.
\end{align*}
See \cite{K,Suter,Xia} for details. 
\end{remark}

Finally, let us show that the tensor supercategory $(\mathcal{SW}(m)\mathchar`-{\rm mod},\boxtimes)$ is rigid. Note that the tensor supercategory $(\mathcal{SW}(m)\mathchar`-{\rm mod},\boxtimes)$ is weakly rigid. In fact, we have the following proposition.
\begin{proposition}[\cite{AL,HLZ2,Xu}]
\label{UV}
Let $\mathcal{C}$ be a vertex tensor (super)category with the unit $\bm{1}$.
Given $U,V\in \mathcal{C}$, assume $V^*\in \mathcal{C}$, where $V^*$ is the contragredient of $V$. Then, we have a natural isomorphism
\begin{align*}
{\rm Hom}_{\mathcal{C}}(U,V)\simeq {\rm Hom}_{\mathcal{C}}(U\boxtimes V^*,\bm{1}).
\end{align*}
\end{proposition}

\begin{theorem}
The braided tensor supercategory $(\mathcal{SW}(m)\mathchar`-{\rm mod},\boxtimes)$ is rigid. For any $M\in \mathcal{SW}(m)\mathchar`-{\rm mod}$, we have $M^\vee=M^*$, where $M^\vee$ is the dual of $M$.
\begin{proof}
From Propositions \ref{fusionXX} and \ref{projstr}, all simple and projective modules are rigid and self-dual. By the structure of the projective modules, we see that all indecomposable modules $M\in\SW(m)\mathchar`-{\rm mod}$ except the simple modules and the projective modules satisfy exact sequence
\begin{eqnarray*}
0\rightarrow L\rightarrow M\rightarrow N\rightarrow 0
\end{eqnarray*}
such that $L$ and $N$ are direct sum of simple modules. Then, from Proposition \ref{Etingof}, we see that $M$ is rigid. Since 
\begin{align*}
{\rm Hom}_{\mathcal{SW}(m)}(M\boxtimes M^*,X_1)={\rm Hom}_{\mathcal{SW}(m)}(M,M)\simeq \mathbb{C}
\end{align*}
from Proposition \ref{UV}, we have $M^\vee=M^*$.
\end{proof}
\end{theorem}

\subsection{Fusion rings}
\label{fusion ring}
In this subsection, following the argument in \cite[Subsection 5.3]{TWFusion}, we introduce two fusion rings $P(\mathcal{SW}(m))$ and $K(\mathcal{SW}(m))$, and determine their structure.

Let $\mathbb{I}_{m}$ be the set consisting of all simple modules $X_s (1\leq s\leq 2m+1)$ and all projective modules $P_s (1\leq s\leq 2m)$.
We introduce the free abelian group $P(\mathcal{SW}(m))$ of rank $4m+1$ generated by the elements of $\mathbb{I}_{m}$:
\begin{align*}
P(\mathcal{SW}(m))=\bigoplus_{s=1}^{2m+1}\mathbb{Z}[X_{s}]_P\oplus \bigoplus_{s=1}^{2m}\mathbb{Z}[P_s]_{P}.
\end{align*}
Then, from the results presented in the previous subsection, we can define the structure of a commutative ring on $P(\mathcal{SW}(m))$ such that the product as a ring is given by
\begin{align*}
[M_1]_P\cdot[M_2]_P=[M_1\boxtimes M_2]_P
\end{align*}
for $M_1,M_2\in \mathbb{I}_{m}$, where we extend the symbol $[\bullet]_P$ as follows
\begin{align*}
\Bigl[\bigoplus^n_{i\geq 1}N_i\Bigr]_P=\bigoplus^n_{i\geq 1}[N_i]_P
\end{align*}
for any $N_i\in \mathbb{I}_{m}$ and any $n\in \mathbb{Z}_{\geq 1}$. 

The operator
\begin{align*}
X=X_{2}\boxtimes -
\end{align*}
define a $\mathbb{Z}$-linear endomorphism of $P(\mathcal{SW}(m))$. Thus $P(\mathcal{SW}(m))$ is a module over $\mathbb{Z}[X]$.
We define the following $\mathbb{Z}[X]$-module map
\begin{align*}
\psi:\mathbb{Z}[X]&\rightarrow P(\mathcal{SW}(m)),\\
f(X)&\mapsto f(X)\cdot[X_1]_P.
\end{align*}
Before examining the action of $\mathbb{Z}[X]$ on $P(\mathcal{SW}(m))$, we introduce the following Chebyshev polynomials.
\begin{definition}
We define {\rm Chebyshev polynomials} $U_n(A)$, $n=0,1,\dots\in \mathbb{Z}[A]$ recursively 
\begin{align*}
U_0(A)=1,\ \ \ \ \ \ \ U_1(A)=A,\ \ \ \ \ \ \ U_{n+1}(A)=AU_n(A)-U_{n-1}(A).
\end{align*}
\end{definition}
\begin{remark}
The coefficient of the leading term of any Chebyshev polynomial $U_n(A)$ is 1. Thus we have
\begin{align*}
\mathbb{Z}[A]=\bigoplus_{n=0}^\infty\mathbb{Z}U_n(A).
\end{align*}
\end{remark}

From the results of previous subsection, we obtain the following lemma.
\begin{lemma}
\label{prop01}
\begin{enumerate}
\item For $s=1,\dots,2m+1$, we have
\begin{align*}
[X_s]_P=U_{s-1}(X)[X_1]_{P}.
\end{align*}
\item For $s=1,\dots,2m$, we have
\begin{align*}
[P_s]_{P}=(U_{2m+s}(X)+U_{2m-s}(X))[X_1]_P.
\end{align*}
\item We have the relation
\begin{align*}
U_{4m+1}(X)[X_1]_P=2U_{2m}(X)[X_1]_P.
\end{align*}
\end{enumerate}
\end{lemma}
From this lemma, we obtain the following theorem.
\begin{theorem}
\label{Fusion-ring}
The $\mathbb{Z}[X]$-module map $\psi$ is surjective and the kernel of $\psi$ is given by the ideal
$
{\rm ker}\psi=\langle U_{4m+1}(X)-2U_{2m}(X)\rangle.
$
\begin{proof}
By Lemma \ref{prop01}, we see that $\psi$ is surjective. We define the following ideal of $\mathbb{Z}[X]$
\begin{align*}
I=\langle U_{4m+1}(X)-2U_{2m}(X)\rangle.
\end{align*}
By the third statement in Lemma \ref{prop01}, we see that $I$ is contained in ${\rm ker}\psi$. It is easy to see that the dimension of $\mathbb{Z}[X]/I$ is $4m+1$. Therefore we obtain ${\rm ker}\psi=I$.\\
\end{proof}
\end{theorem}

Next, let us state the results for the Grothendieck fusion ring of $\mathcal{SW}(m)$. We introduce the rank $2m+1$ Grothendieck group
\begin{equation*}
K(\mathcal{SW}(m))=\bigoplus_{s=1}^{2m+1}\mathbb{Z}[X_s]_{K}.
\end{equation*}
From the results presented in the previous subsection, we see that $K(\mathcal{SW}(m))$ has the structure of a commutative ring whose unit object is $[X_1]_K$. 
The operator
$
X=X_{2}\boxtimes -
$
define a $\mathbb{Z}$-linear endomorphism of $K(\mathcal{SW}(m))$. Thus $K(\mathcal{SW}(m))$ is a module over $\mathbb{Z}[X]$.
Then we can define the following $\mathbb{Z}[X]$-module map
\begin{align*}
\phi:\mathbb{Z}[X]&\rightarrow K(\mathcal{SW}(m)),\\
f(X)&\mapsto f(X)\cdot[X_1]_K.
\end{align*}
Similar to the arguments in the case of $P(\mathcal{SW}(m))$, we obtain the following proposition.
\begin{proposition}
\label{Fusion-ring}
The $\mathbb{Z}[X]$-module map $\phi$ is surjective and the kernel of $\phi$ is given by the ideal
$
{\rm ker}\phi=\langle U_{2m+1}(X)-U_{2m-1}(X)-2\rangle.
$
\end{proposition}

\begin{remark}
In \cite[Subsection 5.3]{TWFusion}, a non-semisimple fusion ring $P(\mathcal{W}_p)$ of the triplet algebra $\mathcal{W}_{p}$ is introduced.
As in the case of $P(\mathcal{SW}(m))$, $P(\mathcal{W}_p)$ is defined by adding a ring structure determined from the tensor product to the free abelian group generated from all simple and projective $\mathcal{W}_p$-modules. As shown in \cite{TWFusion}, $P(\mathcal{W}_{p})$ is isomorphic to the quotient ring
\begin{align*}
\frac{\mathbb{Z}[X,Y]}{\langle Y^2-1,U_{2p-1}(X)-2YU_{p-1}(X)\rangle}.
\end{align*}
The $Y$ variable corresponds to the simple current of $\mathcal{W}_p$ and the $X$ variable to a simple $\mathcal{W}_p\mathchar`-$module which has a weight two Virasoro null vector (according to the notation in \cite{FF3,NT,TWFusion}, $Y$ to $\mathcal{X}^-_1$ and $X$ to $\mathcal{X}^+_2$).
For this ring, setting $p=2m+1$ and $Y=1$, we have a quotient ring isomorphic to $P(\mathcal{SW}(m))$ (see Theorem \ref{Fusion-ring}).
Similarly, we can obtain the Grothendieck ring $K(\mathcal{SW}(m))$ as a quotient of the Grothendieck ring $K(\mathcal{W}_{2m+1})$ determined by \cite{TWFusion}.
In \cite{AM2}, it is shown that the characters of the simple $\mathcal{SW}(m)$-modules can be expressed in terms of the characters of the simple $\mathcal{W}_{2m+1}$-modules. From these results, $\mathcal{SW}(m)\mathchar`-{\rm mod}$ and $\mathcal{W}_{2m+1}\mathchar`-{\rm mod}$ seem to be closely related at the level of tensor categories.
\end{remark}
\section*{{Acknowledgement}}
We would like to thank Koji Hasegawa for reading the manuscript and for valuable comments.
We would also like to thank Yuto Moriwaki and Shigenori Nakatsuka for stimulating discussions.

\vspace{10mm}
\ \ \ \ H.~Nakano, \textsc{Advanced Mathematical Institute, Osaka Metropolitan University, Osaka 558-8585, Japan}\par\nopagebreak
  \textit{E-mail address} : \texttt{hiromutakati@gmail.com}


\begin{thebibliography}{99}
\bibitem[AM1]{AM}
     D. Adamovi\'{c} and A. Milas,
     ``On the triplet vertex algebra $W(p)$,''
     \textit{Advances in Mathematics} $\bold{217}$ (2008) 2664-2699.
\bibitem[AM2]{AM20}
D. Adamovi\'{c} and A. Milas, 
``The N= 1 triplet vertex operator superalgebras: twisted sector,''
\textit{SIGMA. Symmetry, Integrability and Geometry: Methods and Applications}, \textbf{4} (2008) 087.     
\bibitem[AM3]{AM2}
     D. Adamovi\'{c} and A. Milas,
     ``The $N=1$ triplet vertex operator superalgebras,''
     \textit{Communications in Mathematical Physics}, $\bold{288}$ (2009) 225-270.
\bibitem[AM4]{AM3}
     D. Adamovi\'{c} and A. Milas,
     ``The structure of Zhu's algebras for certain W$\mathchar`-$algebras,''
     \textit{Advances in Mathematics} $\bold{227}$ (2011) 2425-2456.
\bibitem[ALSW]{AL}
R. Allen, S. Lentner, C. Schweigert and S. Wood, 
           ``Duality structures for module categories of vertex operator algebras and the Feigin Fuchs boson,''
            (2021), \href{https://arxiv.org/abs/2107.05718}arXiv:2107.05718.
\bibitem[Ba]{Barron}                  
K. D. Barron,
``N=1 Neveu-Schwarz vertex operator superalgebras over Grassmann algebras and with odd formal variables,'' 
 (1999), \href{https://arxiv.org/abs/math/9910007}{arXiv:math/9910007[math.QA]}.        
\bibitem[Be]{Belavin}
    V. A. Belavin,
    ``On the N= 1 super Liouville four-point functions,''
    \textit{Nuclear Physics B}, \textbf{798} (2008) 423-442. 
\bibitem[BPZ]{BPZ}    
A. A. Belavin, A. M. Polyakov and A. B. Zamolodchikov,
``Infinite conformal symmetry in two-dimensional quantum field theory,''
 \textit{Nuclear Physics B}, \textbf{241} (1984) 333-380. 
\bibitem[BS]{BA}
     L. Benoit and Y. Saint-Aubin,
     ``Fusion and the Neveu-Schwarz singular vectors,''
     \textit{International Journal of Modern Physics A} \textbf{9} (1994) 547.     
\bibitem[BMRW]{W} 
          O. Blondeau-Fournier, P. Mathieu, D. Ridout and S. Wood,
          ``Superconformal minimal models and admissible Jack polynomials,''
          \textit{Advances in Mathematics} \textbf{314} (2017) 71-123.
\bibitem[CRR]{CRR}
     M. Canagasabey, J. Rasmussen and D. Ridout,
     ``Fusion rules for the $N=1$ superconformal logarithmic minimal models I: The Neveu-Schwarz sector,''
     \textit{Journal of Physics A: Mathematical and Theoretical}, \textbf{48} (2015) 415402.
\bibitem[CGNS]{CGNS}
    T. Creutzig, N. Genra, S. Nakatsuka and R. Sato, 
    ``Correspondences of categories for subregular W-algebras and principal W-superalgebras,''
    \textit{Communications in Mathematical Physics}, \textbf{393} (2022) 1-60.  
\bibitem[CKL]{CKL}
     T. Creutzig, S. Kanade and A. R. Linshaw, 
     ``Simple current extensions beyond semi-simplicity,''   
     \textit{Commun. Contemp. Math}. \textbf{22} (2020), no. 1, 1950001, 49 pp. 
\bibitem[CKM]{CKM}
     T. Creutzig, S. Kanade and R. McRae, 
     ``Tensor categories for vertex operator superalgebra extensions,''
     (2017), \href{https://arxiv.org/abs/1705.05017}{arXiv:1705.05017[math.QA]}.    
\bibitem[CLR]{CLR}
     T. Creutzig, S. Lentner and M. Rupert, 
     ``An algebraic theory for logarithmic Kazhdan-Lusztig correspondences,''
     (2023), \href{https://arxiv.org/abs/2306.11492}{arXiv:2306.11492[math.QA]}.
\bibitem[CMY1]{CMY}
     T. Creutzig, R. McRae and J. Yang, 
     ``On ribbon categories for singlet vertex algebras,'' 
     \textit{Communications in Mathematical Physics}, \textbf{387} (2021) 865-925.
\bibitem[CMY2]{CMY2}
     T. Creutzig, R. McRae and J. Yang, 
     ``Tensor structure on the Kazhdan–Lusztig category for affine $\mathfrak{gl}(1|1)$,''
     \textit{International Mathematics Research Notices}, \textbf{16} (2022) 12462-12515.  
\bibitem[CMOY]{CMOY}
     T. Creutzig, R. McRae, F. Orosz Hunziker and J. Yang, 
     ``$N=1$ super Virasoro tensor categories,''
     (2024), \href{https://arxiv.org/abs/2412.18127}{arXiv:2412.18127[math.QA]}.  
\bibitem[CR]{CR}
     T. Creutzig and D. Ridout,
     ``Logarithmic conformal field theory: Beyond an introduction,''
     \textit{J. Phys.} \textbf{A46} (2013) 494006.
\bibitem[DF1]{DF1}
     V. S. Dotsenko and V. A. Fateev, 
     ``Conformal algebra and multipoint correlation functions in 2D statistical models,''
     \textit{Nucl. Phys.} \textbf{B240} (1984) 312-348. 
\bibitem[DF2]{DF2}
     V. S. Dotsenko and V. A. Fateev, 
     ``Four-point correlation functions and the operator algebra in 2D conformal invariant theories with central charge $c\leq 1$,''
     \textit{Nucl. Phys.} \textbf{B251} (1985) 691-734.
\bibitem[ESNO]{Etingof}
    P. Etingof, G. Shlomo, D. Nikshych and V. Ostrik, 
     \textit{Tensor Categories}. 
      Number volume 205 in Mathematical Surveys and Monographs. American Mathematical Society 2015.
\bibitem[FGST1]{FGST1}
     B. L. Feigin, A. M. Gainutdinov, A. M. Semikhatov and I. Y. Tipunin,
     ``The Kazhudan-Lusztig correspondence for the representation category of the triplet $W$-algebra in logarithmic conformal field theories,''
     \textit{Theoret. and Math. Phys.} $\bold{148}$ (2006) 1210-1235.
\bibitem[FGST2]{FF2}
     B. L. Feigin, A. M. Gainutdinov, A. M. Semikhatov and I. Y. Tipunin,
     ``Modular group representations and fusion in logarithmic conformal field theories and in the quantum group center,''
     \textit{Comm. Math. Phys.} \textbf{265} (2006) 47–93.
\bibitem[FGST3]{FF3}
     B. L. Feigin, A. M. Gainutdinov, A. M. Semikhatov and I. Y. Tipunin,
     ``Kazhdan-Lusztig correspondence for the representation category of the triplet W-algebra in logarithmic CFT,''
     \textit{Theor.Math.Phys.} \textbf{148} (2006) 1210-1235; \textit{Teor.Mat.Fiz.} \textbf{148}  (2006) 398-427.
\bibitem[FGST4]{FFL}
     B. L. Feigin, A. M. Gainutdinov, A. M. Semikhatov and I. Y. Tipunin,
     ``Logarithmic extensions of minimal models: characters and modular transformation,''
     \textit{Nucl. Phys.} \textbf{B757} (2006) 303-343. 
\bibitem[Fe]{Felder}
     G. Felder,
     ``BRST approach to minimal models,''
     \textit{Nuclear Phys.} \textbf{B317} (1989) 215-236.
\bibitem[FS]{FG}
     G. Felder and R. Silvotti, 
     ``Conformal blocks of minimal models on a Riemann surface,'' 
     \textit{Commun. Math. Phys.} \textbf{144} (1992) 17-40.  
\bibitem[FFHST]{FJ}
      J. Fjelstad, J. Fuchs, S. Hwang, A. M. Semikhatov and I. Y. Tipunin,
     ``Logarithmic conformal field theories via logarithmic deformation,''
     \textit{Nucl. Phys.} \textbf{B633} (2002) 379-413.
\bibitem[Fo1]{Forrester0}
       P. J. Forrester, 
       ``Integration formulas and exact calculations in the Calogero–Sutherland model,''
        \textit{Modern Physics Letters B}, \textbf{9} (1995) 359-371.
\bibitem[Fo2]{Forrester}
        P. J. Forrester, 
        \textit{Log-gases and random matrices (LMS-34)}. 
        Princeton University Press 2010.
\bibitem[FLM]{FLM}
I. Frenkel, J. Lepowsky and A. Meurman,
\textit{Vertex operator algebras and the Monster}.
Academic press 1989. 
\bibitem[GK1]{GK1}
M. R. Gaberdiel and H. G. Kausch,
        ``Indecomposable fusion products,''
        \textit{Nucl. Phys.} \textbf{B477} (1996) 293-318.
\bibitem[GK2]{GK}
M. R. Gaberdiel and H. G. Kausch,
        ``A Local Logarithmic Conformal Field Theory,''
        \textit{Nucl. Phys.} \textbf{B538} (1999) 631-658.
\bibitem[GN]{GN}
T. Gannon and C. Negron, 
``Quantum $ {\rm SL} (2)$ and logarithmic vertex operator algebras at $(p,1)$-central charge,''
(2021), 
\href{https://arxiv.org/abs/2104.12821}{arXiv:2104.12821[math.QA]}.     
\bibitem[Gu]{Gu}
V. Gurarie,
      ``Logarithmic Operators in Conformal Field Theory,''
      \textit{Nucl. Phys.} \textbf{B410} (1993) 535-549.
\bibitem[Ha]{Haraoka}
     Y. Haraoka,
     \textit{Linear Differential Equations in the Complex Domain}.
     Lecture Notes in Mathematics Vol. 2271, Springer, New York 2020.     
\bibitem[Hu]{H}
     Y. Z. Huang,
     ``Cofiniteness conditions, projective covers and the logarithmic tensor product theory,''
     \textit{J. Pure Appl. Algebra}, \textbf{213} (2009) 458-475.
\bibitem[HLZ1]{HLZ1}
     Y. Z. Huang,  J. Lepowsky and L. Zhang,
     ``Logarithmic tensor category theory for generalized modules for a conformal vertex algebra, I: Introduction and strongly graded algebras and their generalized modules,''
     \textit{Proceedings of a Workshop Held at Beijing International Center for Mathematical Research. Springer Berlin Heidelberg}, (2014) 169-248.
\bibitem[HLZ2]{HLZ2}
     Y. Z. Huang,  J. Lepowsky and L. Zhang,
     ``Logarithmic tensor category theory for generalized modules for a conformal vertex algebra, I\hspace{-1pt}I: Logarithmic formal calculus and properties of logarithmic intertwining operators,'' (2010), \href{https://arxiv.org/abs/1012.4196}{arXiv:1012.4196[math.QA]}.
\bibitem[HLZ3]{HLZ3}
     Y. Z. Huang,  J. Lepowsky and L. Zhang,
     ``Logarithmic tensor category theory
     for generalized modules for a conformal vertex algebra, III: Intertwining maps
     and tensor product bifunctors,''
      (2010), \href{https://arxiv.org/abs/1012.4197}{arXiv:1012.4197[math.QA]}.
\bibitem[HLZ4]{HLZ4}
     Y. Z. Huang,  J. Lepowsky and L. Zhang,
     ``Logarithmic tensor category theory for
generalized modules for a conformal vertex algebra, IV: Constructions of tensor
product bifunctors and the compatibility conditions,'' (2010),
\href{https://arxiv.org/abs/1012.4198}{arXiv:1012.4198[math.QA]}. 
\bibitem[HLZ5]{HLZ5}
     Y. Z. Huang,  J. Lepowsky and L. Zhang,
     ``Logarithmic tensor category theory
for generalized modules for a conformal vertex algebra, V: Convergence condition for intertwining maps and the corresponding compatibility condition,'' (2010),
\href{https://arxiv.org/abs/1012.4199}{arXiv:1012.4199[math.QA]}.
\bibitem[HLZ6]{HLZ6}
    Y. Z. Huang,  J. Lepowsky and L. Zhang,
     ``Logarithmic tensor category theory
for generalized modules for a conformal vertex algebra, VI: Expansion condition, associativity of logarithmic intertwining operators, and the associativity isomorphisms,'' (2010), \href{https://arxiv.org/abs/1012.4202}{arXiv:1012.4202[math.QA]}.
\bibitem[HLZ7]{HLZ7}
     Y. Z. Huang,  J. Lepowsky and L. Zhang,
     ``Logarithmic tensor category theory for generalized modules for a conformal vertex algebra, VII: Convergence
and extension properties and applications to expansion for intertwining maps,'' (2011),
\href{https://arxiv.org/abs/1110.1929}{arXiv:1110.1929[math.QA]}.
\bibitem[HLZ8]{HLZ8}
     Y. Z. Huang,  J. Lepowsky and L. Zhang,
      ``Logarithmic tensor category theory
for generalized modules for a conformal vertex algebra, VIII: Braided tensor
category structure on categories of generalized modules for a conformal vertex
algebra,'' (2011), \href{https://arxiv.org/abs/1110.1931}{arXiv:1110.1931[math.QA]}.
\bibitem[HM]{HM}
Y. Z. Huang and A. Milas, 
``Intertwining operator superalgebras and vertex tensor categories for superconformal algebras, I,''
\textit{Communications in Contemporary Mathematics}, \textbf{4} (2002) 327-355.
\bibitem[IK1]{IK1}
     K. Iohara and Y. Koga, 
     ``Representation theory of Neveu-Schwarz and Ramond algebras I: Verma modules,''
     \textit{Adv. Math}., \textbf{178} (2003) 1-65.
\bibitem[IK2]{IK2}
     K. Iohara and Y. Koga, 
     ``Representation theory of Neveu-Schwarz and Ramond Algebra I\hspace{-1pt}I: Fock modules,''
     \textit{Ann. Inst. Fourier, Grenoble}, 53, 6 (2003) 1755-1818.
\bibitem[JPR]{JPR}
     M. Jeng, G. Piroux and P. Ruelle,
     ``Height variables in the Abelian sandpile model: scaling fields and correlations,''
     \textit{J. Stat. Mech}, \textbf{0610} (2006) P105. 
\bibitem[JS]{JS}
     A. Joyal and R. Street,
     ``Braided tensor categories,''
     \textit{Adv. Math}. \textbf{102} (1993) 20-78.   
\bibitem[KW]{KW}
      V. Kac and W. Wang, 
      ``Vertex Operator Superalgebras and Their,''
      \textit{Mathematical aspects of conformal and topological field theories and quantum groups}, 175, 161 (1994).       
\bibitem[KR]{Kanade}
     S. Kanade and D. Ridout, 
     ``NGK and HLZ: Fusion for physicists and mathematicians,''
\textit{Affine, Vertex and W-algebras} (2019) 135-181.
\bibitem[Ka]{Ka}
     H. G. Kausch,
     ``Extended conformal algebras generated by multiplet of primary fields,''
     \textit{Phys. Lett. B}, \textbf{259} (1991) 448-455
\bibitem[KL]{KL}
     D. Kazhdan and G. Lusztig, 
     ``Tensor structures arising from affine Lie algebras. IV,''
     \textit{Journal of the American Mathematical Society}, \textbf{7} (1994) 383-453.
\bibitem[K\"{u}]{K}
      J. K\"{u}lshammer,
      ``Representation type and Auslander-Reiten theory of Frobenius-Lusztig kernels,''
      \textit{Ph.D. thesis, Christian-Albrechts-Universit\"{a}t zu Kiel} (2012).  
\bibitem[Li]{Li}
       H. Li,
       ``The physics superselection principle in vertex operator algebra theory,''
       \textit{Journal of Algebra} \textbf{196} (1997), no. 2, 436-457.       
\bibitem[MR]{MR}
       P. Mathieu and D. Ridout,
       ``From Percolation to Logarithmic Conformal Field Theory,''
       \textit{Phys. Lett.} \textbf{B657} (2007) 120-129.       
\bibitem[MY]{McRae}
      R. McRae and J. Yang,
      ``Structure of Virasoro tensor categories at central charge $13-6p-6p^{-1}$ for integers $p>1$,'' (2021),
      \href{https://arxiv.org/abs/2011.02170}{arXiv:2011.02170[math.QA]}.
\bibitem[NT]{NT}
      K. Nagatomo and A. Tsuchiya, 
      ``The Triplet Vertex Operator Algebra $W(p)$ and Restricted Quantum Group at Root of Unity,''
      \textit{Exploring new structures and natural constructions in mathematical physics. Vol. 61. Mathematical Society of Japan}, (2011) 1-50.
\bibitem[Nah]{Na}
      W. Nahm, 
      ``Quasirational fusion products,''
      \textit{Int. J. Mod. Phys.} \textbf{B8} (1994) 3693-3702. 
\bibitem[Nak]{Nak}   
        H. Nakano, ``The category of modules of the triplet W-algebras associated to the Virasoro minimal models,'' 
        \textit{the Doctor Thesis, Mathematical Institute, Tohoku University} (2023).          
\bibitem[Ni]{Ni}
      A. Nigro, 
      ``Integrals of Motion for Critical Dense Polymers and Symplectic Fermions,''
      \textit{J. Stat. Mech.} \textbf{0910} (2009) P10007.      
\bibitem[PR]{PR}
      P. A. Pearce and J. Rasmussen,
     ``Solvable critical dense polymers,''
     \textit{J. Stat. Mech}. \textbf{0702} (2007) P015.
\bibitem[RS]{RS}
    N. Read and H. Saleur,
    ``Associative-algebraic approach to logarithmic conformal field theories,''
     \textit{Nucl. Phys}. \textbf{B777} (2007) 316-351.      
\bibitem[Ri]{R}
    D. Ridout,
    ``On the Percolation BCFT and the Crossing Probability of Watts,''
    \textit{Nucl. Phys}, \textbf{B810} (2009) 503-526.       
\bibitem[RW]{RW}
    D. Ridout and S. Wood,
    ``Modular transformations and Verlinde formulae for logarithmic $(p_+, p_-)$-models,''
    \textit{Nucl. Phys}. \textbf{B880} (2014) 175-202.       
\bibitem[Su1]{S}
        E. Sussman,
       ``The singularities of Selberg--and Dotsenko--Fateev-like integrals,'' (2023), 
       \href{https://arxiv.org/abs/2301.03750}{arXiv:2301.03750[math.QA]}.
\bibitem[Su2]{S2}
      E. Sussman,
      ``The regularization of Dotsenko–Fateev integrals,'' 
      \textit{Letters in Mathematical Physics}, \textbf{113} (2023) 1-29.
\bibitem[Sut]{Suter}    
      R. Suter,
      ``Modules over $\mathfrak{U}_{q}(\mathfrak{sl}_2)$,'' 
      \textit{Communications in Mathematical Physics}, \textbf{163} (1994) 359-393.  
\bibitem[TK]{TK}
     A. Tsuchiya and Y. Kanie, 
     ``Fock space representations of the Virasoro algebra - Intertwining operators,''
     \textit{Publ. RIMS, Kyoto Univ}. \textbf{22} (1986) 259-327.
\bibitem[TW1]{TWFusion}
      A. Tsuchiya and S. Wood,
      ``The tensor structure on the representation category of the $\mathcal{W}_{p}$ triplet algebra,''
      \textit{J. Phys. A} \textbf{46} (2013) 445203.
\bibitem[TW2]{TW}
	A. Tsuchiya and S. Wood,
	``On the extended W-algebra of type $sl_2$ at positive rational level,''
	\textit{International Mathematics Research Notices} 2015.14 (2015) 5357-5435.
\bibitem[Xi]{Xia}
     J. Xiao,
     ``Restricted representations of $U(sl(2))$-quantizations,''
     \textit{Algebra Colloquium}, \textbf{1} (1994) 56-66. 
\bibitem[Xu]{Xu}
    X. P. Xu,
    ``Intertwining operators for twisted modules of a colored vertex operator superalgebra,'' 
    \textit{Journal of Algebra}, \textbf{175} (1995) 241-273.    
\bibitem[Zh]{Zhu}
     Y. Zhu,
     ``Modular invariance of characters of of vertex operator algebras,''
     \textit{J. Amer. Math. Soc}. $\bold{9}$ (1996) 237-302.


\end{thebibliography}
\end{document}